\documentclass[review]{elsarticle}

\usepackage{hyperref}

\journal{Finite Elements in Analysis and Design}

\usepackage{graphicx}
\usepackage{amsmath}
\usepackage{amsfonts}
\usepackage{amssymb}
\usepackage{bm}
\usepackage{placeins} 									
\usepackage{tabularx}
\newcolumntype{Y}{>{\centering\arraybackslash}X}
\def\tabulkywidth{125mm}

\usepackage[percent]{overpic}
\usepackage{arydshln}
\usepackage{verbatim}
\usepackage{float}
\usepackage{multirow}  
\usepackage[bottom]{footmisc} 
\usepackage{pgfplots} 
\pgfplotsset{compat=1.5}

\def\PART#1#2{\frac{\partial #1}{\partial #2}}
\def\ODT#1#2{\frac{\textrm{d} #1}{\textrm{d} #2}}
\def\T{\text{T}}
\def\n{\text{n}}
\def\w{\text{w}}
\def\s{\text{s}}
\def\b{\text{b}}
\def\m{\text{m}}
\def\e{\text{e}}
\def\dA{\text{d}A}
\def\rf{\text{ref}}
\def\US{US-ATFHS8~}

\definecolor{col1}{RGB}{0,150,0}

\makeatletter
\newcommand*{\Biggg}[1]{{\hbox{$\left#1\vbox to28\p@{}\right.\n@space$}}}
\newcommand*{\Bigggg}[1]{{\hbox{$\left#1\vbox to42\p@{}\right.\n@space$}}}
\makeatother

\def\arraystretch{0.8} 

\begin{document}
	
	\begin{frontmatter}
		
		\title{DKMQ24 shell element with improved membrane behaviour}
		
		
		\author[TUaddress]{V\'{\i}t\v{e}zslav \v{S}tembera$\!$\corref{cor1}}
		\ead{vitastembera@hotmail.com}
		\author[TUaddress]{Josef F\"{u}ssl}
		\address[TUaddress]{Institute for Mechanics of Materials and Structures, Vienna University of Technology, Karlsplatz 13, 1040 Vienna, Austria}
		\cortext[cor1]{Corresponding author}

		

		


		\begin{abstract}
			An approach to improve the membrane behaviour of the four-node shell element with 24 degrees of freedom DKMQ24 proposed by Katili et al. (2015) is presented. This improvement is based on a different approximation of drilling rotations, based on Allman's shape functions. Further, the element formulation is enhanced by the use of selective reduced integration of shear terms and by the proportional scaling of the penalty constant to the shell thickness. Additionaly, nodal moment corrections for the approximation of distributed normal loads are proposed. This element, called DKMQ24$_2$+, was tested on nine standard benchmark problems for the case of linear elasticity. The element is free of shear and membrane locking and shows lower numerical errors for membrane and membrane-dominated problems.
		\end{abstract} 
		
		
		
	\end{frontmatter}
	

	\section*{Nomenclature}
	
	\subsubsection*{Abbreviations}
	\begin{tabular}{ll}
		c.s.           & coordinate system \\
		dofs           & degrees of freedom \\
		Eq.            & equation \\
	\end{tabular} 
	
	\subsubsection*{Latin symbols}
	\begin{tabular}{ll}
		$a_i(r,s)$                                    & $i$-th basis function at point $(r,s)$ \\
		$c$                                           & penalty constant for the DKMQ24 element\\
		$c_1$                                         & penalty constant for the DKMQ24$_1$+ element\\
		$c_2$                                         & penalty constant for the DKMQ24$_2$+ element\\
		$E$                                           & Young's modulus \\
		$g$                                           & gravitational acceleration \\
		$G$                                           & shear modulus \\	
		$h(r,s)$                                      & shell thickness\\
		$\bm{n}_i$                                    & normal vector at node $i$ \\
		$\mathfrak{q}_x,\mathfrak{q}_y$               & shear internal forces \\
		$\bm{q}^{\T}=[\bm{\varphi}^{\T},~\bm{u}^{\T}]$& vector of unknown rotations and displacements \\
		$\bm{l}_5,\bm{l}_6,\bm{l}_7,\bm{l}_8$         & element side in-plane normal vectors \\
	\end{tabular} 
	\newpage
	\begin{tabular}{ll}
		$m_x, m_y, m_{xy}$                            & internal moments \\
		$n_x, n_y, n_{xy}$                            & membrane (in-plane) internal forces \\
		$\bm{n}_5,\bm{n}_6,\bm{n}_7,\bm{n}_8$         & element side normal vectors \\
		$r,s,t$                                       & coordinates in the element reference c.s.  \\
		$\bm{t}_5,\bm{t}_6,\bm{t}_7,\bm{t}_8$         & element side directional vectors\\	$\bm{T}(r,s)$                                 & transformation matrix from the global c.s. \\
		& to the local c.s. at point $(r,s)$ \\
		$u_x,u_y,u_z$                                 & displacements in the local c.s.  \\    
		$u_X,u_Y,u_Z$                                 & displacements in the global c.s.  \\
		$\bm{v}_1(r,s),\bm{v}_2(r,s),\bm{v}_3(r,s)\equiv\bm{n}(r,s)$   & local orthonormal basis at point $(r,s) $\\
		$x,y,z$                                       & coordinates in the local c.s.  \\
		$X,Y,Z$                                       & coordinates in the global c.s.  \\
		$\bm{X}_r(r,s),\bm{X}_s(r,s)$                 & covariant basis vectors at point $(r,s) $\\
		$\bm{X}^r(r,s),\bm{X}^s(r,s)$                 & dual basis vectors at point $(r,s) $\\
	\end{tabular} 
	
	\subsubsection*{Greek symbols}
	\begin{tabular}{ll}
		$\beta_x$                                         &  rotation of the $z$-axis towards the $x$-axis\\
		$\beta_y$                                         &  rotation of the $z$-axis towards the $y$-axis\\
		$\gamma_{xz}, \gamma_{yz}$                        &  transversal shear strains\\ 
		$\varepsilon_x, \varepsilon_y, \gamma_{xy}$       &  membrane (in-plane) strains\\ 
		$\kappa_x, \kappa_y, \kappa_{xy}$                 &  bending strains\\ 
		$\nu$											  &  Poisson's ratio\\	
		$\varphi_x,\varphi_y,\varphi_z$                   &  rotations in the local c.s.  \\
		$\varphi_X,\varphi_Y,\varphi_Z$                   &  rotations in the global c.s.  \\
		$\psi=\frac{1}{2}\left(\PART{u_y}{x}-\PART{u_x}{y}\right)$  &  in-plane rotation of the displacement field \\
	\end{tabular}
	

	\section{Introduction}
	Finite element modeling of shell structures is of great practical importance in engineering practice. Therefore there is a need for an efficient shell element, which would have an optimal convergence rate for different kinds of problems (membrane dominated, bending dominated and mixed shell problems), which would be at the same time free of shear locking, free of membrane locking, without zero energy modes and which would keep the optimal convergence rate also on distorted meshes. Especially, it is desired to find an optimal low-order shell element, which would provide the highest benefit-cost ratio when used in practical engineering applications. To find such an element is a tremendous task, not solved until today.

	Let us restrict ourselves in this paper to 4-nodal quadrangle shell elements, based on the linear approximation of the transversal shear (the Reissner-Mindlin concept), which use both translational and rotational degrees of freedom in a node. Moreover, let us concentrate on the case of small deformations, however, an extension to the case of large deformations, keeping the small strain assumption, can be easily made in the corotational approach \cite{Felippa, BattiniC, Tang}, in which no change of the element formulation is needed.
	
	In order to construct a quadrangle shell element, it is possible to either consider a flat element geometry, by which however it is necessary to use the warp correction term suggested by Taylor \cite{RigidLinkCorrectionTaylor, RigidLinkCorrectionVan} in the element formulation, or a warped element geometry. Formulations of flat quadrangles can be found in \cite{Groenwald, Tang, NguyenVan, Shang,  Wang}. Formulations using a non-flat geometry include the QUAD4 element suggested by MacNeal \cite{MacNeal} and the famous MITC4 (Mixed Interpolation Tensorial Components) element proposed by Dvorkin and Bathe \cite{MITC4}, which was further extended to the so called MITC4+ element \cite{MITC4Plus, MITC4Plus2}. Another extension of the MITC4 element called MITC4S was introduced by Niemi \cite{Niemi}. As alternative, shells can be modelled using solid-shell elements, as for example the US-ATFHS8 \cite{Huang}, which is a 8-nodal element with 24-dofs (displacement dofs only). Other recent shell formulations can be found in Gruttmannn \cite{Gruttmann}, Kim \cite{Kim} or Katili \cite{Katili2015}. The latter one introduces the DKMQ24 shell element, based on the Naghdi/Reissner/Mindlin shell theory with six degrees of freedom per node. This shell element is an extension of the DKMQ (Discrete Kirchhoff-Mindlin Quadrangle) plate element suggested by the same author, which proved to have good numerical properties on a wide set of problems \cite{DKMQ}. The DKMQ24 shell element has full coupling of membrane-bending energy, is free of shear locking, free of membrane locking and proved to converge on a wide set of standard benchmark problems \cite{Katili2015}. It was successfully applied to orthotropic multilayered shells \cite{KatiliComposite} and to a beam analysis \cite{DKMQ24Beams}. 
	
	The DKMQ24 shell element can be, due to its good numerical properties, considered as a promising candidate for an efficient low order quadrangle shell element. However, despite its numerical efficiency, it suffers from one drawback -- the drilling rotations are independent of the remaining degrees of freedom. The stiffness matrix of the DKMQ24 shell element has a full rank due to stabilization, however, the drilling rotations do not enrich the displacement interpolation space while increasing computational cost at the same time. One possibility to avoid this inefficiency would be the use of five degrees of freedom per node only, which was proposed by Irpanni \cite{Katili5dofs}. However, from a practical point of view six degrees of freedom at a node is of benefit in order to easily connect more arbitrary oriented shell (or beam) elements to a single node. 
	
	The aim of our contribution is to modify the formulation of the DKMQ24 element in order to overcome this imperfection. We propose the incorporation of drilling rotations with the help of Allman's shape functions \cite{Allman} in the element formulation. This enriches the in-place displacement field with incomplete quadratic functions and improves the membrane behaviour of the element. The element formulation contains a penalty term, forcing the drilling rotations to be equal to an in-plane rotation of the displacement field. This term is multiplied by the penalty coefficient, which is either constant or scaled to the shell thickness $h$. A further improvement of the element formulation is achieved by the selective reduced integration scheme of shear terms, as described in Section 3.4. Moreover, nodal moment corrections for the case of distributed normal loading are derived in Section 3.5. The modified element formulation is named either DKMQ24$_1$+ (for a constant penalty coefficient) or DKMQ24$_2$+ (for a penalty coefficient scaled to the shell thickness $h$). The formulation was tested on nine well-established benchmark problems, which show that it is free of shear locking, free of membrane locking, converges in all considered test cases, and provides a reduced numerical error in the membrane-dominated problems.
	
	The practical implementation of the DKMQ24, DKMQ24$_1$+ and DKMQ24$_2$+ shell elements was done in the C\# programming language using the Microsoft .NET Framework technology as a part of the finite element program femCalc \cite{femCalc}. 
	
	In Section 2 we give a full and detailed derivation of the DMKQ24 shell element introduced in \cite{DKMQ}. In Section 3 modifications of the DKMQ24 element, which construct the DKMQ24$_1$+ and DKMQ24$_2$+ elements, are introduced. The dependence of the relative error at the selected benchmarks on the penalty constant of the DKMQ24, DKMQ24$_1$+ and DKMQ24$_2$+ elements is discussed in Section 4. In Section 5, convergence behaviour of the DKMQ24, DKMQ24$_1$+ and DKMQ24$_2$+ elements on pure membrane, pure bending and on mixed membrane-bending benchmark problems is tested and discussed. Closing remarks are given in Section 6.
	\section{Formulation of the DKMQ24 shell element}
	\label{DKMQ24}
	In this Section a detailed derivation of the DKMQ24 shell element, introduced in \cite{Katili2015}, is given. The reason for this decision is that in reference \cite{Katili2015} some parts of the derivation are not presented in detail.  
	
	
	\subsection{Element characterization}
	The presented shell element is four-nodal with six unknowns at each node (three displacements $[u_X,u_Y,u_Z]$ and three rotations $[\varphi_X, \varphi_Y, \varphi_Z]$ in the global c.s.). The element is considered to be warped and considers full coupling between membrane and bending terms. Linear elasticity theory is assumed throughout the whole article.
	
	\subsection{Definition of rotations}
\begin{figure} 
		\centering
		\begin{overpic}[width=0.6\textwidth]{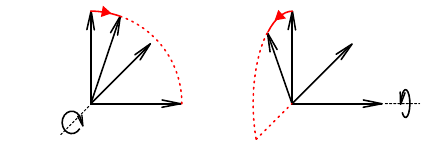} 
			\put(104,8.15){$\varphi_x$}		
			\put(41,4.25){$x$}
			\put(36,19.5){$y$}
			\put(20.25,32.5){$z$}
			\put(28.9,31){$z'$}
			\put(24,33.75){$\beta_x$}
			\put(9.75,-2){$\varphi_y$}
			\put(89,4.25){$x$}		
			\put(84.5,19.5){$y$}
			\put(68.75,32.75){$z$}
			\put(61.25,26.75){$z'$}
			\put(60.25,32.5){$-\beta_y$}
		\end{overpic}
		\caption{Definition of rotations}
		\label{betaphi2}
	\end{figure}
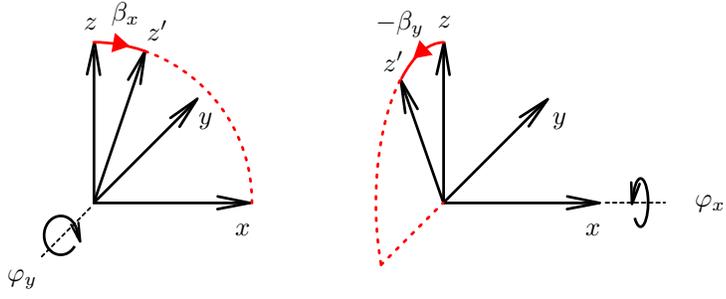

	Rotations $\varphi_x, \varphi_y, \varphi_z$ are defined using the standard definition, i.e. as rotations about axes $x,y,z$ using the right-hand-rule. Rotations $\varphi_x, \varphi_y$ can be equivalently defined as rotations of axis $z$ towards the axis $x$ or $y$ respectively, denoted by $\beta_x, \beta_y$. The relation between these definitions is the following (see Figure \ref{betaphi2}):
	\begin{align}
	\beta_x&=\varphi_y \\
	\beta_y&=-\varphi_x. 
	\end{align}
	We use rotations $\beta_x, \beta_y$ during the derivation of the element, because it is more convenient. However, the final formulation of the DKMQ24 shell element is given in the standard nodal rotations $\varphi_x, \varphi_y, \varphi_z$.
	
	\FloatBarrier
	
	\subsection{Coordinate system} 
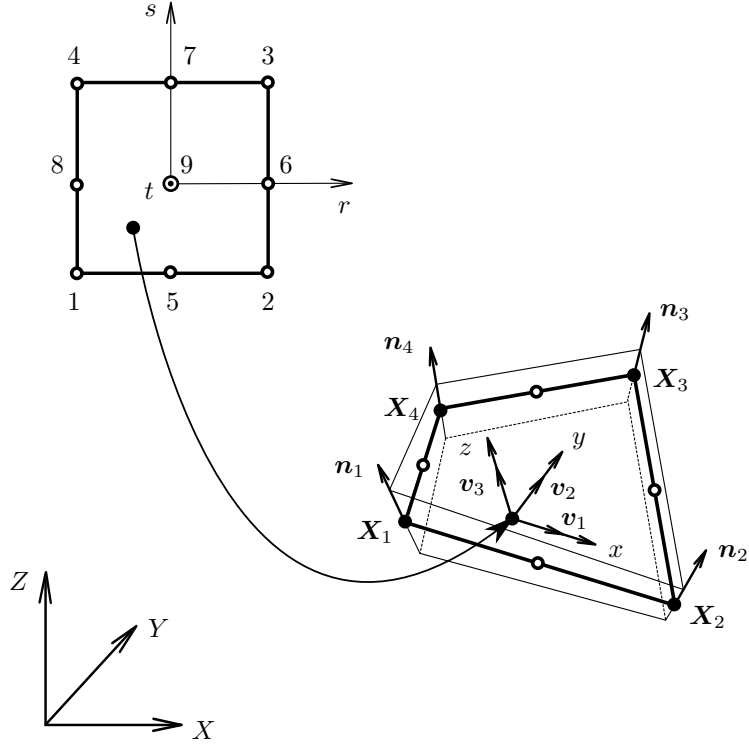
\begin{figure} 
		\centering
		\begin{overpic}[width=0.55\textwidth]{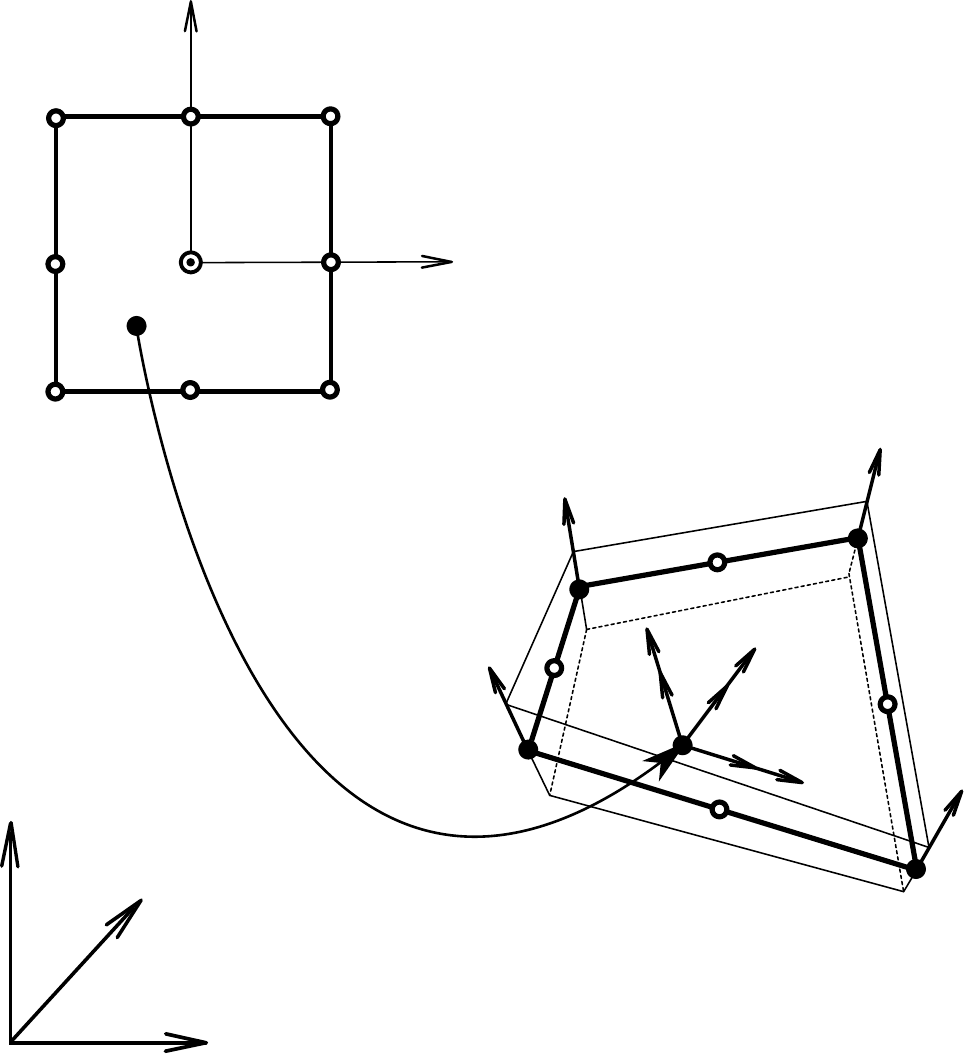} 
			\put(21,-1){$X$}
			\put(15,13){$Y$}
			\put(-4,19.5){$Z$}
			\put(4,57.75){$1$}
			\put(17.5,57.75){$5$}
			\put(30.5,57.75){$2$}
			\put(4, 91.5){$4$}
			\put(19.8,91.5){$7$}                         
			\put(30.5,91.5){$3$}
			\put(1.7,76.5){$8$}
			\put(19.4,76.5){$9$}
			\put(33,76.5){$6$}
			\put(41,71){$r$}		
			\put(14.5,98){$s$}		
			\put(14.5,73){$t$}		
			\put(78,24){$x$}		
			\put(73,39.5){$y$}		
			\put(57.5,38){$z$}		
			\put(71.5,28.3){$\bm{v}_1$}		
			\put(70,32.5){$\bm{v}_2$}		
			\put(57.5,33){$\bm{v}_3$}				
			\put(40.5,35.5){$\bm{n}_1$}		
			\put(93,24){$\bm{n}_2$}		
			\put(85,57){$\bm{n}_3$}		
			\put(47.25,52.5){$\bm{n}_4$}		
			\put(43.5,26.5){$\bm{X}_1$}		
			\put(89,14.75){$\bm{X}_2$}		
			\put(84,47.25){$\bm{X}_3$}		
			\put(47.2,43.25){$\bm{X}_4$}		
		\end{overpic}
		\caption{The reference and local configuration for the DKMQ24 shell element. The numbering follows the Gmsh standard \cite{gmsh}.} 
		\label{referencelocal}
	\end{figure}

	The geometry of an isoparametric curvilinear shell element in the reference configuration is described by 
	\begin{align}
	\bm{X}(r,s,t)=\sum_{i=1}^4 a_i(r,s) \left[\bm{X}_i+t\frac{h}{2}\bm{n}_i\right], \label{eq3}
	\end{align}
	where $\bm{X}_i=[X_i,Y_i,Z_i]$ is the location of the $i$-th node in the reference (undeformed) configuration, $h$ is a shell thickness and $\bm{n}_i$ is a given normal vector at node $i$, see Figure \ref{referencelocal}. The shell thickness $h$ in Eq. (\ref{eq3}) is assumed to be constant in derivation of the strain matrices (similar to the approach used in \cite{Lakshminarayana}). The actual thickness and its variation within the element is introduced into the constitutive relations only (Eq. (\ref{const1}, \ref{Dm}, \ref{const2})). 
	
	The basis functions are summarized in Table 1. At this point we can calculate a covariant vector basis with respect to all parameters 
	\begin{align}
	\bm{X}_{,r}(r,s,t)&=\sum_{i=1}^4 a_{i,r}(r,s) \left[\bm{X}_i+t\frac{h}{2}\bm{n}_i\right], \\
	\bm{X}_{,s}(r,s,t)&=\sum_{i=1}^4 a_{i,s}(r,s) \left[\bm{X}_i+t\frac{h}{2}\bm{n}_i\right],  \\
	\bm{X}_{,t}(r,s,t)&=\sum_{i=1}^4 a_{i}(r,s) \frac{h}{2} \bm{n}_i,
	\end{align}
	where we have used the notation $\bullet_{,r}:=\PART{\bullet}{r},~\bullet_{,s}:=\PART{\bullet}{s}$. These vectors are in general neither orthogonal nor orthonormal. Let us denote
	\begin{align}
	\bm{X}_{r}(r,s)&:=\bm{X}_{,r}(r,s,0)=\sum_{i=1}^N a_{i,r}(r,s)\bm{X}_i, \\
	\bm{X}_{s}(r,s)&:=\bm{X}_{,s}(r,s,0)=\sum_{i=1}^N a_{i,s}(r,s)\bm{X}_i. 
	\end{align}
	Let us define a normal vector at any point $(r,s)$ by
	\begin{align}
	\bm{v}_3(r,s)&=\frac{\bm{X}_{r}(r,s)\times\bm{X}_{s}(r,s)}{||\bm{X}_{r}(r,s)\times\bm{X}_{s}(r,s)||},
	\end{align}
	with the help of which we define normal vectors at nodes $\bm{n}_i=\bm{v}_3(r_i,s_i),~i=1,\dots,4$. These vectors are used in Eq. (\ref{eq3}). The remaining orthonormal vectors $\bm{v}_1(r,s)$ and $\bm{v}_2(r,s)$ can be obtained simply by
	\begin{align}
	\bm{v}_1(r,s) &= \frac{\bm{X}_{r}(r,s)}{||\bm{X}_{r}(r,s)||},\\
	\bm{v}_2(r,s) &= \frac{\bm{v}_3(r,s)\times\bm{v}_1(r,s)}{||\bm{v}_3(r,s)\times\bm{v}_1(r,s)||},
	\end{align}
	where the vectors $\bm{v}_1,~\bm{v}_2$ and $\bm{v}_3$ comprise the local coordinate system. If an orthotropic material of the shell is considered, the vector $\bm{v}_1$ is usually chosen to be collinear with the direction of fibers.
	
	\begin{center}
		\begin{table}
			\centering
			\begin{tabular}{c|c|c}
				node $i$  & $[r_i,s_i]$ & $a_i(r,s)$                 \rule[-2mm]{0mm}{6mm} \\ 
				\hline
				$1$       & $[ -1, -1]$ & $\frac{1}{4}(1-r)(1-s)$    \rule[-2mm]{0mm}{6mm} \\
				$2$       & $[~~1, -1]$ & $\frac{1}{4}(1+r)(1-s)$    \rule[-2mm]{0mm}{6mm} \\
				$3$       & $[~~1,~~1]$ & $\frac{1}{4}(1+r)(1+s)$    \rule[-2mm]{0mm}{6mm} \\
				$4$       & $[ -1,~~1]$ & $\frac{1}{4}(1-r)(1+s)$    \rule[-2mm]{0mm}{6mm} \\
			\end{tabular} \nonumber
			\caption{Basis functions}
		\end{table}
	\end{center}
	
	\FloatBarrier
	
	\subsection{Transformation between global and local coordinate systems} 
	The local orthonormal coordinate system $\bm{v}_1(r,s),\bm{v}_2(r,s),\bm{v}_3(r,s)$ forms the transformation matrix $\bm{T}(r,s)$
	\begin{align}
	\bm{T}(r,s)=\left[\begin{array}{c}
	\bm{v}^{\T}_1(r,s) \\ \bm{v}^{\T}_2(r,s) \\ \bm{v}^{\T}_3(r,s)
	\end{array}
	\right],~~\bm{T}^{-1}(r,s)=\bm{T}^{\T}(r,s)=[\bm{v}_1(r,s) ~ \bm{v}_2(r,s) ~ \bm{v}_3(r,s)].
	\end{align}
	The same transformation rules hold for coordinates, displacements and rotations
	\begin{align}
	\left[
	\begin{array}{c}
	x \\ y \\ z
	\end{array}
	\right]=\bm{T}
	\left[
	\begin{array}{c}
	X \\ Y \\ Z
	\end{array}
	\right],~ 
	\left[
	\begin{array}{c}
	u_x \\ u_y \\ u_z
	\end{array}
	\right]=\bm{T}
	\left[
	\begin{array}{c}
	u_X \\ u_Y \\ u_Z
	\end{array}
	\right],~ 
	\left[
	\begin{array}{c}
	\varphi_x \\ \varphi_y \\ \varphi_z
	\end{array}
	\right]=\bm{T}
	\left[
	\begin{array}{c}
	\varphi_X \\ \varphi_Y \\ \varphi_Z
	\end{array}
	\right], \label{eq122}
	\end{align}
	where upper-case letters refer to the reference coordinate system and lower-case letters refer to the local coordinate system. The matrix which comprises the orthonormal basis is orthogonal, therefore inversion of above relations is trivial:
	\begin{align}
	\left[
	\begin{array}{c}
	X \\ Y \\ Z
	\end{array}
	\right]=\bm{T}^{\T}
	\left[
	\begin{array}{c}
	x \\ y \\ z
	\end{array}
	\right],~ 
	\left[
	\begin{array}{c}
	u_X \\ u_Y \\ u_Z
	\end{array}
	\right]=\bm{T}^{\T}
	\left[
	\begin{array}{c}
	u_x \\ u_y \\ u_z
	\end{array}
	\right],~ 
	\left[
	\begin{array}{c}
	\varphi_X \\ \varphi_Y \\ \varphi_Z
	\end{array}
	\right]=\bm{T}^{\T}
	\left[
	\begin{array}{c}
	\varphi_x \\ \varphi_y \\ \varphi_z
	\end{array}
	\right].
	\end{align}
	
	\FloatBarrier
	\subsection{Displacements} 
	Let us denote the orthonormal basis vectors $\bm{v}_1(r,s),\bm{v}_2(r,s),\bm{v}_3(r,s)$ evaluated at $i$-th node simply by $\bm{v}_1^i,\bm{v}_2^i,\bm{v}_3^i$. The displacement field of the DKMQ24 shell element is then defined by 
	\begin{align}
	\left[
	\begin{array}{c}
	u_X \\ u_Y \\ u_Z
	\end{array}
	\right](r,s,t)=&\sum_{i=1}^4a_i(r,s)\left[
	\left[
	\begin{array}{c}
	u_{iX} \\ u_{iY} \\ u_{iZ}
	\end{array}
	\right]+t\frac{h}{2}\bm{E}_i
	\left[
	\begin{array}{c}
	\varphi_{iX} \\ \varphi_{iY} \\ \varphi_{iZ}
	\end{array}
	\right]
	\right]+\nonumber\\
	&t\frac{h}{2}\sum_{k=5}^8a_k(r,s)\triangle\beta_{tk}\bm{t}_k,~~t\in[-1,1], \label{eq15}
	\end{align}
	where
	\begin{align}
	\bm{E}_i=-\text{Spin}(\bm{v}^i_3)=\left[
	\begin{array}{ccc}
	~~0         &~~v^i_{3Z} & -v^i_{3Y}\\
	-v^i_{3Z}  &~~0        &~~v^i_{3X}\\
	~~v^i_{3Y}  & -v^i_{3X} &~~0       \\
	\end{array}
	\right].
	\end{align}
	
	\begin{center}
		\begin{table}
			\centering
			\begin{tabular}{c|c|c}
				node $k$  & $[r_k,s_k]$ & $a_k(r,s)$                 \rule[-2mm]{0mm}{6mm} \\ 
				\hline
				$5$       & $[~~0, -1]$ & $\frac{1}{2}(1-r^2)(1-s)$    \rule[-2mm]{0mm}{6mm} \\
				$6$       & $[~~1,~~0]$ & $\frac{1}{2}(1+r)(1-s^2)$    \rule[-2mm]{0mm}{6mm} \\
				$7$       & $[~~0,~~1]$ & $\frac{1}{2}(1-r^2)(1+s)$    \rule[-2mm]{0mm}{6mm} \\
				$8$       & $[ -1,~~0]$ & $\frac{1}{2}(1-r)(1-s^2)$    \rule[-2mm]{0mm}{6mm} \\
			\end{tabular} \nonumber
			\caption{The supplementary basis functions for midside rotations $\triangle \beta_{tk}$ on side $k=5,\dots,8$}		
		\end{table}
	\end{center}
	
	We have used the definition of the spin matrix $\text{Spin}(\bm{v})$ taken from \cite{Felippa}. The supplementary quadratic basis functions $a_5, a_6, a_7, a_8$, which enrich the bending behaviour of the element, are summarized in Table 2. 
	
\begin{figure} 
		\centering
		\begin{overpic}[width=0.80\textwidth]{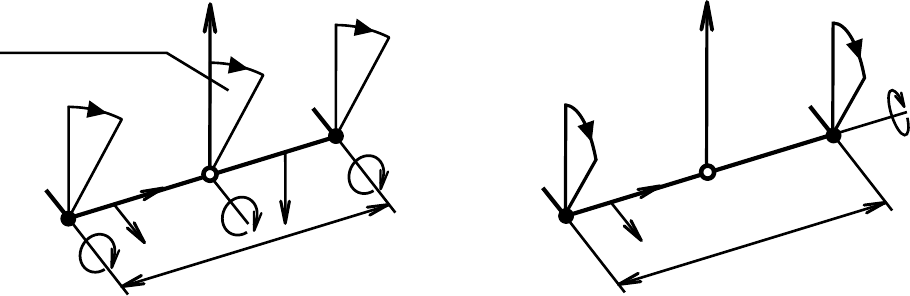} 
			\put(-0.5,28.5){$\frac{\beta_{ti}+\beta_{tj}}{2}+\triangle\beta_{tk}$}
			\put(8.5,23){$\beta_{ti}$}
			\put(38,32.25){$\beta_{tj}$}	
			\put(32.9,9.5){$\bm{q}$}		
			\put(24,32.75){$\bm{n}_k$}		
			\put(5,5){$i$}		
			\put(39.75,17){$j$}		
			\put(16.5,4){$m$}		
			\put(17,12.85){$t$}		
			\put(28.5,1.75){$L_k$}		
			\put(62,23){$\beta_{mi}$}
			\put(91.5,32){$\beta_{mj}$}	
			\put(79,32.75){$\bm{n}_k$}		
			\put(60,5){$i$}		
			\put(94.75,15.75){$j$}		
			\put(71.5,4){$m$}		
			\put(72,12.85){$t$}		
			\put(83.5,1.75){$L_k$}		
		\end{overpic}
		\caption{Approximation of rotations $\beta_t,\beta_m$ on the $k$-th element side}
		\label{betatm}
	\end{figure}

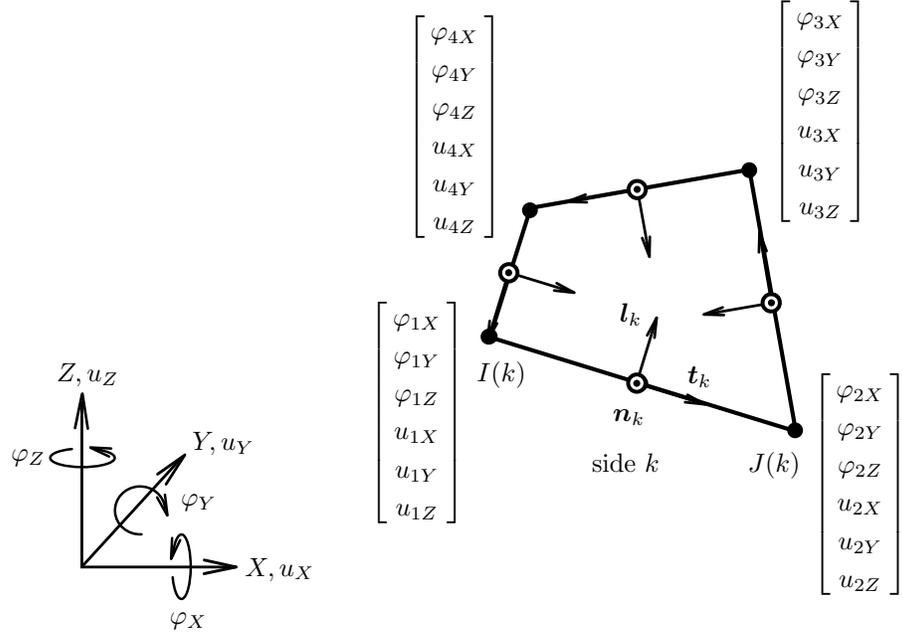
\begin{figure} 
		\centering
		\begin{overpic}[width=0.62\textwidth]{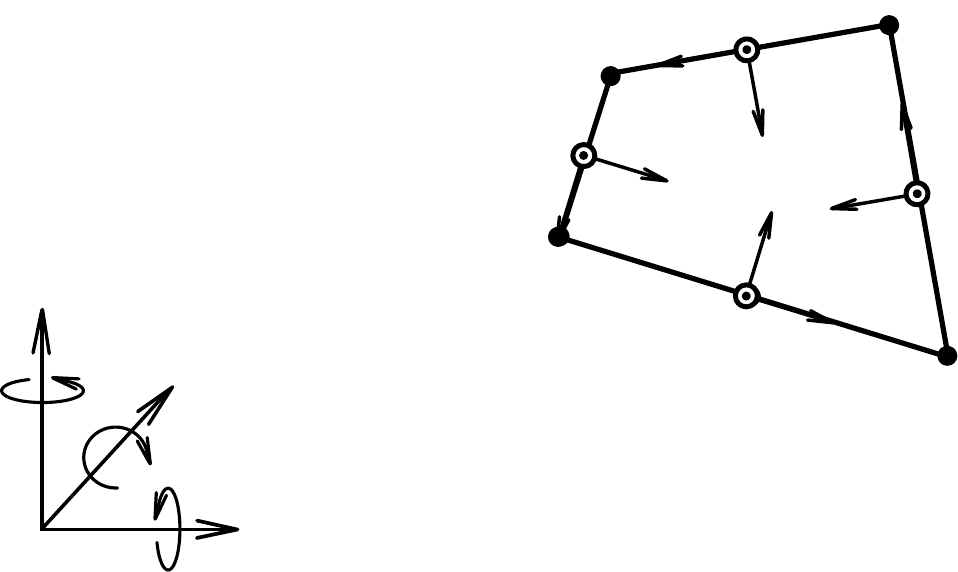} 
			\put(74.8,23.8){$\bm{n}_k$}		
			\put(72,17){side $k$}		
			\put(84.7,28.8){$\bm{t}_k$}		
			\put(75.7,37){$\bm{l}_k$}
			\put(56.77,29){$I(k)$}		
			\put(92.7,17){$J(k)$}						
			\put(26,3.25){$X,u_X$}		
			\put(19,20){$Y,u_Y$}		
			\put(1,29){$Z,u_Z$}		
			\put(16,-3.2){$\varphi_X$}		
			\put(17.5,12.5){$\varphi_Y$}		
			\put(-5,18.5){$\varphi_Z$}		
			\put(42.5,24){$\arraycolsep=2.2pt\def\arraystretch{0.8}\left[\begin{array}{c}
				\varphi_{1X}\\
				\varphi_{1Y}\\
				\varphi_{1Z}\\
				u_{1X}\\
				u_{1Y}\\
				u_{1Z}\\
				\end{array}\right]$}		
			\put(101.3,14.5){$\arraycolsep=2.2pt\def\arraystretch{0.8}\left[\begin{array}{c}
				\varphi_{2X}\\
				\varphi_{2Y}\\
				\varphi_{2Z}\\
				u_{2X}\\
				u_{2Y}\\
				u_{2Z}\\
				\end{array}\right]$}		
			\put(96,64){$\arraycolsep=2.2pt\def\arraystretch{0.8}\left[\begin{array}{c}
				\varphi_{3X}\\
				\varphi_{3Y}\\
				\varphi_{3Z}\\
				u_{3X}\\
				u_{3Y}\\
				u_{3Z}\\
				\end{array}\right]$}		
			\put(47.5,62){$\arraycolsep=2.2pt\def\arraystretch{0.8}\left[\begin{array}{c}
				\varphi_{4X}\\
				\varphi_{4Y}\\
				\varphi_{4Z}\\
				u_{4X}\\
				u_{4Y}\\
				u_{4Z}\\
				\end{array}\right]$}		
		\end{overpic}
		\caption{Nodal unknowns and tangential and normal vectors on side $k$. $I(k)$ the starting node and $J(k)$ is the ending node of side $k$. }
		\label{nodalunknowns}
	\end{figure}

	The second term in Eq. (\ref{eq15}) represents the displacement increment caused by nodal rotations. Let us derive this term for the $i$-th node: 
	
	\begin{align}
	\left[
	\begin{array}{c}
	\triangle u_{iX} \\
	\triangle u_{iY} \\
	\triangle u_{iZ} \\
	\end{array}
	\right]=\,\,&z\,
	[\bm{v}^i_1~~\bm{v}^i_2]
	\left[
	\begin{array}{c}
	\beta_{ix} \\ \beta_{iy}
	\end{array}
	\right]
	=z\,
	[-\bm{v}^i_2~~\bm{v}^i_1~~\bm{0}]
	\left[
	\begin{array}{c}
	\varphi_{ix} \\ \varphi_{iy} \\ \varphi_{iz}
	\end{array}
	\right]\stackrel{(\ref{eq122})}{=}\nonumber \\
	&z\!
	\left[
	\begin{array}{ccc}
	-v^i_{2X} & v^i_{1X} & 0 \\
	-v^i_{2Y} & v^i_{1Y} & 0 \\
	-v^i_{2Z} & v^i_{1Z} & 0 \\
	\end{array}
	\right]
	\underbrace{\left[
		\begin{array}{ccc}
		v^i_{1X} & v^i_{1Y} & v^i_{1Z} \\
		v^i_{2X} & v^i_{2Y} & v^i_{2Z} \\
		v^i_{3X} & v^i_{3Y} & v^i_{3Z} \\
		\end{array}
		\right]}_{\bm{T}(r_i, s_i)}
	\left[
	\begin{array}{c}
	\varphi_{iX} \\ \varphi_{iY} \\ \varphi_{iZ}
	\end{array}
	\right]= \nonumber \\
	& z\!
	\underbrace{\left[
		\begin{array}{ccc}
		0         & ~~v^i_{3Z} &  -v^i_{3Y} \\
		-v^i_{3Z} &   0         & ~~v^i_{3X} \\
		~~v^i_{3Y} &  -v^i_{3X} &   0         \\
		\end{array}
		\right]}_{\bm{E}_i} \left[
	\begin{array}{c}
	\varphi_{iX} \\ \varphi_{iY} \\ \varphi_{iZ}
	\end{array}
	\right],~~z=t\frac{h}{2},~~t\in[-1,1],~~i=1,2,3,4,
	\end{align}
	where we have used the definition of the vector multiplication $\bm{v}_3=\bm{v}_1\times\bm{v}_2$. 
	The last term in Eq. (\ref{eq15}) contains supplementary degrees of freedom $\triangle\beta_{tk}$ (see Figure \ref{betatm}). The supplementary degrees of freedom enrich the rotation approximation in the tangential $t$-direction, where the quadratic approximation is considered. In the perpendicular $m$-direction to side $k$ only a linear approximation of rotation $\beta_m$ is used:
	\begin{align}
	\beta_t(t)&=\left(1-\frac{t}{L_k}\right)\beta_{ti}+\frac{t}{L_k}\beta_{tj}+4\frac{t}{L_k}\left(1-\frac{t}{L_k}\right)\triangle\beta_{tk}, \\
	\beta_m(t)&=\left(1-\frac{t}{L_k}\right)\beta_{mi}+\frac{t}{L_k}\beta_{mj}. \label{beta}
	\end{align}
	The supplementary degrees of freedom $\triangle\beta_{tk}$ will be eliminated later on the element level, so that the only element unknowns remain three rotations and three displacements at each node (see Figure \ref{nodalunknowns}).
	
	\FloatBarrier
	
	\FloatBarrier 
	
	\subsection{Jacobi matrix}
	The derivative of the nodal location in the global coordinate system (\ref{eq3}) with respect to parameters $r,s,t$ yields
	\begin{align}
	\bm{J}=\PART{(X, Y, Z)}{(r,s,t)} &= \left[\bm{X}_r,\bm{X}_s,\bm{m}\right]+t\left[\bm{m}_r,\bm{m}_s,\bm{0}\right],
	\end{align}
	where
	\begin{align}
	\bm{X}_r &= \sum_{i=1}^4 a_{i,r}(r,s) \bm{X}_i, \nonumber \\
	\bm{X}_s &= \sum_{i=1}^4 a_{i,s}(r,s) \bm{X}_i, \nonumber \\
	\bm{m} &:= \frac{h}{2}\sum_{i=1}^4 a_{i}(r,s)\bm{n}_i, \nonumber \\
	\bm{m}_r &:= \frac{h}{2}\sum_{i=1}^4 a_{i,r}(r,s)\bm{n}_i, \nonumber \\
	\bm{m}_s &:= \frac{h}{2}\sum_{i=1}^4 a_{i,s}(r,s)\bm{n}_i. \nonumber  
	\end{align}
	The Jacobi matrix required for the calculation of the membrane stiffness term yields
	\begin{align}
	\overline{\bm{J}}:=\lim_{t \to 0} \PART{(r,s,t)}{(X, Y, Z)} &= \lim_{t \to 0} \bm{J}^{-1} = \left[\lim_{t \to 0} \bm{J}\right]^{-1} \!\!\!\!\!= \left[\bm{X}_r,~~\bm{X}_s,~~\bm{m}\right]^{-1} = \left[
	\begin{array}{c}
	(\bm{X}^r)^{\T} \\ (\bm{X}^s)^{\T} \\ \bm{m}^{\T} / ||\bm{m}||_2^2
	\end{array}
	\right].
	\end{align}
	For the calculation of the bending stiffness term the following term is also needed: $\lim_{t \to 0} \ODT{}{t} \PART{(r,s,t)}{(X, Y, Z)}$. By recalling the matrix identity 
	\begin{align}
	\ODT{\bm{A}^{-1}}{t}=-\bm{A}^{-1}\ODT{\bm{A}}{t}\bm{A}^{-1},
	\end{align}
	which follows from the derivation of the identity matrix $\bm{I}=\bm{A}\bm{A}^{-1}$, we get
	\begin{align}
	\overline{\bm{J}}'&:=\lim_{t \to 0} \ODT{}{t} \PART{(r,s,t)}{(X, Y, Z)} = \lim_{t \to 0} \ODT{\bm{J}^{-1}}{t}  = 
	\lim_{t \to 0} \left[-\bm{J}^{-1} \ODT{\bm{J}}{t} \bm{J}^{-1}\right] = \nonumber \\
	&-\overline{\bm{J}} \left[\lim_{t \to 0} \ODT{\bm{J}}{t}\right] \overline{\bm{J}} = 
	-\overline{\bm{J}} \left[\bm{m}_r,~~\bm{m}_s,~~\bm{0}\right] \overline{\bm{J}},
	\end{align}
	where upper indices denote dual vectors\footnote{For the set of linearly independent vectors $\bm{u}_1,\dots,\bm{u}_N$ one can uniquely define a set of dual vectors $\bm{u}^1,\dots,\bm{u}^N$, which satisfy $\bm{u}_i^{\T}\bm{u}^j=\delta_{ij},~i,j=1,\dots,N$, where $\delta_{ij}$ is the Kronecker delta symbol.}. Rows of the $\overline{\bm{J}}'$ matrix are explicitly denoted by vectors $\bm{o}_r^{\T}, \bm{o}_s^{\T}, \bm{o}_t^{\T}$:
	\begin{align}
	\overline{\bm{J}}'=
	\left[
	\begin{array}{c}
	\bm{o}_r^{\T}\\
	\bm{o}_s^{\T}\\
	\bm{o}_t^{\T}\\
	\end{array}
	\right].
	\end{align}
	
	\subsection{Displacement derivatives}
	The derivative of the displacement field (\ref{eq15}) with respect to parameters $r,s,t$ yields
	\begin{align}
	\PART{(u_X, u_Y, u_Z)}{(r,s,t)}&=\left[\bm{u}_r,~\bm{u}_s,~\bm{\Phi}+\bm{\beta}\right]+t\left[\bm{\Phi}_r+\bm{\beta}_r,~\bm{\Phi}_s+\bm{\beta}_s,~\bm{0}\right],
	\end{align}
	where
	\begin{align}
	\bm{u}_r&:=\sum_{i=1}^4a_{i,r}(r,s)\left[
	\begin{array}{c}
	u_{iX} \\ u_{iY} \\ u_{iZ}
	\end{array}
	\right],\nonumber \\
	\bm{u}_s&:=\sum_{i=1}^4a_{i,s}(r,s)\left[
	\begin{array}{c}
	u_{iX} \\ u_{iY} \\ u_{iZ}
	\end{array}
	\right],\nonumber \\
	\bm{\Phi}&:=\frac{h}{2}\sum_{i=1}^4a_{i}(r,s)\bm{E}_i\left[
	\begin{array}{c}
	\varphi_{iX} \\ \varphi_{iY} \\ \varphi_{iZ}
	\end{array}
	\right],\nonumber \\
	\bm{\Phi}_r&:=\frac{h}{2}\sum_{i=1}^4a_{i,r}(r,s)\bm{E}_i\left[
	\begin{array}{c}
	\varphi_{iX} \\ \varphi_{iY} \\ \varphi_{iZ}
	\end{array}
	\right],\nonumber \\
	\bm{\Phi}_s&:=\frac{h}{2}\sum_{i=1}^4a_{i,s}(r,s)\bm{E}_i\left[
	\begin{array}{c}
	\varphi_{iX} \\ \varphi_{iY} \\ \varphi_{iZ}
	\end{array}
	\right],\nonumber \\
	\bm{\beta}  &:=\frac{h}{2}\sum_{k=5}^8a_{k}(r,s)\triangle\beta_{tk}\bm{t}_k, \nonumber \\
	\bm{\beta}_r&:=\frac{h}{2}\sum_{k=5}^8a_{k,r}(r,s)\triangle\beta_{tk}\bm{t}_k, \nonumber \\
	\bm{\beta}_s&:=\frac{h}{2}\sum_{k=5}^8a_{k,s}(r,s)\triangle\beta_{tk}\bm{t}_k. \nonumber 
	\end{align}
	
	\subsection{Membrane and bending strain matrices}
	The vector of element unknowns $\bm{q}$ of size $24\times1$ is composed in the following way
	\begin{align}
	\bm{q}=\left[
	\begin{array}{c}
	\bm{\varphi}\\
	\bm{u}\\
	\end{array}
	\right]=\def\arraystretch{0.6}\left[
	\begin{array}{c}
	\varphi_{1X} \\
	\vdots \\
	\varphi_{1Y} \\
	\vdots \\
	\varphi_{1Z} \\
	\vdots \\
	u_{1X} \\
	\vdots \\
	u_{1Y} \\
	\vdots \\
	u_{1Z} \\
	\vdots \\
	\end{array}
	\right]. \label{qvector}
	\end{align}
	The membrane and bending strain matrices can be calculated in the local coordinate system using the strain definition 
	\begin{align}
	\bm{\varepsilon}&=\left[
	\begin{array}{c}
	\varepsilon_x \\ \varepsilon_y \\ \gamma_{xy}
	\end{array}
	\right]=\lim_{t \to 0}
	\left[
	\begin{array}{c}
	\PART{u_x}{x} \\
	\PART{u_y}{y} \\
	\PART{u_x}{y} + \PART{u_y}{x}\\
	\end{array}
	\right],\nonumber \\
	\bm{\kappa}&=\left[\begin{array}{c}
	\kappa_x \\ \kappa_y \\ \kappa_{xy}
	\end{array}
	\right]=\lim_{t \to 0}\left[
	\ODT{}{z}
	\left[
	\begin{array}{c}
	\PART{u_x}{x} \\
	\PART{u_y}{y} \\
	\PART{u_x}{y} + \PART{u_y}{x}\\
	\end{array}
	\right]\right]=\frac{2}{h}\lim_{t \to 0}\left[
	\ODT{}{t}
	\left[
	\begin{array}{c}
	\PART{u_x}{x} \\
	\PART{u_y}{y} \\
	\PART{u_x}{y} + \PART{u_y}{x}\\
	\end{array}
	\right]\right],\nonumber \\
	&z=\frac{h}{2}t. 
	\end{align}
	Using the chain rule we get
	\begin{align}
	\PART{u_x}{x} &= \underbrace{\PART{u_x}{(u_X, u_Y, u_Z)}}_{\bm{v}_1^{\T}(r,s)}\PART{(u_X, u_Y, u_Z)}{(r,s,t)}\PART{(r,s,t)}{(X,Y,Z)}\underbrace{\PART{(X,Y,Z)}{x}}_{\bm{v}_1(r,s)}, \\
	\PART{u_x}{y} &= \underbrace{\PART{u_x}{(u_X, u_Y, u_Z)}}_{\bm{v}_1^{\T}(r,s)}\PART{(u_X, u_Y, u_Z)}{(r,s,t)}\PART{(r,s,t)}{(X,Y,Z)}\underbrace{\PART{(X,Y,Z)}{y}}_{\bm{v}_2(r,s)},
	\end{align}
	and analogously for the remaining derivatives $\PART{u_y}{x}$ and $\PART{u_y}{y}$. The membrane strain $\bm{\varepsilon}$ yields
	\begin{align}
	\bm{\varepsilon}=\left[
	\begin{array}{c}
	\varepsilon_x \\ \varepsilon_y \\ \gamma_{xy}
	\end{array}
	\right]&=\bm{B}_{\m,3\times24}\, \bm{q},  \label{eps0}
	\end{align}
	where
	\begin{align}
	\bm{B}_{\m,3\times24}&=
	\left[
	\begin{array}{c}
	E_{11} \\ E_{22} \\ E_{12} + E_{21}
	\end{array}
	\right], \label{Bm} \\
	E_{IJ}&=\left[\bm{v}_I^{\T}\bm{u}_r,~~\bm{v}_I^{\T}\bm{u}_s,~~\bm{v}_I^{\T}(\bm{\Phi} + \bm{\beta})~\right]
	\left[
	\begin{array}{c}
	(\bm{X}^r)^{\T}\bm{v}_J \\ 
	(\bm{X}^s)^{\T}\bm{v}_J \\ \bm{0}
	\end{array}
	\right]=\big[\bm{0}^{\T},~~
	\bm{0}^{\T},~~
	\bm{0}^{\T},~~
	\nonumber\\
	\bm{v}_{IX}&\underbrace{\left[\left((\bm{X}^r)^{\T}\bm{v}_J\right)a_{j,r}+\left((\bm{X}^s)^{\T}\bm{v}_J\right)a_{j,s}\right]}_{j=1,2,3,4},~\bm{v}_{IY}\underbrace{\left[\left((\bm{X}^r)^{\T}\bm{v}_J\right)a_{j,r}+\left((\bm{X}^s)^{\T}\bm{v}_J\right)a_{j,s}\right]}_{j=1,2,3,4},~\nonumber \\
	\bm{v}_{IZ}&\underbrace{\left[\left((\bm{X}^r)^{\T}\bm{v}_J\right)a_{j,r}+\left((\bm{X}^s)^{\T}\bm{v}_J\right)a_{j,s}\right]}_{j=1,2,3,4},~~I,J=1,2, \label{eps} 
	\end{align}
	where the columns are indexed by index $j$. The bending strain $\bm{\kappa}$ yields
	\begin{align}
	\bm{\kappa}=\left[
	\begin{array}{c}
	\kappa_x \\ \kappa_y \\ \kappa_{xy}
	\end{array}
	\right]=\bm{B}_{\b\varphi, 3\times24} \bm{q}+\bm{B}_{\b\triangle\varphi, 3\times4} \bm{\triangle\beta},\label{Bb}
	\end{align}
	where 
	\begin{align}
	\bm{B}_{\b\varphi, 3\times24}&=\left[
	\begin{array}{c}
	K_{11} \\ K_{22} \\ K_{12} + K_{21}
	\end{array}
	\right],\nonumber\\
	K_{IJ}&=\frac{2}{h}\left[\bm{v}_I^{\T}\bm{u}_r,~~\bm{v}_I^{\T}\bm{u}_s,~~\bm{v}_I^{\T}\bm{\Phi}~\right]
	\underbrace{\left[
		\begin{array}{c}
		(\bm{o}_r)^{\T}\bm{v}_J \\ (\bm{o}_s)^{\T}\bm{v}_J \\ (\bm{o}_t)^{\T}\bm{v}_J
		\end{array}
		\right]}_{=0\text{~for~planar~elements}}+\nonumber\\
	&~~~~\frac{2}{h}\left[\bm{v}_I^{\T}\bm{\Phi}_r,~~\bm{v}_I^{\T}\bm{\Phi}_s,~~0\right]
	\left[
	\begin{array}{c}
	(\bm{X}^r)^{\T}\bm{v}_J \\ 
	(\bm{X}^s)^{\T}\bm{v}_J \\ 
	\bm{m}^{\T}\bm{v}_J / ||\bm{m}||_2^2
	\end{array}
	\right]= \nonumber \\
	\Bigg[&(\bm{v}_I)^{\T}\underbrace{\bm{E}_{j,1}\left[\left((\bm{X}^r)^{\T}\bm{v}_J\right)a_{j,r}+\left((\bm{X}^s)^{\T}\bm{v}_J\right)a_{j,s}+\left((\bm{o}_t)^{\T}\bm{v}_J\right)a_j\right]}_{j=1,2,3,4},~\nonumber\\ 
	&(\bm{v}_I)^{\T}\underbrace{\bm{E}_{j,2}\left[\left((\bm{X}^r)^{\T}\bm{v}_J\right)a_{j,r}+\left((\bm{X}^s)^{\T}\bm{v}_J\right)a_{j,s}+\left((\bm{o}_t)^{\T}\bm{v}_J\right)a_j\right]}_{j=1,2,3,4},\nonumber \\
	&(\bm{v}_I)^{\T}\underbrace{\bm{E}_{j,3}\left[\left((\bm{X}^r)^{\T}\bm{v}_J\right)a_{j,r}+\left((\bm{X}^s)^{\T}\bm{v}_J\right)a_{j,s}+\left((\bm{o}_t)^{\T}\bm{v}_J\right)a_j\right]}_{j=1,2,3,4},~\nonumber\\
	&\frac{2}{h}\bm{v}_{I,X}\underbrace{\left[\left((\bm{o}_r)^{\T}\bm{v}_J\right)a_{j,r}+\left((\bm{o}_s)^{\T}\bm{v}_J\right)a_{j,s}\right]}_{j=1,2,3,4},\nonumber \\
	&\frac{2}{h}\bm{v}_{I,Y}\underbrace{\left[\left((\bm{o}_r)^{\T}\bm{v}_J\right)a_{j,r}+\left((\bm{o}_s)^{\T}\bm{v}_J\right)a_{j,s}\right]}_{j=1,2,3,4},\nonumber \\
	&\frac{2}{h}\bm{v}_{I,Z}\underbrace{\left[\left((\bm{o}_r)^{\T}\bm{v}_J\right)a_{j,r}+\left((\bm{o}_s)^{\T}\bm{v}_J\right)a_{j,s}\right]}_{j=1,2,3,4},~~I,J=1,2, \label{eqNormalDerivatives}
	\end{align}
	\begin{align}
	\bm{B}_{\b\triangle\varphi, 3\times4}&=\left[
	\begin{array}{c}
	L_{11} \\ L_{22} \\ L_{12} + L_{21}
	\end{array}
	\right], \nonumber\\
	L_{IJ}&=\frac{2}{h}\left[0,~~0,~~\bm{v}_I^{\T}\bm{\beta}~\right]
	\underbrace{\left[
		\begin{array}{c}
		(\bm{o}_r)^{\T}\bm{v}_J \\ (\bm{o}_s)^{\T}\bm{v}_J \\ (\bm{o}_t)^{\T}\bm{v}_J
		\end{array}
		\right]}_{=0\text{~for~planar~elements}}+\nonumber\\
	&~~~~\frac{2}{h}\left[\bm{v}_I^{\T}\bm{\beta}_r,~~\bm{v}_I^{\T}\bm{\beta}_s,~~0\right]
	\left[
	\begin{array}{c}
	(\bm{X}^r)^{\T}\bm{v}_J \\ 
	(\bm{X}^s)^{\T}\bm{v}_J \\ 
	\bm{m}^{\T}\bm{v}_J / ||\bm{m}||_2^2
	\end{array}
	\right]=\nonumber\\
	&~~~~\underbrace{(\bm{v}_I)^{\T}\bm{t}_k\left[((\bm{X}^r)^{\T}\bm{v}_J)a_{k,r}+((\bm{X}^s)^{\T}\bm{v}_J)a_{k,s}+((\bm{o}_t)^{\T}\bm{v}_J)a_k\right]}_{k=5,6,7,8},~~I,J=1,2. \label{kappa}
	\end{align}
	where the columns are indexed by index $j$ and $k$ respectively (no summation over these indices is considered). %
	\subsection{Strain energy and constitutive relations}
	The expression for the strain energy of the element reduces to the standard form:
	\begin{align}
	\Pi_{\text{e}}&=\frac{1}{2}\bm{q}^{\T}\bm{K}_{\text{e}}\bm{q},
	\end{align}
	where the vector of unknowns $\bm{q}$ is given in Eq. (\ref{qvector}).
	The stiffness matrix of an element $e$ is given by the formula
	\begin{align}
	\bm{K}_{\text{e}}&=\int_{\e} \left[\bm{B}^{\T}\bm{D}\bm{B}+\bm{B}_{\s}^{\T}\bm{D}_{\s}\bm{B}_{\s}\right]\,\text{d}A+\bm{K}_{\text{stab}}, \label{Ke}
	\end{align}
	where
	\begin{align}
	\bm{B}&=\left[
	\begin{array}{c}
	\bm{B}_{\b}\\
	\bm{B}_{\m}\\
	\end{array}
	\right],~\bm{D} = \left[
	\begin{array}{cc}
	\bm{D}_{\b}&\bm{0}\\
	\bm{0}&\bm{D}_{\m}\\
	\end{array}
	\right],\label{eq33}
	\end{align}
	and $A$ is an area of the element $e$, $\bm{B}_{\b}$ is the bending strain matrix defined later in Eq. (\ref{Bbb}), $\bm{B}_{\m}$ is the membrane strain matrix defined in Eq. (\ref{Bm}) and $\bm{K}_{\text{stab}}$ is a matrix stabilizing drilling rotations specified later in Eq. (\ref{Kdrill}). The elasticity matrices in Eq. (\ref{eq33}) equal in case of single-layered shells made of isotropic material to
	\begin{align}
	\bm{D}_{\b}&= \frac{h^3(r,s)}{12}\left[
	\begin{array}{ccc}
	\frac{E}{1-\nu^2}&\frac{\nu E}{1-\nu^2}&0\\
	\frac{\nu E}{1-\nu^2}&\frac{E}{1-\nu^2}&0\\
	0&0&G\\
	\end{array}
	\right],
	\label{const1} \\
	\bm{D}_{\m}&= h(r,s)\left[
	\begin{array}{ccc}
	\frac{E}{1-\nu^2}&\frac{\nu E}{1-\nu^2}&0\\
	\frac{\nu E}{1-\nu^2}&\frac{E}{1-\nu^2}&0\\
	0&0&G\\
	\end{array}
	\right], \label{Dm} \\
	\bm{D}_{\s}&=\frac{5}{6}h(r,s)\left[
	\begin{array}{cc}
	G & 0 \\
	0 & G \\
	\end{array}
	\right]. \label{const2} 
	\end{align}
	Actual thickness and its variation within the element is accounted for in the calculation of $\bm{D}_{\b},~\bm{D}_{\m},~\bm{D}_{\s}$ matrices. The strain resultants evaluated in the Gauss quadrature points are
	\begin{align}
	\left[
	\begin{array}{c}
	\bm{\kappa}\\
	\bm{\varepsilon}\\
	\end{array}
	\right]=\left[
	\begin{array}{c}
	\kappa_{x}\\
	\kappa_{y}\\
	\kappa_{xy}\\
	\varepsilon_{x}\\
	\varepsilon_{y}\\
	\gamma_{xy}\\
	\end{array}
	\right]=\bm{B}\bm{q},~~\bm{\gamma}=\left[
	\begin{array}{c}
	\gamma_{xz}\\
	\gamma_{yz}\\
	\end{array}
	\right]=\bm{B}_{\s}\bm{q},~
	\end{align}
	and analogously the stress resultants evaluated in the Gauss quadrature points are
	\begin{align}
	\left[
	\begin{array}{c}
	\bm{m}\\
	\bm{n}\\
	\end{array}
	\right]=\left[
	\begin{array}{c}
	m_{x}\\
	m_{y}\\
	m_{xy}\\
	n_{x}\\
	n_{y}\\
	n_{xy}\\
	\end{array}
	\right]=\bm{D}\bm{B}\bm{q},~~\bm{\mathfrak{q}}=\left[
	\begin{array}{c}
	\mathfrak{q}_{x}\\
	\mathfrak{q}_{y}\\
	\end{array}
	\right]=\bm{D}_{\s}\bm{B}_{\s}\bm{q}.
	\end{align}
	The stress resultants $\bm{m}$, $\bm{n}$ are defined in the local coordinate system defined by vectors $\bm{v}_1,~\bm{v}_2$. If it is desired to get the stress resultants in other coordinate system $\bm{v}_1',~\bm{v}_2'$, we use the transformation
	\begin{align}
	\left[\begin{array}{c}
	m_{x'}\\
	m_{y'}\\
	m_{x'y'}\\
	n_{x'}\\
	n_{y'}\\
	n_{x'y'}
	\end{array}\right]&=
	\left[\begin{array}{cc}
	\bm{\mathcal{T}}  & \bm{I}      \\
	\bm{I}            & \bm{\mathcal{T}}  
	\end{array}\right]
	\left[\begin{array}{c}
	m_x\\
	m_y\\
	m_{xy}\\
	n_x\\
	n_y\\
	n_{xy}
	\end{array}\right],~~\bm{\mathcal{T}}=\left[\begin{array}{ccc}
	T^2_{11}     & T^2_{12}     & 2T_{12}T_{11}\\
	T^2_{21}     & T^2_{22}     & 2T_{22}T_{21}\\
	T_{21}T_{11} & T_{22}T_{12} & T_{22}T_{11}+T_{21}T_{12}
	\end{array}\right],
	\end{align}
	where $T_{ij}=(\bm{v}_i')^{\T}\bm{v}_j,~i=1,2$. The bending strain matrix $\bm{B}_{\b}$ and the shear strain matrix $\bm{B}_{\s}$ are derived to depend also on supplementary degrees of freedom $\triangle\beta_{tk}$. These supplementary degrees of freedom are then eliminated with the help of the matrix $\bm{A}_{\n,4 \times 24}$, which relates them to nodal degrees of freedom as follows:
	\begin{align}
	\bm{\triangle\beta}&=\left[\begin{array}{c}
	\triangle\beta_{t5}\\
	\triangle\beta_{t6}\\
	\triangle\beta_{t7}\\
	\triangle\beta_{t8}\\
	\end{array}\right]=\bm{A}_{\n,4 \times 24}\bm{q}.
	\end{align}
	The assembly of the bending stiffness matrix is then given by
	\begin{align}
	\bm{\kappa}&=\bm{B}_{\b,3 \times 24}\bm{q}, \label{Bbb} 
	\end{align}
	where
	\begin{align}
	\bm{B}_{\b,3 \times 24} = \bm{B}_{\b\varphi, 3 \times 24} + \bm{B}_{\b\triangle \beta, 3 \times 4} \bm{A}_{\n,4 \times 24}. 
	\end{align}
	The shear stiffness matrix takes the form
	\begin{align}
	\bm{B}_{\s, 2 \times 24} &= \bm{B}_{\s\triangle \beta, 2 \times 4} \bm{A}_{\n,4 \times 24}. \label{eq4b}
	\end{align}
	
	\subsection{Evaluation of shear matrix} 
	Matrix $\bm{B}_{\s\triangle \beta, 2 \times 4}$, which relates the transversal shear strains 
	to the supplementary rotations $\triangle \beta_{tk}$, is identical to the same matrix used in the DKMQ plate element \cite{DKMQ}. We recap the derivation for completeness.
	
\begin{figure} 
		\centering
		\begin{overpic}[width=0.78\textwidth]{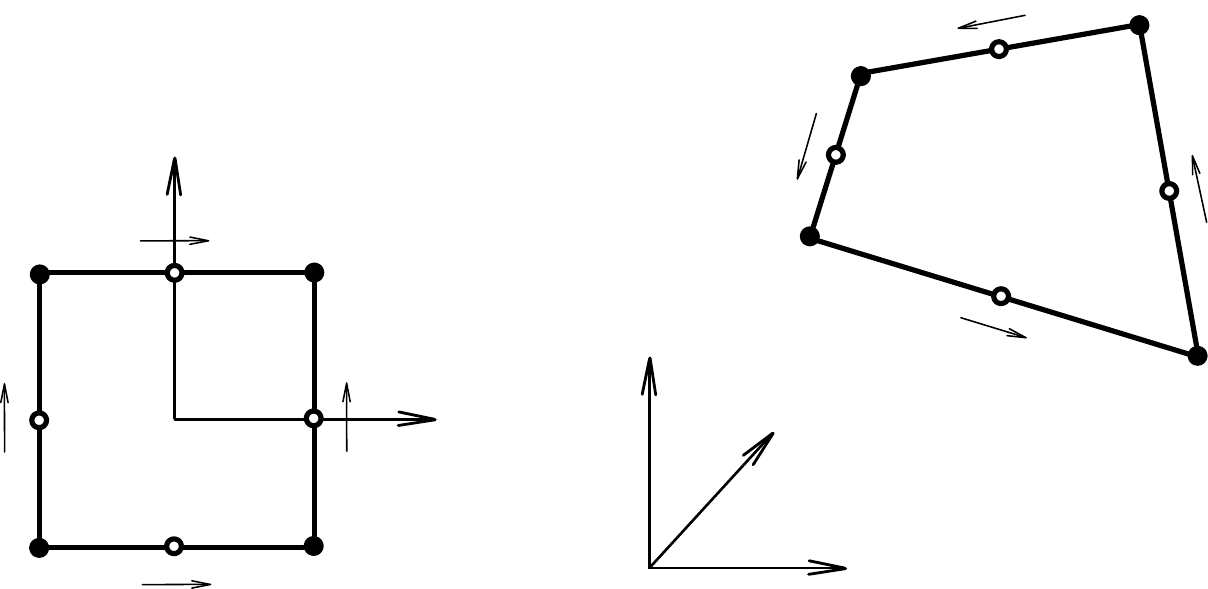} 
			\put(4.5,5.5){$1$}
			\put(13.5,5.5){$5$}
			\put(23,5.5){$2$}
			\put(4.5,13){$8$}
			\put(23,11){$6$}
			\put(4.5,23){$4$}
			\put(12,23){$7$}
			\put(23,23){$3$}
			\put(70,30){$1$}
			\put(95,22){$2$}
			\put(92,42){$3$}
			\put(72.5,39.5){$4$}	
			\put(82.5,26.5){$5$}
			\put(94,32){$6$}
			\put(82,41){$7$}
			\put(71.25,35){$8$}	
			\put(35,10.7){$r$}
			\put(11.5,34){$s$}
			\put(71.5,0.5){$X$}
			\put(65,12.7){$Y$}
			\put(52.65,21.2){$Z$}
			\put(4.45,-0.5){$\overline{\gamma}_{rz5}$}
			\put(30,17){$\overline{\gamma}_{sz6}$}
			\put(16.5,31.2){$\overline{\gamma}_{rz7}$}
			\put(-7,14){$\overline{\gamma}_{sz8}$}
			\put(79,17.5){$\overline{\gamma}_{tz5}$} 		
			\put(101.5,32.75){$\overline{\gamma}_{tz6}$}	
			\put(80.5,49.5){$\overline{\gamma}_{tz7}$}    
			\put(59.5,36.5){$\overline{\gamma}_{tz8}$}		
		\end{overpic}
		\caption{Shear strains in the reference and local configurations} 
		\label{DKMQshear}
	\end{figure}
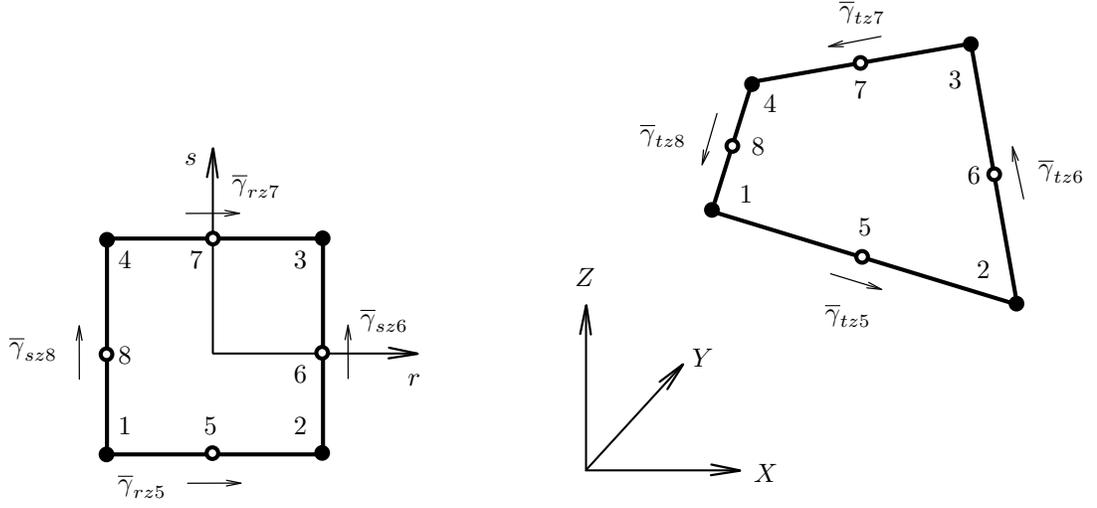

	Let us evaluate an average transversal shear strain on each side. Taking an arbitrary side $k$, the directional parameter is denoted by $t$. The transversal shear strain in direction $t$ has the form
	\begin{align}
	\overline{\gamma}_{tz}&=\frac{\mathfrak{q}_t}{D_{\s}}.
	\end{align}
	where $\mathfrak{q}_t$ is a shear internal force and $D_{\s}$ is the shear stiffness, which in case of single-layered isotropic shells equals to $\frac{5}{6}Gh$, where $G$ is the shear modulus and $h$ is an average shell thickness on side $k$. The stress equilibrium yields
	\begin{align}
	\mathfrak{q}_{t}&=m_{t,t}+m_{tm,m}, 
	\end{align}
	and Hooke's law has the form
	\begin{align}
	\left[
	\begin{array}{c}
	m_t \\ m_m \\ m_{tm}
	\end{array}
	\right]&=
	\left[
	\begin{array}{ccc}
	D_{\b} & * & * \\ 
	* & * & * \\
	* & * & * \\
	\end{array}
	\right]
	\left[
	\begin{array}{c}
	\beta_{t,t} \\ \beta_{m,m} \\ \beta_{t,m}+\beta_{m,t}
	\end{array}
	\right],
	\end{align}
	where $D_{\b}$ is an average bending stiffness on the side $k$, which in case of single-layered isotropic shell equals to $\frac{Eh^3}{12(1-\nu^2)}$, where $E$ is Young's modulus and $h$ the shell thickness.
	We now consider a linear approximation of angle $\beta$ in normal $m$-direction and quadratic in tangential $t$-direction according to Eq. (\ref{beta}). 
	Both approximations do not depend on the $m$-direction, which immediately yields
	\begin{align}
	\overline{\gamma}_{tz} &= \frac{D_{\b}}{D_{\s}} \beta_{t,tt}(t) = -\frac{2}{3} \Phi_k \triangle\beta_{tk},   \label{eq2220}
	\end{align}
	where
	\begin{align}
	\Phi_k=\frac{12}{L_k^2} \frac{D_{\b}}{D_{\s}}=\frac{12}{5(1-\nu)}\left(\frac{h}{L_k}\right)^2.
	\end{align}
	We see that the transversal shear approximation $\overline{\gamma}_{tz}$ is constant on the side.
	If we liked to transform these shear strains back to the reference plane, we would have to incorporate the determinant of the transformation and also the parametrization sign change, see to Figure \ref{DKMQshear}:
	\begin{align}
	\overline{\gamma}_{\xi z, 5} &= (\det J)_{5} \overline{\gamma}_{tz, 5}=+\frac{L_{5}}{2}\left(-\frac{2}{3}\Phi_{5}\triangle\beta_{t5} \right), \\
	\overline{\gamma}_{\eta z, 6} &= (\det J)_{6} \overline{\gamma}_{tz, 6}=+\frac{L_{6}}{2}\left(-\frac{2}{3}\Phi_{6}\triangle\beta_{t6} \right), \\
	\overline{\gamma}_{\xi z, 7} &= (\det J)_{7} \overline{\gamma}_{tz, 7}=-\frac{L_{7}}{2}\left(-\frac{2}{3}\Phi_{7}\triangle\beta_{t7} \right), \\
	\overline{\gamma}_{\eta z, 8} &= (\det J)_{8} \overline{\gamma}_{tz, 8}=-\frac{L_{8}}{2}\left(-\frac{2}{3}\Phi_{8}\triangle\beta_{t8} \right).
	\end{align}
	As in the MITC4 element \cite{MITC4} the transverse shear strains are approximated with
	\begin{align}
	\overline{\gamma}_{rz} &= \frac{1}{2}(1-s)\overline{\gamma}_{tz5}+\frac{1}{2}(1+s)\overline{\gamma}_{tz7}, \\
	\overline{\gamma}_{sz} &= \frac{1}{2}(1-r)\overline{\gamma}_{tz8}+\frac{1}{2}(1+r)\overline{\gamma}_{tz6}, 
	\end{align}
	where
	\begin{align}
	\left[
	\begin{array}{c}
	\overline{\gamma}_{xz} \\ \overline{\gamma}_{yz}
	\end{array}
	\right]&=\bm{J}^{-1}
	\left[
	\begin{array}{c}
	\overline{\gamma}_{rz} \\ \overline{\gamma}_{sz}
	\end{array}
	\right],
	\end{align}
	so we finally get 
	\begin{align}
	&\left[
	\begin{array}{c}
	\overline{\gamma}_{xz} \\ \overline{\gamma}_{yz}
	\end{array}
	\right]=\bm{B}_{\s\triangle \beta, 2 \times 4}
	\left[
	\begin{array}{c}
	\triangle \beta_{t5}  \rule[-2mm]{0mm}{6mm}\\
	\triangle \beta_{t6}  \rule[-2mm]{0mm}{6mm}\\
	\triangle \beta_{t7}  \rule[-2mm]{0mm}{6mm}\\
	\triangle \beta_{t8}  \rule[-2mm]{0mm}{6mm}\\
	\end{array}
	\right],~\\
	&\bm{B}_{\s\triangle \beta, 2 \times 4}=\frac{1}{6}\left[
	\begin{array}{cccc}
	\!\!\!\!-J^{-1}_{11}(1\!-\!s)L_5\Phi_5 & 
	\!-J^{-1}_{12}(1\!+\!r)L_6\Phi_6 & 
	\!J^{-1}_{11}(1\!+\!s)L_7\Phi_7 & 
	\!J^{-1}_{12}(1\!-\!r)L_8\Phi_8 \!\!\rule[-2mm]{0mm}{6mm}\\
	\!\!\!\!-J^{-1}_{21}(1\!-\!s)L_5\Phi_5 & 
	\!-J^{-1}_{22}(1\!+\!r)L_6\Phi_6 & 
	\!J^{-1}_{21}(1\!+\!s)L_7\Phi_7 & 
	\!J^{-1}_{22}(1\!-\!r)L_8\Phi_8  \!\!\rule[-2mm]{0mm}{6mm}
	\end{array}
	\right]. \nonumber
	\end{align}
	
	\FloatBarrier
	\subsection{Evaluation of matrix \textbf{A$_{\text{n}}$}} 
	We use the Hu-Washizu functional according to which
	\begin{align}
	\int_0^{L_k} (\gamma_{tz}-\overline{\gamma}_{tz})\,\text{d}t=0,
	\end{align}
	where $\gamma_{tz}$ is a shear approximated from displacement and $\overline{\gamma}_{tz}$ is an average shear approximated by Eq. (\ref{eq2220}). We have
	\begin{align}
	\gamma_{tz} &= \PART{u_z}{t}+\beta_t=\nonumber\\
	&\frac{u_{jz}-u_{iz}}{L_k} + \left(1-\frac{t}{L_k}\right)\beta_{ti}+\frac{t}{L_k}\beta_{tj}+4\frac{t}{L_k}\left(1-\frac{t}{L_k}\right)\triangle\beta_{tk},\\
	\overline{\gamma}_{tz} &= -\frac{2}{3}\Phi_k\triangle\beta_{tk},\\
	0&=\int_0^{L_k} (\gamma_{tz}-\overline{\gamma}_{tz})\,\text{d}t=u_{jz}-u_{iz}+\frac{L_k}{2}\left(\beta_{ti}+\beta_{tj}\right) + \frac{2}{3} L_k(1+\Phi_k) \triangle \beta_{tk}. \nonumber
	\end{align}
	By rewriting the scalar unknowns to vector ones we get  
	\begin{align}
	\bm{n}_k^{\T}[\bm{u}_i-\bm{u}_j]-\frac{L_k}{2}\bm{l}_k^{\T}\left(\bm{\varphi}_{i}+\bm{\varphi}_{j}\right) &=\frac{2}{3} L_k(1+\Phi_k) \triangle \beta_{tk}. \label{eq44}
	\end{align}
	In the matrix form the Eq. (\ref{eq44}) takes the form
	\begin{align}
	&\underbrace{\left[\bm{A}_{\w}^{\text{I}}~\bm{A}_{\w}^{\text{II}}~\bm{A}_{\w}^{\text{III}}~\bm{A}_{\w}^{\text{IV}}~\bm{A}_{\w}^{\text{V}}~\bm{A}_{\w}^{\text{VI}}\right]}_{\bm{A}_{\w, 4 \times 24}}
	\bm{q}=\\
	&\underbrace{\left[
		\begin{array}{cccc}
		\frac{2}{3}L_5(1\!+\!\Phi_5) &&& \\
		& \frac{2}{3}L_6(1\!+\!\Phi_6) && \\
		&& \frac{2}{3}L_7(1\!+\!\Phi_7) & \\
		&&& \frac{2}{3}L_8(1\!+\!\Phi_8) \\
		\end{array}
		\right]}_{\bm{A}_{\triangle \beta,4 \times 4}}
	\underbrace{\left[
		\begin{array}{c}
		\triangle\beta_{t5} \\
		\triangle\beta_{t6} \\
		\triangle\beta_{t7} \\
		\triangle\beta_{t8} \\
		\end{array}
		\right]}_{\bm{\triangle\beta}},\nonumber\\
	\bm{d}_{k}&=-\frac{L_k}{2}\bm{l}_k,~k=5,6,7,8, \nonumber\\
	\bm{n}_k&=\frac{\bm{n}_i+\bm{n}_j}{2},~k=5,6,7,8,\nonumber\\
	\bm{A}_{\w}^{\text{I}} &=\left[\begin{array}{cccc}
	d_{5X} & d_{5X} & 0 & 0  \rule[-2mm]{0mm}{6mm}\\
	0 & d_{6X} & d_{6X} & 0  \rule[-2mm]{0mm}{6mm}\\
	0 & 0 & d_{7X} & d_{7X}  \rule[-2mm]{0mm}{6mm}\\     
	d_{8X} & 0 & 0 & d_{8X}  \rule[-2mm]{0mm}{6mm}\\
	\end{array}\right],~\bm{A}_{\w}^{\text{II}}=\left[\begin{array}{cccc}
	d_{5Y} & d_{5Y} & 0 & 0  \rule[-2mm]{0mm}{6mm}\\
	0 & d_{6Y} & d_{6Y} & 0  \rule[-2mm]{0mm}{6mm}\\
	0 & 0 & d_{7Y} & d_{7Y}  \rule[-2mm]{0mm}{6mm}\\     
	d_{8Y} & 0 & 0 & d_{8Y}  \rule[-2mm]{0mm}{6mm}\\
	\end{array}\right], \nonumber\\
	\bm{A}_{\w}^{\text{III}}&=\left[\begin{array}{cccc}
	d_{5Z} & d_{5Z} & 0 & 0  \rule[-2mm]{0mm}{6mm}\\
	0 & d_{6Z} & d_{6Z} & 0  \rule[-2mm]{0mm}{6mm}\\
	0 & 0 & d_{7Z} & d_{7Z}  \rule[-2mm]{0mm}{6mm}\\     
	d_{8Z} & 0 & 0 & d_{8Z}  \rule[-2mm]{0mm}{6mm}\\
	\end{array}\right],~\bm{A}_{\w}^{\text{IV}}= 
	\left[\begin{array}{cccccccccccc}
	n_{5X} &\!\!\!-n_{5X} & 0 & 0    \rule[-2mm]{0mm}{6mm}\\
	0 & n_{6X} &\!\!\!-n_{6X} & 0    \rule[-2mm]{0mm}{6mm}\\
	0 & 0 & n_{7X} &\!\!\!-n_{7X}    \rule[-2mm]{0mm}{6mm}\\     
	\!\!\!-n_{8X}& 0 & 0 & n_{8X}    \rule[-2mm]{0mm}{6mm}\\
	\end{array}\right], \nonumber\\
	\bm{A}_{\w}^{\text{V}}&=	\left[\begin{array}{cccccccccccc}
	n_{5Y} &\!\!\!-n_{5Y} & 0 & 0    \rule[-2mm]{0mm}{6mm}\\
	0 & n_{6Y} &\!\!\!-n_{6Y} & 0    \rule[-2mm]{0mm}{6mm}\\
	0 & 0 & n_{7Y} &\!\!\!-n_{7Y}    \rule[-2mm]{0mm}{6mm}\\     
	\!\!\!-n_{8Y}& 0 & 0 & n_{8Y}    \rule[-2mm]{0mm}{6mm}\\
	\end{array}\right],~\bm{A}_{\w}^{\text{VI}}=	\left[\begin{array}{cccccccccccc}
	n_{5Z} &\!\!\!-n_{5Z} & 0 & 0    \rule[-2mm]{0mm}{6mm}\\
	0 & n_{6Z} &\!\!\!-n_{6Z} & 0    \rule[-2mm]{0mm}{6mm}\\
	0 & 0 & n_{7Z} &\!\!\!-n_{7Z}    \rule[-2mm]{0mm}{6mm}\\     
	\!\!\!-n_{8Z}& 0 & 0 & n_{8Z}    \rule[-2mm]{0mm}{6mm}\\
	\end{array}\right], \nonumber
	\end{align}
	where $\bm{n}_k$ is the normal vector corresponding to the side $k$ starting with node $i$ and ending with node $j$. 
	Finally, the supplementary rotations $\triangle\beta_{tk}$ can be evaluated by the following equation
	\begin{align}
	\bm{\triangle\beta}=\bm{A}_{\n, 4 \times 24}
	\bm{q},~~\bm{A}_{\n, 4 \times 24}=\bm{A}_{\triangle \beta,4 \times 4}^{-1}\bm{A}_{\w, 4 \times 24}.
	\end{align}
	
	\subsection{Drilling rotations stabilization}
	The shell element must be stabilized with respect to drilling rotations $\varphi_z$ in order to avoid a singular stiffness matrix, which appears if a node connects coplanar elements. We follow the stabilization procedure given in \cite{Katili2015}. The added energy, which stabilizes drilling rotation has the form 
	\begin{align}
	\Pi_\text{stab} &= \frac{c}{2} \left[Gh\int_{\e} \,\varphi_{z}^2(r,s)\,\text{d}A + \frac{Eh^3}{12}\int_{\e} \left[\varphi_{z,x}^2(r,s) + \varphi_{z,y}^2(r,s)\right]\,\text{d}A\right], \label{phiz}
	\end{align}
	where $c$ is the penalty constant. The value $c=0.001$ is considered in our implementation of the DKMQ24 element as well as in \cite{Katili2015}. The corresponding stiffness matrix in the global coordinate system has the form:
	\begin{align}
	& K_{\text{stab};8+i,8+j} =c \int_{-1}^1 \int_{-1}^1 \Big\{Gh\,a_i(r,s) a_j(r,s) + \frac{Eh^3}{12}\big[a_{i,x}(r,s)a_{j,x}(r,s)+\nonumber \\ &a_{i,y}(r,s)a_{j,y}(r,s)\big]\Big\}\,
	(\bm{n}_i)^{\T}\bm{n}_j\,\text{det}J~\text{d}r\text{d}s,~~i,j=1,\dots,4. \label{Kdrill}
	\end{align}
	The first term is called the MacNeal stabilization and has rank 1, the second has rank 3. It can be proven that such a matrix removes all zero energy modes connected with drilling rotations \cite{Katili2015}.
	
	
	\subsection{Numerical integration}
	\label{integration}
	Elements are integrated by means of the standard $2\times 2$ Gauss quadrature except of the stabilization term (\ref{Kdrill}), for which a reduced one-point integration rule is used.

	\subsection{Load vector}
	\label{LoadVectorSection}
	The work of external forces for an arbitrary element $e$ in the discretized form yields
	\begin{align}
	\Pi_{\text{ext,e}}=\bm{f}^{\T}\bm{q},~\bm{f}=\left[
	\begin{array}{c}
	m_{1X} \\
	\vdots \\
	m_{1Y} \\
	\vdots \\
	m_{1Z} \\
	\vdots \\
	f_{1X} \\
	\vdots \\
	f_{1Y} \\
	\vdots \\
	f_{1Z} \\
	\vdots \\
	\end{array}\right].
	\end{align}
	In order to evaluate the consistent load vector $\bm{f}$ let us evaluate the deformation of the middle surface using Eq. (\ref{eq15}) by setting $t=0$, which yields
	\begin{align}
	\left[
	\begin{array}{c}
	u_X \\ u_Y \\ u_Z
	\end{array}
	\right](r,s,0)=&\sum_{i=1}^4a_i(r,s)
	\left[
	\begin{array}{c}
	u_{iX} \\ u_{iY} \\ u_{iZ}
	\end{array}
	\right]. \label{middSurf1}
	\end{align}
	Let us further consider a general loading case given by tangential pressures $p_x(r,s)$ and $p_y(r,s)$, by the normal pressure $p(r,s)$ and by the element loading force $\bm{F}_{\e}$. The consistent nodal force load vector in the global coordinate system is given by the integration of the basis functions given by Eq. (\ref{middSurf1}) and yields
	\begin{align}
	\bm{f}_i=&\int_{\e} \bm{P}(r,s)a_i(r,s)\,\dA=
	\sum_{q=1}^4  \bm{P}(r_q,s_q) a_i(r_q,s_q)w_q \det J(r_q,s_q), ~i=1,2,3,4,
	\end{align}
	where 
	\begin{align}
	\bm{P}(r,s)=&p_x(r,s)\bm{v}_1(r,s)+p_y(r,s)\bm{v}_2(r,s)-p(r,s)\bm{n}(r,s)+\frac{\bm{F}_{\e}}{A}, \label{Prs}
	\end{align}
	where $\bm{v}_1(r,s), \bm{v}_2(r,s), \bm{v}_3(r,s)\equiv\bm{n}(r,s)$ are the basis vectors, $A$ is the element area and $w_q,~q=1,2,3,4$ are the weights and $r_q,\,s_q,~q=1,2,3,4$ are the quadrature points of the Gauss quadrature $2\times 2$.

	\section{Improvements to the DKMQ24 shell element}
	
	\label{NewModif}
	
	In this Section we introduce modifications of the DKMQ24 shell element, derived in Section 2. These modifications apply both to the DKMQ24$_1$+ and DKMQ24$_2$+ shell elements. The difference between these two variants is explained in Section \ref{penaltySection}.
	
	\subsection{The treatment of drilling rotations} In the original derivation drilling rotations are stabilized according to formula (\ref{phiz}), which causes the local element stiffness matrix to have a full rank even for the coplanar element configuration. However, this treatment keeps drilling rotations independent from the remaining degrees of freedom and as a consequence the element behaviour is not improved although we pay for these additional degrees of freedom by an increase of computational costs. Therefore, we use the approach proposed in \cite{OPQ}, where the drilling rotations are considered as derivatives of the in-plane displacement field. The space of basis functions is enriched by incomplete quadratic polynomials in this case. 
	
	Let us start with calculation of a rotation around a normal vector $\mathbf{n}$ at an arbitrary point $(r,s)$ 
	\begin{align}
	\varphi_z(r,s)=\bm{n}^{\T} \left[
	\begin{array}{c}
	\varphi_{X}(r,s) \\ \varphi_{Y}(r,s) \\ \varphi_{Z}(r,s)
	\end{array}
	\right]. \label{phiz3}
	\end{align}
	The displacement field $\bm{u}=\def\arraystretch{0.5}\left[
	\begin{array}{c}
	u_X(r,s) \\ u_Y(r,s) \\ u_Z(r,s)
	\end{array}
	\right]$ given by formula (\ref{eq15}) is then enriched by the following term
	\begin{align}
	\triangle\bm{u}=&\sum_{k=5}^8 a_k(r,s) \frac{L_k}{8}\left(\varphi_{J(k)z}-\varphi_{I(k)z}\right)\left(-\bm{l}_k\right){=}_{(\ref{phiz3})}\nonumber \\
	&\sum_{k=5}^8 a_k(r,s) \, \bm{l}_k \frac{L_k}{8}\bm{n}^{\T} \left(\left[
	\begin{array}{c}
	\varphi_{I(k)X} \\ \varphi_{I(k)Y} \\ \varphi_{I(k)Z}
	\end{array}
	\right]-\left[
	\begin{array}{c}
	\varphi_{J(k)X} \\ \varphi_{J(k)Y} \\ \varphi_{J(k)Z}
	\end{array}
	\right]\right), \label{eq223} 
	\end{align}
	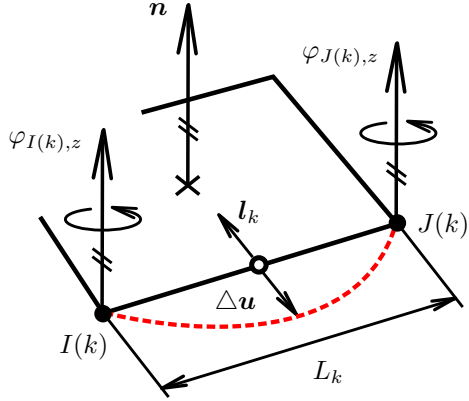
\begin{figure}[H]
		\centering
		\begin{overpic}[width=0.35\textwidth]{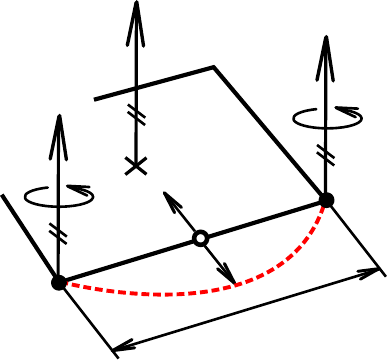} 
			\put(26,90){$\bm{n}$}
			\put(-7,60){$\varphi_{I(k),z}$}
			\put(62,80){$\varphi_{J(k),z}$}
			\put(47,42){$\bm{l}_k$}
			\put(41.5,21){$\triangle\bm{u}$}
			\put(5,10.5){$I(k)$}		
			\put(89,39){$J(k)$}		
			\put(64.5,4.5){$L_k$}		
		\end{overpic}
		\caption{Enrichment of the in-plane displacement field with the help of drilling rotations on side $k$} 
		\label{Picturephiz}
	\end{figure}
	\FloatBarrier
	where the starting node $I(k)$ and the ending node $J(k)$ for the side $k=5,6,7,8$ are given by 
	\begin{align}
	\begin{tabular}{c||c|c}
	$k$ & $I(k)$ & $J(k)$ \\
	\hline
	5 & 1 & 2 \\
	6 & 2 & 3 \\
	7 & 3 & 4 \\
	8 & 4 & 1 \\
	\end{tabular} \nonumber 
	\end{align}
	and $\bm{l}_k$ is the in-plane normal vector corresponding to the side $k$, see Figure \ref{nodalunknowns}, given by
	$\bm{l}_k=\frac{\bm{n}_k\times\bm{t}_k}{||\bm{n}_k\times\bm{t}_k||}$, where $\bm{n}_k=\frac{\bm{n}_{I(k)}+\bm{n}_{J(k)}}{2},~\bm{t}_k=\frac{P_{J(k)} - P_{I(k)}}{|P_{J(k)} - P_{I(k)}|}$, see Figure \ref{Picturephiz}. 
	The membrane strain field $\bm{\varepsilon}=\def\arraystretch{0.5}\left[
	\begin{array}{c}
	\varepsilon_x(r,s) \\ \varepsilon_y(r,s) \\ \gamma_{xy}(r,s)
	\end{array}
	\right]$ given by formula (\ref{eps0}) is enriched by the term 
	$\def\arraystretch{0.5}\left[
	\begin{array}{c}
	F_{11} \\ F_{22} \\ F_{12} + F_{21}
	\end{array}
	\right],$
	where
	\begin{align}
	F_{KL}=&\,\bm{v}_K^{\T}(r,s)\PART{(\triangle u_X, \triangle u_Y, \triangle u_Z)}{(r,s,t)}\PART{(r,s,t)}{(X,Y,Z)}\bm{v}_L(r,s)=\nonumber\\
	&\sum_{k=5}^8 \frac{L_k}{8}(\bm{v}_K^{\T}\bm{l}_k)\left[\left((\bm{X}^r)^{\T}\bm{v}_L\right)a_{k,r}(r,s)+\left((\bm{X}^s)^{\T}\bm{v}_L\right)a_{k,s}(r,s)\right] N(k), ~~K,L=1,2, 
	\end{align}
	where
	\begin{align}
	N(k)&:=\bm{n} ^{\T} \left(\left[
	\begin{array}{c}
	\varphi_{I(k)X} \\ \varphi_{I(k)Y} \\ \varphi_{I(k)Z}
	\end{array}
	\right]-\left[
	\begin{array}{c}
	\varphi_{J(k)X} \\ \varphi_{J(k)Y} \\ \varphi_{J(k)Z}
	\end{array}
	\right]\right) = \nonumber \\
	&n_X\mathbf{e}_{I(k)}-n_X\mathbf{e}_{J(k)}+
	n_Y\mathbf{e}_{4+I(k)}-n_Y\mathbf{e}_{4+J(k)}+
	n_Z\mathbf{e}_{8+I(k)}-n_Z\mathbf{e}_{8+J(k)},
	\end{align}
	where the vector $\mathbf{e}_i$ is a zero vector with 1 at index $i$.
	The bending strain field $\bm{\kappa}$ given by Eq. (\ref{Bb}) stays unchanged. To conclude this Section, let us remark that the DKMQ24 shell element has improved bending behaviour, due to quadratic enrichment of the displacement field with the help of out-of-plane rotations $\varphi_x, \varphi_y$ and in the same manner, the DKMQ24$_1$+ and DKMQ24$_2$+ shell elements have improved membrane behaviour due to quadratic enrichment of the displacement field with the help of drilling rotations $\varphi_z$.
	
	\subsection{The static condensation}
	We enrich the displacement field by a bubble mode according to the approach given in \cite{OPQ}. Two new degrees of freedom $u_1,u_2$ corresponding to the mid-element node $9$ are considered:
	\begin{align}
	\bm{u}_{\n}(r,s) &=a_9(r,s) \left[u_1\bm{v}_1(0,0)+u_2\bm{v}_2(0,0)\right],~~a_9(r,s)=(1-r^2)(1-s^2).
	\end{align}
	The corresponding strain vector has the form
	\begin{align}
	\bm{\varepsilon}_{\n}&=\bm{B}_{\n,3\times2}\left[\begin{array}{c}
	u_1\\
	u_2
	\end{array}\right]. 
	\end{align}
	These additional degrees of freedom are statically condensed out on the element level, therefore the shell element has again 24 degrees of freedom. The static condensation is based on a standard procedure applied on the membrane part only (Appendix A). The statically condensed strain matrix $\bm{B}$ takes the form
	\begin{align}
	\bm{B}=\left[\begin{array}{c}
	\bm{B}_{\b}\\
	\bm{B}_{\m}-\bm{B}_{\n,3\times2}\bm{K}_{\n\n}^{-1}(\bm{K}_{\m\n})^{\T}
	\end{array}\right],
	\end{align}
	where
	\begin{align}
	\bm{K}_{\m\n}&=\int_{\text{e}} (\bm{B}_{\m,3\times24})^{\T} \bm{D}_{\m} \bm{B}_{\n,3\times2}~\text{d}A, \\
	\bm{K}_{\n\n}&=\int_{\text{e}} (\bm{B}_{\n,3\times2})^{\T} \bm{D}_{\m} \bm{B}_{\n,3\times2}~\text{d}A.
	\end{align}
	
	\subsection{Penalty term}
	\label{penaltySection}
	The drilling rotations stabilization term (\ref{phiz}), used in the case of the DKMQ24 shell element, is not used in case of the DKMQ24$_1$+ and DKMQ24$_2$+ shell elements. However, a penalty term which forces the drilling rotations $\varphi_z(r,s)$ to be equal to a rotation of the displacement field $\psi$ at the same point, derived by Hughes and Brezzi \cite{OPQ}, is considered. This energy term takes the form
	\begin{align}
	\Pi_{\text{stab}} &= \frac{c_1Gh}{2}\int_{\e} \left[\psi(r,s)-\varphi_z(r,s)\right]^2\,\text{d}A,  
	\label{phiz2}
	\end{align}
	where
	\begin{align}
	\psi(r,s)         &=\frac{1}{2}\left(\PART{u_y}{x}(r,s)-\PART{u_x}{y}(r,s)\right), \\
	\PART{u_y}{x}(r,s)&=\sum_{i=1}^4 a_{i,x}(r,s)\, (\bm{v}_2(r,s))^{\T} \left[\begin{array}{c}u_{iX}\\u_{iY}\\u_{iZ}\end{array}\right], \\
	\PART{u_x}{y}(r,s)&=\sum_{i=1}^4 a_{i,y}(r,s)\, (\bm{v}_1(r,s))^{\T} \left[\begin{array}{c}u_{iX}\\u_{iY}\\u_{iZ}\end{array}\right],\\
	a_{i,x}(r,s)&=a_{i,r}(r,s)(\bm{X}^r(r,s))^{\T}\bm{v}_1(r,s)+a_{i,s}(r,s)(\bm{X}^s(r,s))^{\T}\bm{v}_1(r,s),\\
	a_{i,y}(r,s)&=a_{i,r}(r,s)(\bm{X}^r(r,s))^{\T}\bm{v}_2(r,s)+a_{i,s}(r,s)(\bm{X}^s(r,s))^{\T}\bm{v}_2(r,s),\\
	\varphi_z(r,s)&=\sum_{i=1}^4 a_{i}(r,s)\, (\bm{v}_3(r,s))^{\T} \left[\begin{array}{c}\varphi_{iX}\\\varphi_{iY}\\\varphi_{iZ}\end{array}\right],\label{phiz2b}
	\end{align}
	and $c_1$ is the penalty constant and is taken as $c_1=0.15$ in our implementation. The corresponding stiffness matrix in the global coordinate system has the form:
	\begin{align}
	& K_{\text{stab};i,j} = c_1 Gh\int_{-1}^1 \int_{-1}^1 \, g_{i}(r,s) \, g_{j}(r,s) \,\text{det}J~\text{d}r\text{d}s,~~i,j=1,\dots,24, \label{Kdrill2}  
	\end{align}
	where the vector $\bm{g}$ has length 24 and takes the form
	\begin{align}
	\bm{g}^{\T}(r,s)=\Bigg[&-a_1(r,s)v_{3X}(r,s),\dots,-a_4(r,s)v_{3X}(r,s), \nonumber \\
	&-a_1(r,s)v_{3Y}(r,s),\dots,-a_4(r,s)v_{3Y}(r,s), \nonumber \\
	&-a_1(r,s)v_{3Z}(r,s),\dots,-a_4(r,s)v_{3Z}(r,s), \nonumber \\
	&\frac{a_{1,x}(r,s)v_{2X}(r,s)-a_{1,y}(r,s)v_{1X}(r,s)}{2},\dots,~\frac{a_{4,x}(r,s)v_{2X}(r,s)-a_{4,y}(r,s)v_{1X}(r,s)}{2}, \nonumber \\
	&\frac{a_{1,x}(r,s)v_{2Y}(r,s)-a_{1,y}(r,s)v_{1Y}(r,s)}{2},\dots,~\frac{a_{4,x}(r,s)v_{2Y}(r,s)-a_{4,y}(r,s)v_{1Y}(r,s)}{2}, \nonumber \\
	&\frac{a_{1,x}(r,s)v_{2Z}(r,s)-a_{1,y}(r,s)v_{1Z}(r,s)}{2},\dots,~\frac{a_{4,x}(r,s)v_{2Z}(r,s)-a_{4,y}(r,s)v_{1Z}(r,s)}{2}\Bigg].
	\end{align}
	The one point integration rule is used to integrate the stiffness matrix (\ref{Kdrill2}) in order to relax the prescribed condition (\ref{phiz2}).
	
	We use either a constant penalty parameter $c_1$ (let us call this element variant DKMQ24$_1$+) or we propose the following scaling
	\begin{align}
	c_1=c_2\frac{h}{\sqrt{A}} \label{penalty2} 
	\end{align}
	where $h$ is the element thickness and $A$ is the element area. This makes the scaling dependent only on the local information and independent of length units. Moreover, this scaling increase the penalty constant for thicker and smaller elements (let us call this element variant DKMQ24$_2$+).
	
	\subsection{Numerical integration of shear terms}
	By using the following explicit notation of the shear strain matrix rows
	\begin{align}
	\bm{B}_{\s}=\left[\begin{array}{c}
	\bm{b}_{\s x}^{\T} \\
	\bm{b}_{\s y}^{\T} \\
	\end{array}
	\right],
	\end{align}
	we can decompose the element stiffness matrix $\bm{K}_{\text{e}}$, given by Eq. (\ref{Ke}), to the form
	\begin{align}
	\bm{K}_{\text{e}}&=\int_{\e}\bm{B}^{\T}\bm{D}\bm{B}\,\text{d}A+\frac{5}{6}G\int_{\e}h\bm{b}_{\s x}\bm{b}_{\s x}^{\T}\,\text{d}A+\frac{5}{6}G\int_{\e}h\bm{b}_{\s y}\bm{b}_{\s y}^{\T}\,\text{d}A+\bm{K}_{\text{stab}}, \label{Ke2}
	\end{align}
	where the stabilisation stiffness matrix $\bm{K}_{\text{stab}}$ is given by Eq. (\ref{Kdrill2}). We have numerically tested different quadrature rules for each term in Eq. (\ref{Ke2}) and we suggest the reduced integration of shear terms, summarized in Table 3. The proposed reduced integration provides smaller numerical errors when compared to the $2\times 2$ quadrature rule, mainly at the pinched cylinder benchmark in variant $h=0.3\,\text{m}$ (Section \ref{b4}). Also, the stiffness matrix assembling is slightly faster, because the number of quadrature points for the shear terms is reduced to half. 
	
	\begin{center}
		\begin{table}
			\centering
			\begin{tabular}{cc}
				Term      & Gauss integration rule            \rule[-2mm]{0mm}{6mm} \\ 
				\hline
				$\int_{\e}\bm{B}^{\T}\bm{D}\bm{B}\,\text{d}A$       & $2\times2$     \rule[-2mm]{0mm}{6mm} \\
				$\frac{5}{6}G\int_{\e}h\bm{b}_{\s x}\bm{b}_{\s x}^{\T}\,\text{d}A$       &  $1\times2$  \rule[-2mm]{0mm}{6mm} \\
				$\frac{5}{6}G\int_{\e}h\bm{b}_{\s y}\bm{b}_{\s y}^{\T}\,\text{d}A$       &  $2\times1$    \rule[-2mm]{0mm}{6mm} \\
				$\bm{K}_{\text{stab}}$       &  $1\times1$    \rule[-2mm]{0mm}{6mm} \\
			\end{tabular} \nonumber
			\label{tabshear}
			\caption{Quadrature rules used for the stiffness matrix assembling}
		\end{table}
	\end{center}
	
	
	

	\subsection{Nodal moment corrections}
	\label{LoadVectorSection2}
	Note that due to absence of rotations in Eq. (\ref{middSurf1}) the load vector does not contain any nodal moment corrections. The reason is the Mindlin–Reissner hypothesis of the formulation, which separates high order terms of the displacement to the rotations. The real (physical) displacement of the element middle shell surface can be, analogously to Eq. (\ref{eq223}), considered in the following way
	\begin{align}
	u_z(r,s,0)=\bm{n}^{\T}
	\left[
	\begin{array}{c}
	u_X \\ u_Y \\ u_Z
	\end{array}
	\right]=&\sum_{i=1}^4a_i(r,s)
	\bm{n}^{\T}\left[
	\begin{array}{c}
	u_{iX} \\ u_{iY} \\ u_{iZ}
	\end{array}
	\right]+\sum_{k=5}^8a_k(r,s)\frac{L_k}{8}\bm{l}_k^{\T}
	\left(\left[
	\begin{array}{c}
	\varphi_{I(k)X} \\ \varphi_{I(k)Y} \\ \varphi_{I(k)Z}
	\end{array}
	\right]-\left[
	\begin{array}{c}
	\varphi_{J(k)X} \\ \varphi_{J(k)Y} \\ \varphi_{J(k)Z}
	\end{array}
	\right]\right). \label{middSurf2}
	\end{align}
	Using the normal projection of the pressure and force loading introduced in Eq. (\ref{Prs}), denoted by $\lambda(r,s)$, given by formula
	\begin{align}
	\lambda(r,s):=\bm{n}^{\T}(r,s)\bm{P}(r,s)=-p(r,s)+\frac{\bm{n}^{\T}(r,s)\bm{F}_{\e}}{A},
	\end{align}
	and considering the work of the normal forces $\delta W = \lambda(r,s)\delta u_z(r,s)$, the second term in Eq. (\ref{middSurf2}) yields the nodal moment corrections corresponding to the distributed load described by Eq. (\ref{Prs}) of the form
	\begin{align}
	\bm{m}_i=&\int_{\e} \lambda(r,s)\left[\sum_{k=5}^8a_k(r,s)\frac{L_k}{8} 
	\left[
	\begin{array}{c}
	l_{kX}\left(\delta_{iI(k)}-\delta_{iJ(k)}\right)\\
	l_{kY}\left(\delta_{iI(k)}-\delta_{iJ(k)}\right)\\
	l_{kZ}\left(\delta_{iI(k)}-\delta_{iJ(k)}\right)\\
	\end{array}
	\right]\right] \,\dA= \nonumber\\
	&\sum_{q=1}^4  \lambda(r_q,s_q)\, w_q \det J(r_q,s_q) \left[\sum_{k=5}^8a_k(r_q,s_q)\frac{L_k}{8} 
	\left[
	\begin{array}{c}
	l_{kX}\left(\delta_{iI(k)}-\delta_{iJ(k)}\right)\\
	l_{kY}\left(\delta_{iI(k)}-\delta_{iJ(k)}\right)\\
	l_{kZ}\left(\delta_{iI(k)}-\delta_{iJ(k)}\right)\\
	\end{array}
	\right]\right], ~i=1,2,3,4,
	\end{align}
	where $w_q,~q=1,2,3,4$ are the weights and $r_q,\,s_q,~q=1,2,3,4$ are the quadrature points of the Gauss quadrature $2\times 2$. In the component form we finally get
	\begin{align}
	\left[
	\begin{array}{c}
	m_{1X} \\
	m_{2X} \\
	m_{3X} \\
	m_{4X} \\
	m_{1Y} \\
	m_{2Y} \\
	m_{3Y} \\
	m_{4Y} \\
	m_{1Z} \\
	m_{2Z} \\
	m_{3Z} \\
	m_{4Z} \\
	\end{array}\right]=\frac{1}{8}\sum_{q=1}^4  \lambda(r_q,s_q)\, w_q \det J(r_q,s_q)
	\left[
	\begin{array}{c}
	a_5(r_q,s_q)L_5l_{5X} - a_8(r_q,s_q)L_8l_{8X}\\
	a_6(r_q,s_q)L_6l_{6X} - a_5(r_q,s_q)L_5l_{5X}\\
	a_7(r_q,s_q)L_7l_{7X} - a_6(r_q,s_q)L_6l_{6X}\\
	a_8(r_q,s_q)L_8l_{8X} - a_7(r_q,s_q)L_7l_{7X}\\
	a_5(r_q,s_q)L_5l_{5Y} - a_8(r_q,s_q)L_8l_{8Y}\\
	a_6(r_q,s_q)L_6l_{6Y} - a_5(r_q,s_q)L_5l_{5Y}\\
	a_7(r_q,s_q)L_7l_{7Y} - a_6(r_q,s_q)L_6l_{6Y}\\
	a_8(r_q,s_q)L_8l_{8Y} - a_7(r_q,s_q)L_7l_{7Y}\\
	a_5(r_q,s_q)L_5l_{5Z} - a_8(r_q,s_q)L_8l_{8Z}\\
	a_6(r_q,s_q)L_6l_{6Z} - a_5(r_q,s_q)L_5l_{5Z}\\
	a_7(r_q,s_q)L_7l_{7Z} - a_6(r_q,s_q)L_6l_{6Z}\\
	a_8(r_q,s_q)L_8l_{8Z} - a_7(r_q,s_q)L_7l_{7Z}\\
	\end{array}
	\right].
	\end{align}
	These nodal moment corrections are used in the shell elements DKMQ24$_1$+ and DKMQ24$_2$+.
	\FloatBarrier
	\section{Dependence on penalty constants}
	In this Section, the dependence of the relative error at the selected benchmark problems on the penalty constant $c$ (in case of the DKMQ24 shell element), on the penalty constant $c_1$ (in case of the DKMQ24$_1$+ shell element) and the penalty constant $c_2$ (in case of the DKMQ24$_2$+ shell element) is investigated. In general, the benchmark problems can be divided in three categories: The first category contains problems whose error is higher for higher values of the penalty constant. The second category, which is in our case presented only by Raasch's hook, is opposite -- a sufficient value of the penalty constant is needed in order to get a precise solution. The last and the largest category of benchmark problems exhibits no or negligible dependence on the penalty constant. We have selected some representative benchmark problems, so that all three categories are included in our comparison. 
	
	Let us start with the DKMQ24 shell element in our implementation, see Fig. \ref{fig:DKMQ0}. Raasch's hook allows only a small range of the penalty constant, approx. $c \in [10^{-2},10^{-1}]$, in which the relative error is acceptible. On the other hand, the hemispherical shell benchmark with $16\times 16$ mesh belongs to the first category -- the penalty constant needs to satisfy $c < 10^{-4}$ in order to get acceptible accuracy at this benchmark. In conclusion, there is no satisfactory choice of the penalty constant $c$ for this element. Let us now focus on the DKMQ24$_1$+ shell element, see Fig. \ref{fig:DKMQ1}. Raasch's hook benchmark problem needs the penalty constant $c_1$ to be higher then approx. $10^{-1}$ to get a reasonable error. The twisted beam benchmark problem (with setting $h=0.0032\,\text{m}, F_Y>0$) provides the lower error for the lower values of the penalty constant $c_1$. Therefore, the value $c_1=0.15$, used in our calculation is near the optimal value. Let us now focus on the shell element DKMQ24$_2$+, see Fig. \ref{fig:DKMQ2}. The situation here is qualitatively similar to the situation of the shell element DKMQ24$_1$+, only the relative errors are lower, see Fig. \ref{fig:DKMQ2}. Raasch's hook benchmark problem needs the penalty constant $c_2$ to be higher then approx. $0.05$. The twisted beam benchmark problem (with setting $h=0.0032\,\text{m}, F_Y>0$)  needs the penalty constant $c_2$ not to be higher then approx. $0.1$. Therefore, the value $c_2=0.1$ used in our calculation is near the optimal value.
	
	By comparison of the DKMQ24$_1$+ and DKMQ24$_2$+ shell elements at all benchmark problems given in Section 5, we see that they provide almost an identical error except for two cases: the twisted beam benchmark problem ($h=0.0032$) and the Raasch's hook benchmark problem, at both of which the DKMQ24$_2$+ element is superior. Therefore these two benchmarks were selected into our comparison in this Section. Moreover, we therefore recommend the DKMQ24$_2$+ shell element over the DKMQ24$_1$+ shell element.

	\begin{figure}
		\centering
		\begin{tikzpicture} 
		\begin{semilogxaxis}[ 
		xlabel= $c~\text{[-]}$, 
		ylabel= \text{relative error [-]}, 
		x tick label style={
			/pgf/number format/1000 sep={}  
		},
		y tick label style={
			/pgf/number format/fixed,       
			/pgf/number format/precision=5  
		},
		scaled y ticks=false, 
		xmin=0.0000001,
		xmax=10,
		ymin=0,
		ymax=0.25,
		legend style={at={(1.02,0.02)},anchor=south west}	
		] 
		\addplot[line width=1.5pt, color=black,mark=none] coordinates { 
			(0.001,  80.0E-2) 
			(0.001, 0.001E-2)
		};
		\addplot[line width=1.5pt, color=blue, mark=triangle*] coordinates {
			(0.0000001, 7.87E-2)
			(0.000001, 7.87E-2)
			(0.00001, 7.87E-2)
			(0.0001, 7.87E-2)
			(0.001, 7.87E-2)
			(0.01, 7.87E-2) 
			(0.1, 7.87E-2)
			(1.0, 7.87E-2)
			(10.0, 7.87E-2)
		};
		\addplot[line width=1.5pt, color=col1, mark=x] coordinates {  
			(0.0000001, 0.69E-2)
			(0.000001, 1.32E-2)
			(0.00001, 1.44E-2)
			(0.0001, 1.45E-2)
			(0.001, 1.45E-2)
			(0.01, 1.45E-2)
			(0.1, 1.45E-2)
			(1.0, 1.45E-2)
			(10.0, 1.45E-2)
		};
		\addplot[line width=1.5pt, color=black, mark=square, dashed] coordinates {  
			(0.0000001, 0.27E-2)
			(0.000001, 0.29E-2)
			(0.00001, 0.29E-2) 
			(0.0001, 0.29E-2)
			(0.001, 0.29E-2)
			(0.01, 0.29E-2)
			(0.1, 0.29E-2)
			(1.0, 0.29E-2)
			(10.0, 0.29E-2)
		};
		\addplot[line width=1.5pt, color=green, mark=square*, mark size=4pt] coordinates {  
			(0.0000001, 155.85E-2)
			(0.000001, 153.02E-2)
			(0.00001, 133.84E-2)
			(0.0001, 69.17E-2)
			(0.001, 17.70E-2)
			(0.01, 4.17E-2)
			(0.1, 4.33E-2)
			(1.0, 40.38E-2)
			(10.0, 81.91E-2)
		};
		\addplot[line width=1.5pt, color=purple, mark=+, mark size=5pt] coordinates {    
			(0.0000001, 0.12E-2)
			(0.000001, 0.83E-2)
			(0.00001, 3.25E-2)
			(0.0001, 5.08E-2)
			(0.001, 9.33E-2)
			(0.01, 36.15E-2) 
			(0.1, 82.24E-2)
			(1.0, 96.79E-2)
			(10.0, 99.01E-2)
		};
		\addplot[line width=1.5pt, color=yellow, mark=x, mark size=5pt] coordinates { 
			(0.0000001, 7.45E-2)
			(0.000001, 7.45E-2)
			(0.00001, 7.45E-2)
			(0.0001, 7.48E-2)
			(0.001, 7.52E-2)
			(0.01, 7.62E-2) 
			(0.1, 8.46E-2)
			(1.0, 15.45E-2)
			(10.0, 41.84E-2)
		}; 
		\legend{$c=0.001$, Cook's membrane $8\times8$, Twisted beam $h=0.0032\,\text{m};F_Y>0; 2\times12$, Twisted beam $h=0.0032\,\text{m};F_Y>0; 8\times48$, Raasch's hook $136\times 20$, Hemispherical shell $16\times16$, Pinched cylinder $h=0.03\,\text{m}$; dist.mesh; $16\times16$ } 
		\end{semilogxaxis} 
		\end{tikzpicture} 
		\caption{Relative error of the deflection for various benchmark problems calculated by the DKMQ24 element as a function of the penalty parameter $c$}
		\label{fig:DKMQ0}
	\end{figure}
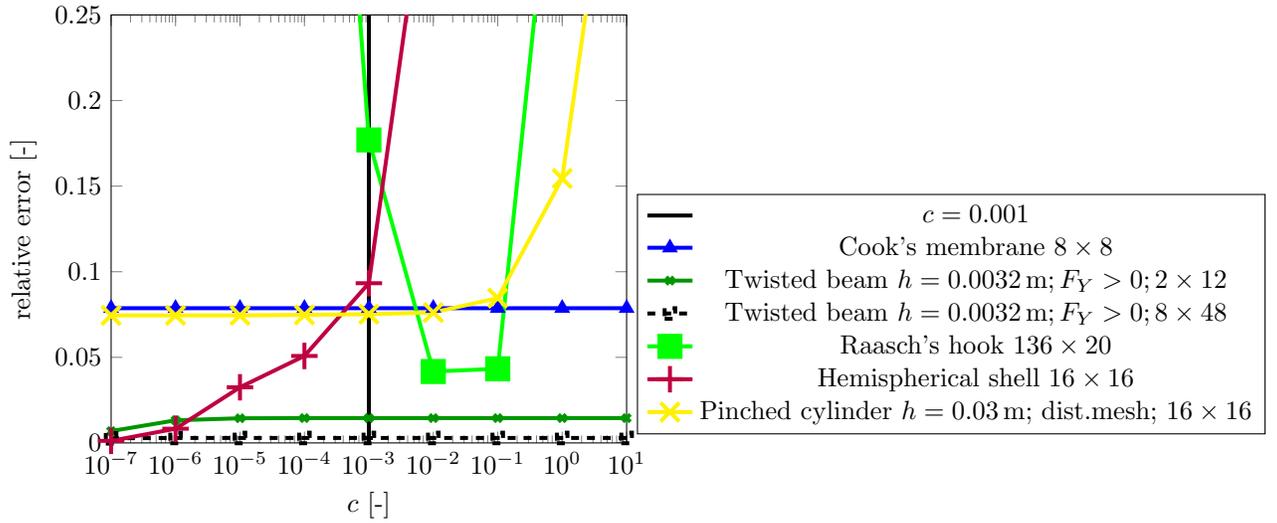
	
	\begin{figure}
		\centering
		\begin{tikzpicture} 
		\begin{semilogxaxis}[ 
		xlabel= $c_1~\text{[-]}$, 
		ylabel= \text{relative error [-]}, 
		x tick label style={
			/pgf/number format/1000 sep={}  
		},
		y tick label style={
			/pgf/number format/fixed,       
			/pgf/number format/precision=5  
		},
		scaled y ticks=false, 
		xmin=0.001,
		xmax = 100000,
		ymin=0,
		ymax=0.25,
		legend style={at={(1.02,0.02)},anchor=south west}	
		] 
		\addplot[line width=1.5pt, color=black,mark=none] coordinates { 
			(0.15,  80.0E-2) 
			(0.15, 0.001E-2)
		};
		\addplot[line width=1.5pt, color=blue, mark=triangle*] coordinates {
			(0.001, 0.77E-2)
			(0.01, 0.78E-2) 
			(0.1, 0.83E-2) 
			(1, 0.93E-2) 
			(10, 1.07E-2)
			(100,  1.16E-2) 
			(1000, 1.19E-2)
			(10000, 1.19E-2)
			(100000,  1.19E-2)
		};
		\addplot[line width=1.5pt, color=col1, mark=x] coordinates {  
			(0.001, 4.13E-2) 
			(0.01, 11.08E-2)
			(0.1, 16.82E-2)
			(1, 18.08E-2)
			(10, 18.27E-2)
			(100,  18.29E-2) 
			(1000, 18.29E-2)
			(10000, 18.29E-2)
			(100000,   18.29E-2) 
		};
		\addplot[line width=1.5pt, color=black, mark=square, dashed] coordinates {  
			(0.001, 0.11E-2)
			(0.01, 0.17E-2)
			(0.1, 0.37E-2)
			(1, 0.82E-2)
			(10, 1.28E-2)
			(100,  1.40E-2) 
			(1000, 1.42E-2)
			(10000, 1.41E-2)
			(100000,   1.34E-2)
		};
		\addplot[line width=1.5pt, color=green, mark=square*, mark size=4pt] coordinates {
			(0.001, 80.26E-2)
			(0.01, 15.96E-2)
			(0.1, 1.56E-2)
			(1, 0.33E-2)
			(10, 0.53E-2)
			(100,  0.55E-2)
			(1000, 0.55E-2)
			(10000, 0.55E-2)
			(100000,   0.55E-2) 
		};
		\addplot[line width=1.5pt, color=purple, mark=+, mark size=5pt] coordinates {    
			(0.001, 0.03E-2)
			(0.01, 0.03E-2)
			(0.1, 0.04E-2)
			(1, 0.07E-2)
			(10, 0.10E-2)
			(100,  0.12E-2) 
			(1000, 0.16E-2)
			(10000, 0.20E-2)
			(100000,   0.21E-2)
		};
		\addplot[line width=1.5pt, color=yellow, mark=x, mark size=5pt] coordinates { 
			(0.001, 7.79E-2)
			(0.01, 7.79E-2)
			(0.1, 7.81E-2)
			(1, 7.86E-2)
			(10, 8.00E-2)
			(100,  8.19E-2) 
			(1000, 8.24E-2)
			(10000, 8.25E-2)
			(100000, 8.25E-2)
		}; 
		\legend{$c_1=0.15$, Cook's membrane $8\times8$, Twisted beam $h=0.0032\,\text{m};F_Y>0; 2\times12$, Twisted beam $h=0.0032\,\text{m};F_Y>0; 8\times48$, Raasch's hook $136\times 20$, Hemispherical shell $16\times16$, Pinched cylinder $h=0.03\,\text{m}$; dist.mesh; $16\times16$ } 
		\end{semilogxaxis} 
		\end{tikzpicture} 
		\caption{Relative error of the deflection for various benchmark problems calculated by the DKMQ24$_1$+ element as a function of the penalty parameter $c_1$}
		\label{fig:DKMQ1}
	\end{figure}
	
	\begin{figure}
		\centering
		\begin{tikzpicture} 
		\begin{semilogxaxis}[ 
		xlabel= $c_2~\text{[-]}$, 
		ylabel= \text{relative error [-]}, 
		x tick label style={
			/pgf/number format/1000 sep={}  
		},
		y tick label style={
			/pgf/number format/fixed,       
			/pgf/number format/precision=5  
		},
		scaled y ticks=false, 
		xmin=0.001,
		xmax = 100000,
		ymin=0,
		ymax=0.25,
		legend style={at={(1.02,0.02)},anchor=south west}	
		] 
		\addplot[line width=1.5pt, color=black,mark=none] coordinates { 
			(0.1,  80.0E-2) 
			(0.1, 0.001E-2)
		};
		\addplot[line width=1.5pt, color=blue, mark=triangle*] coordinates {
			(0.001, 0.77E-2) 
			(0.01, 0.77E-2) 
			(0.1, 0.79E-2) 
			(1.0, 0.85E-2) 
			(10.0, 0.97E-2) 
			(100.0, 1.10E-2) 
			(1000.0, 1.17E-2) 
			(10000.0, 1.19E-2) 
			(100000.0, 1.19E-2) 
		};
		\addplot[line width=1.5pt, color=col1, mark=x] coordinates {  
			(0.001, 0.64E-2) 
			(0.01, 0.97E-2) 
			(0.1, 2.72E-2) 
			(1.0, 8.15E-2) 
			(10.0, 15.42E-2) 
			(100.0, 17.87E-2) 
			(1000.0, 18.23E-2) 
			(10000.0, 18.29E-2) 
			(100000.0, 18.29E-2) 
		};
		\addplot[line width=1.5pt, color=black, mark=square, dashed] coordinates {  
			(0.001, 0.09E-2) 
			(0.01, 0.09E-2) 
			(0.1, 0.11E-2) 
			(1.0, 0.20E-2) 
			(10.0, 0.45E-2) 
			(100.0, 0.95E-2) 
			(1000.0, 1.33E-2) 
			(10000.0, 1.41E-2) 
			(100000.0, 1.42E-2) 
		};
		\addplot[line width=1.5pt, color=green, mark=square*, mark size=4pt] coordinates {  
			(0.001, 55.20E-2) 
			(0.01, 7.48E-2) 
			(0.1, 0.33E-2) 
			(1.0, 0.46E-2) 
			(10.0, 0.54E-2) 
			(100.0, 0.55E-2) 
			(1000.0, 0.55E-2)  
			(10000.0, 0.55E-2)
			(100000.0, 0.55E-2) 
		};
		\addplot[line width=1.5pt, color=purple, mark=+, mark size=5pt] coordinates {    
			(0.001, 0.03E-2) 
			(0.01, 0.03E-2) 
			(0.1, 0.03E-2) 
			(1.0, 0.04E-2) 
			(10.0, 0.06E-2) 
			(100.0, 0.09E-2) 
			(1000.0, 0.12E-2) 
			(10000.0, 0.14E-2) 
			(100000.0, 0.19E-2) 
		};
		\addplot[line width=1.5pt, color=yellow, mark=x, mark size=5pt] coordinates { 
			(0.001, 7.79E-2) 
			(0.01, 7.79E-2) 
			(0.1, 7.79E-2) 
			(1.0, 7.81E-2) 
			(10.0, 7.87E-2) 
			(100.0, 8.01E-2) 
			(1000.0,  8.19E-2) 
			(10000.0, 8.25E-2) 
			(100000.0, 8.25E-2) 
		}; 
		\legend{$c_2=0.1$, Cook's membrane $8\times8$, Twisted beam $h=0.0032\,\text{m};F_Y>0; 2\times12$, Twisted beam $h=0.0032\,\text{m};F_Y>0; 8\times48$, Raasch's hook $136\times 20$, Hemispherical shell $16\times16$, Pinched cylinder $h=0.03\,\text{m}$; dist.mesh; $16\times16$ } 
		\end{semilogxaxis} 
		\end{tikzpicture} 
		\caption{Relative error the deflection for various benchmark problems calculated by the DKMQ24$_2$+ element as a function of the penalty parameter $c_2$}
		\label{fig:DKMQ2}
	\end{figure}
	
	\FloatBarrier
	\section{Convergence tests}

	\label{Benchmarks}
	The behaviour of the DKMQ24, DKMQ24$_1$+ and DKMQ24$_2$+ shell elements has been tested on nine well-established benchmark problems, including one pure membrane benchmark problem suggested by Cook \cite{Cook}, six shell and one plate benchmark problems suggested by MacNeal and Harder \cite{MacNealHarder} and Belytschko \cite{Belytschko}, which were reused by Katili \cite{Katili2015}, and one challenging benchmark proposed by Raasch \cite{Raasch}. We compare the results of the DKMQ24$_1$+ and DKMQ24$_2$+ shell elements with the recent state-of-the-art quadrangle shell/solid-shell elements, which are all summarized in Table 4. In all cases, geometric and material linearity is considered and there is no initial deformation. 	
	
	\begin{center}
		\begin{table}
			\centering
			\begin{tabular}{ccc}
				\hline
				Notation &  Source &  Description \\
				\hline
				\hline
				MITC4+                     & from ref.\cite{MITC4Plus3}, (2017) & - \\
				\!US-ATFHS8                & from ref.\cite{Huang}, (2018) & 8-nodal solid-shell element, without rotational dofs, 24 dofs \\
				DKMQ24 \!\cite{Katili2015} & from ref.\cite{Katili2015}, (2015) & - \\
				DKMQ24                     & our implementation & DKMQ24 element ($c=0.001$ as in \cite{Katili2015})\\
				DKMQ24($c=10^{-7}$)        & our implementation & DKMQ24 element ($c=10^{-7}$)\\
				DKMQ24$_1$+                & our implementation & improved DKMQ24 element, which uses Eq. (\ref{phiz2}) ($c_1=0.15$) \\
				DKMQ24$_2$+                & our implementation & improved DKMQ24 element, which uses Eq. (\ref{penalty2}) ($c_2=0.1$) \\
				\hline
			\end{tabular} \nonumber
			\label{tabEl}
			\caption{Tested elements}
		\end{table}
	\end{center}
	\subsection{Cook's membrane}
	\label{b1}
	\begin{figure}[H]
		\centering
		\begin{overpic}[width=0.33\textwidth]{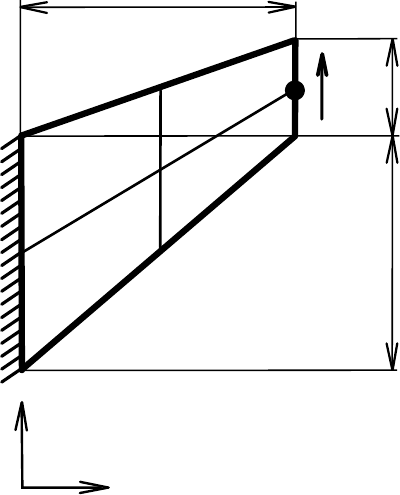} 
			\put(29,102.7){$L_1$}
			\put(81.5,49){$L_2$}
			\put(81.5,81){$L_3$}
			\put(68,81){$F_Z$}
			\put(53,83){$A$}
			\put(23,-2){$X$}
			\put(-2,14){$Z$}
			\put(34,42){thickness $h$}
			\put(100,92){$L_1=48\,\text{m}$}
			\put(100,82){$L_2=44\,\text{m}$}
			\put(100,72){$L_3=16\,\text{m}$}
			\put(100,62){$h=1\,\text{m}$}
			\put(100,52){$E=1\,\text{kPa}$}
			\put(100,42){$\nu=\frac{1}{3}$}
			\put(100,32){$F_Z=1\,\text{N}$}
		\end{overpic}
		\caption{Cook's membrane with mesh $2\times 2$}
		\label{Cook}
	\end{figure}
	\FloatBarrier
	
	The standard Cook's membrane problem tests an element behaviour for the case of in-plane shear loading. The problem also tests the effect of mesh distortion. The membrane is fully fixed at $X=0$ and loaded by a distributed force at $X=L_1$ (see Figure \ref{Cook}). The $Z$-component of the displacement is evaluated at the test point A. The analytical solution for this example is not available, therefore the reference solution is obtained numerically using the DKMQ24$_2$+ element on the refined mesh $192\!\times\!128$. The numerical results are summarized in Table 3. We observe a massive reduction of the relative error in case of the DKMQ24$_1$+ and DKMQ24$_2$+ shell elements in comparison to the MITC4+ and DKMQ24 shell elements. 
	
	\begin{center}
		\begin{table}
			\centering
			\begin{tabularx}{160mm}{ c *{8}{Y} }
				\hline
				Mesh  
				& \multicolumn{2}{c}{MITC4+} 
				& \multicolumn{2}{c}{DKMQ24}  
				& \multicolumn{2}{c}{DKMQ24$_1$+} 
				& \multicolumn{2}{c}{DKMQ24$_2$+}  \rule[-2mm]{0mm}{6mm}  \\
				& $u_Z(A)$ & err. & $u_Z(A)$ & err. & $u_Z(A)$ & err. & $u_Z(A)$ & err.    \\
				& [mm] &  [\%]    & [mm] &  [\%]    & [mm] &  [\%]     & [mm] &  [\%]      \\
				\hline
				\hline
				$2\!\times\!2$   & $11.729$ & $51.06$    & $11.845$ & $50.57$ &  $22.353$ & $6.72$ &  $22.396$ & $6.55$ \\
				$4\!\times\!4$   & $18.266$ & $23.78$    & $18.299$ & $23.64$ &  $23.397$ & $2.37$ &  $23.415$ & $2.30$ \\
				$8\!\times\!8$   & $22.075$ &  $7.89$    & $22.079$ &  $7.87$ &  $23.764$ & $0.84$ &  $23.776$ & $0.79$ \\
				$16\!\times\!16$ & $23.430$ &  $2.23$    & $23.430$ &  $2.23$ &  $23.892$ & $0.30$ &  $23.897$ & $0.29$ \\
				$32\!\times\!32$ & $23.818$ &  $0.62$    & $23.818$ &  $0.61$ &  $23.939$ & $0.11$ &  $23.940$ & $0.11$ \\
				\hline
				\multicolumn{9}{c}{Reference solution: $u_Z(A)=23.965\,\text{mm}$ (DKMQ24$_2$+ element on mesh $192\!\times\!128$)} \\ 
				\hline
			\end{tabularx} \nonumber
			\caption{Deflection $u_Z(A)$ for Cook's membrane problem}
		\end{table}
	\end{center}
	
	\begin{figure}
		\centering
		\begin{tikzpicture} 
		\begin{loglogaxis}[ 
		xlabel=\text{number of elements}, 
		ylabel= $\frac{u_Z(A)}{u_{Z,\rf}(A)}$,  
		x tick label style={
			/pgf/number format/1000 sep={}  
		},
		y tick label style={
			/pgf/number format/fixed,       
			/pgf/number format/precision=5  
		},
		scaled y ticks=false, 
		xmin= 4,
		xmax = 1074,
		legend style={at={(0.02,0.02)},anchor=south west}	
		] 
		\addplot[line width=1.5pt, color=purple, mark=+, mark size=5pt] coordinates {
			(4,    51.06E-2) 
			(16,   23.78E-2) 
			(64,   7.89E-2) 
			(256,  2.23E-2) 
			(1024, 0.62E-2) 
		};
		\addplot[line width=1.5pt, color=blue, mark=triangle*, mark size=4pt, dotted] coordinates {
			(4,     50.57E-2) 
			(16,    23.64E-2) 
			(64,     7.87E-2) 
			(256,    2.23E-2) 
			(1024,   0.61E-2)
		};
		\addplot[line width=1.5pt, color=col1, mark=x, mark size=5pt, dashed] coordinates {  
			(4,      6.72E-2) 
			(16,     2.37E-2) 
			(64,     0.84E-2) 
			(256,    0.30E-2) 
			(1024,   0.105E-2) 
		};
		\addplot[line width=1.5pt, color=green, mark=square*, mark size=4pt] coordinates {  
			(4,     6.55E-2) 
			(16,    2.30E-2) 
			(64,    0.79E-2) 
			(256,   0.29E-2) 
			(1024,  0.105E-2)
		};
		\legend{MITC4+, DKMQ24, DKMQ24$_1$+, DKMQ24$_2$+} 
		\end{loglogaxis} 
		\end{tikzpicture} 
		\caption{Relative error of deflection $u_Z(A)$ for Cook's membrane problem}
		\label{fig:BENCH1}
	\end{figure}
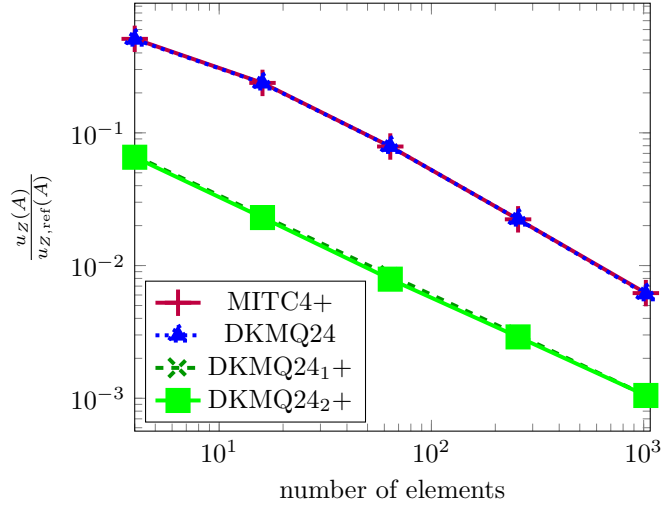

	\subsection{Hemispherical shell}
	\label{b2}
	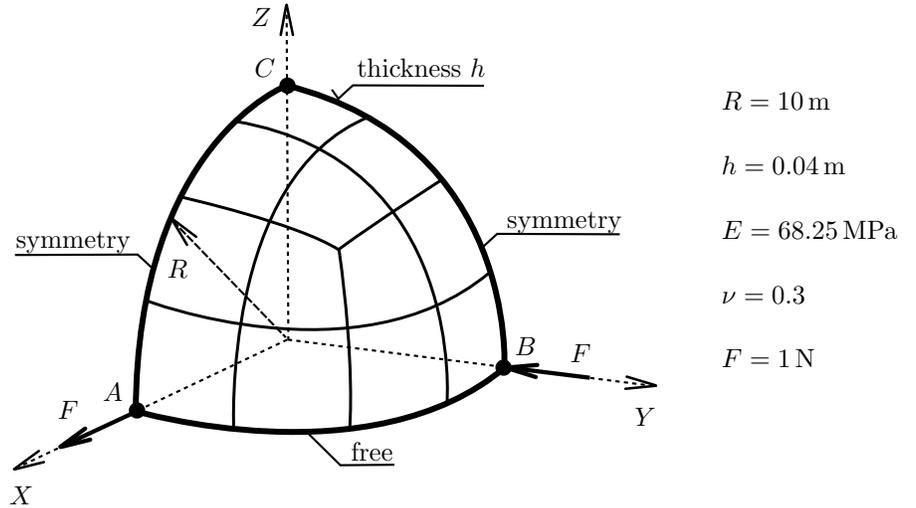
\begin{figure}[H]
		\centering
		\begin{overpic}[width=0.53\textwidth]{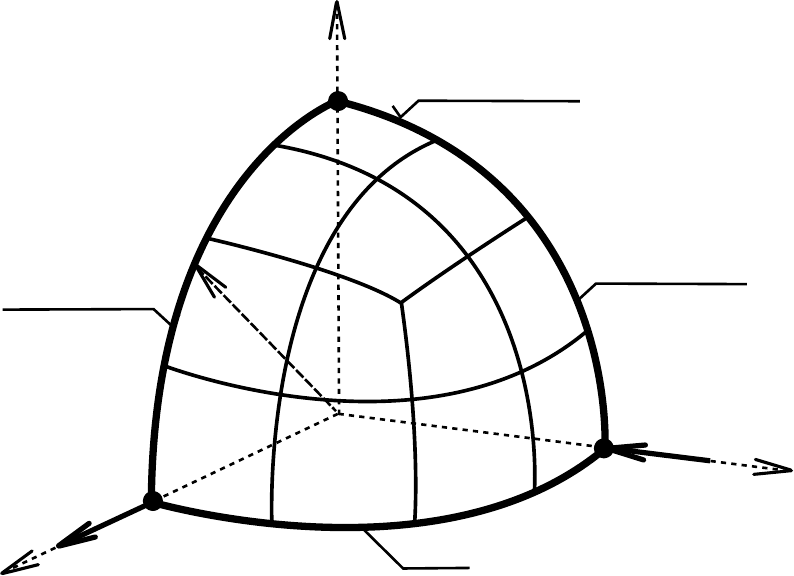} 
			\put(7,8){$F$}
			\put(86.5,16.8){$F$}
			\put(14,10.75){$A$}
			\put(77.9,18){$B$}
			\put(37.5,61.2){$C$}
			\put(-0.5,-5.3){$X$}
			\put(96.5,7){$Y$}
			\put(37,69){$Z$}
			\put(52.5,1.6){free}
			\put(0.5,34.6){symmetry}
			\put(24,30){$R$}
			\put(76.7,37.75){symmetry}
			\put(53.4,61.1){thickness $h$}
			\put(110,56){$R=10\,\text{m}$}
			\put(110,46){$h=0.04\,\text{m}$}
			\put(110,36){$E=68.25\,\text{MPa}$}
			\put(110,26){$\nu=0.3$}
			\put(110,16){$F=1\,\text{N}$}
		\end{overpic}
		\caption{Hemispherical shell with mesh $4\times 4$}
		\label{HemisphericalShell}
	\end{figure}
	\FloatBarrier
	
	The hemispherical shell problem, consisting of a hemisphere loaded by four opposite single forces, is depicted in Figure \ref{HemisphericalShell}. Due to symmetries, only one quarter of this hemisphere needs to be considered. The following symmetry conditions are applied: $u_Y=\varphi_X=\varphi_Z=0$ on the side AC and $u_X=\varphi_Y=\varphi_Z=0$ on the side BC. Additionally, the following boundary condition is applied at point C: $u_Z=0$. The analytical solution to this problem is given by \cite{MacNealHarder} as $u_X(A)=-u_Y(B)=92.400\,\text{mm}$. 
	
	Note that this example is sensitive to the mesh quality and it is also the only benchmark in our set for which the mesh is not exactly specified by its depiction in Figure \ref{HemisphericalShell}. In our case we have used the Laplacian smoothing of internal nodes \cite{Laplacian} until we have converged to a unique mesh for the calculation. 
	
	We report that our implementation of the DKMQ24 shell element gives significantly higher error then the one reported in \cite{Katili2015}. A similar error was obtained for the reduced penalty constant $c=10^{-7}$. 
	In comparison to the DKMQ24 element from \cite{Katili2015} and MITC4+ element, our elements DKMQ24$_1$+ and DKMQ24$_2$+ provide higher errors on coarser meshes and lower errors on finer meshes.
	

	\begin{center}
		\begin{table}
			\centering
			\begin{tabularx}{160mm}{ c *{12}{Y} }
				\hline
				Mesh  
				& \multicolumn{2}{c}{MITC4+}  
				& \multicolumn{2}{c}{DKMQ24 \cite{Katili2015}}  
				& \multicolumn{2}{c}{DKMQ24}  
				& \multicolumn{2}{c}{DKMQ24}   
				& \multicolumn{2}{c}{DKMQ24$_1$+} 
				& \multicolumn{2}{c}{DKMQ24$_2$+}  \rule[-2mm]{0mm}{6mm}\\
				& \multicolumn{2}{c}{}  
				& \multicolumn{2}{c}{}  
				& \multicolumn{2}{c}{}  
				& \multicolumn{2}{c}{($c=10^{-7}$)}   
				& \multicolumn{2}{c}{} 
				& \multicolumn{2}{c}{}  \\			
				& $u_X(A)$ & err. & $u_X(A)$ & err. & $u_X(A)$ & err. & $u_X(A)$ & err. & $u_X(A)$ & err.& $u_X(A)$ & err.\\
				& [mm] &  [\%]    & [mm] &  [\%]    & [mm] &  [\%]     & [mm] &  [\%]    & [mm] &  [\%]  & [mm] &  [\%] \\
				\hline
				\hline
				$4\!\times\!4$ &  $94.802$ & $2.60$ &$78.721$&$14.80$&$49.911$ & $45.98$ & $85.305$ &  $7.68$ &  $40.625$ & $56.03$ & $48.902$ & $47.08$\\
				$8\!\times\!8$ &  $92.585$ & $0.20$ &$93.624$&$1.32$& $74.940$ & $18.90$ & $93.963$ &  $1.69$ &  $88.929$ & $3.76$ &  $89.224$ & $3.44$\\
				$16\!\times\!16$& $92.030$ & $0.40$ &$92.246$&$0.17$& $83.775$ & $9.33$  & $92.289$ &  $0.12$ &  $92.356$ & $0.05$ &  $92.372$ & $0.03$\\
				$32\!\times\!32$& $92.234$ & $0.18$ &-&-&             $86.654$ & $6.22$  & $92.146$ &  $0.28$ &  $92.405$ & $0.01$ &  $92.406$ & $0.01$\\
				\hline
				\multicolumn{13}{c}{Reference solution: $u_X(A)=-u_Y(B)=92.400\,\text{mm}$ (\cite{MacNealHarder})} \\ 
				\hline
			\end{tabularx} \nonumber
			\caption{Deflection $u_X(A)$ for the hemispherical shell problem}
		\end{table}
	\end{center}

	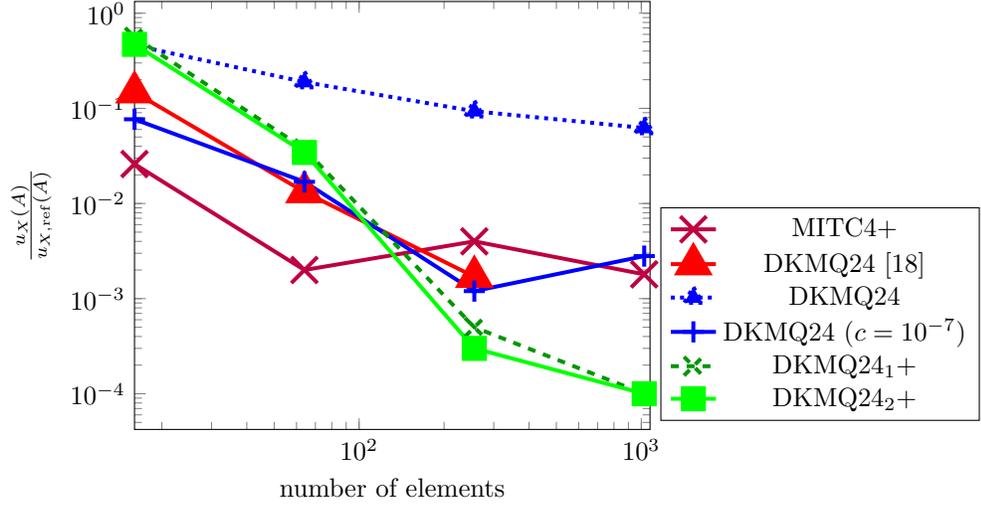
\begin{figure}
		\centering
		\begin{tikzpicture} 
		\begin{loglogaxis}[ 
		xlabel=\text{number of elements},
		ylabel= $\frac{u_X(A)}{u_{X,\rf}(A)}$, 
		x tick label style={
			/pgf/number format/1000 sep={}  
		},
		y tick label style={
			/pgf/number format/fixed,       
			/pgf/number format/precision=5  
		},
		scaled y ticks=false, 
		xmin= 16,
		xmax = 1074,
		legend style={at={(1.02,0.02)},anchor=south west}	
		] 
		\addplot[line width=1.5pt, color=purple, mark=x, mark size=7pt] coordinates {
			(16,  2.60E-2) 
			(64,  0.2E-2) 
			(256, 0.4E-2) 
			(1024,0.18E-2)
		}; 
		\addplot[line width=1.5pt, color=red, mark=triangle*, mark size=7pt] coordinates { 
			(16,  14.80E-2) 
			(64,  1.32E-2) 
			(256, 0.17E-2) 
		}; 
		\addplot[line width=1.5pt, color=blue, mark=triangle*, mark size=4pt, dotted] coordinates { 
			(16,  45.98E-2) 
			(64,  18.90E-2) 
			(256,  9.33E-2) 
			(1024, 6.22E-2)
		};  
		\addplot[line width=1.5pt, color=blue, mark=+, mark size=4pt] coordinates { 
			(16,  7.68E-2) 
			(64,  1.69E-2) 
			(256, 0.12E-2) 
			(1024,0.28E-2)
		}; 
		\addplot[line width=1.5pt, color=col1, mark=x, mark size=5pt, dashed] coordinates { 
			(16,  56.03E-2) 
			(64,  3.76E-2) 
			(256, 0.05E-2) 
			(1024,0.01E-2)
		}; 
		\addplot[line width=1.5pt, color=green, mark=square*, mark size=4pt] coordinates {
			(16,  47.08E-2) 
			(64,  3.44E-2) 
			(256, 0.03E-2) 
			(1024,0.01E-2)
		}; 
		\legend{MITC4+, DKMQ24 \cite{Katili2015}, DKMQ24, DKMQ24 ($c=10^{-7}$), DKMQ24$_1$+, DKMQ24$_2$+}
		\end{loglogaxis} 
		\end{tikzpicture} 
		\caption{Relative error of deflection $u_X(A)$ for the hemispherical shell problem}
		\label{fig:BENCH2}
	\end{figure}

	\subsection{Scordelis-Lo roof}
	\label{b3}
	\begin{figure}[H]
		\centering
		\begin{overpic}[width=0.6\textwidth]{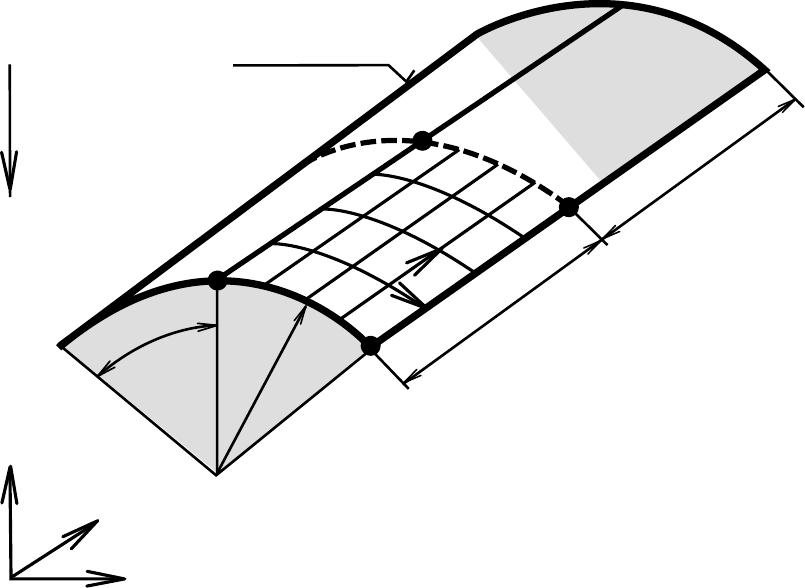} 
			\put(4,57){$g$}
			\put(31,30){$R$}
			\put(64,30){$\frac{L}{2}$}
			\put(87,47.25){$\frac{L}{2}$}
			\put(19.5,26){$\varphi$}
			\put(14,-4.3){$X$}
			\put(10,10){$Y$}
			\put(-4,13.5){$Z$}
			\put(44.4,24.75){$A$}
			\put(69.25,49.5){$B$}
			\put(51.6,58.25){$C$}
			\put(25.5,40.2){$D$}
			\put(47.4,33.25){$x$}
			\put(48.4,40.5){$y$}
			\put(30,66.25){$\text{thickness}~h$}		    
			\put(95,37){$L=6\,\text{m}$}
			\put(95,31){$R=3\,\text{m}$}
			\put(95,25){$h=0.03\,\text{m}$}
			\put(95,19){$\varphi=40^{\circ}$}
			\put(95,13){$E=30000\,\text{MPa}$}
			\put(95,7){$\nu=0$}	
			\put(95,1){$h\rho=625\,\text{kg}\,\text{m}^{-2}$} 	
		\end{overpic}
		\caption{Scordelis-Lo roof with mesh $4\times4$}
		\label{ScordelisLo}
	\end{figure}
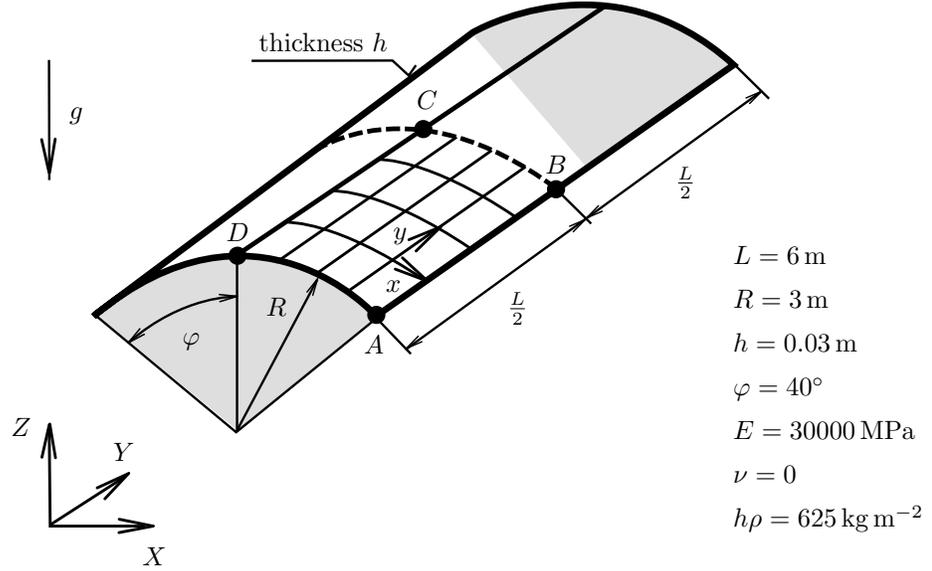
	\FloatBarrier
	
	The Scordelis-Lo roof is a cylindrical surface with radius $R$, length $L$ and thickness $h$, reinforced at positions $Y=0$ and $Y=L$ by two rigid diaphragms (see Figure \ref{ScordelisLo}). Because of symmetry only one quarter is calculated. The shell is loaded by self-weight, where gravitational acceleration is considered to be $g=10\,\text{m}\, \text{s}^{-2}$, which yields the total force acting on the roof $f_Z=-h\rho gA=-12500\pi\,\text{N}$. The following boundary and symmetric conditions are applied: $u_X=u_Z=\varphi_Y=0$ on the side AD, $u_X=\varphi_Y=\varphi_Z=0$ on the side CD and $u_Y=\varphi_X=\varphi_Z=0$ on the side CB. The analytical solution is given in \cite{Scordelis}. We note that two versions of this benchmark, having the same geometry setting differing only in the parameter setting, are given in the literature. The one version is taken from \cite{Scordelis} (described in Fig. \ref{ScordelisLo} and used in Tables 7-10 for comparison with the DKMQ24 \cite{Katili2015} element), the second version taken from \cite{MacNealHarder} (used in Table 11 for comparison with the MITC4+ and \US  elements).
	
	With this benchmark problem we test not only nodal displacements, but also values of the nodal stress-resultants. We point out that nodal stress-resultants are bilinearly extrapolated from the Gauss quadrature points to nodes.
	
	If compared to the DKMQ24 element from \cite{Katili2015}, our modified elements DKMQ24$_1$+ and DKMQ24$_2$+ exhibit lower errors at all test cases. If compared to the MITC4+ and \US elements, our modified elements DKMQ24$_1$+ and DKMQ24$_2$+ exhibit lower errors on both coarser meshes, but higher errors on the finest mesh $16\times 16$.
	

	\begin{center}
		\begin{table}
			\centering
			\begin{tabularx}{142mm}{ c *{8}{Y} }
				\hline
				Mesh  
				& \multicolumn{2}{c}{DKMQ24 \cite{Katili2015}}  
				& \multicolumn{2}{c}{DKMQ24}  
				& \multicolumn{2}{c}{DKMQ24$_1$+} 
				& \multicolumn{2}{c}{DKMQ24$_2$+}  \rule[-2mm]{0mm}{6mm}  \\
				& $u_Z(B)$   & err. &    $u_Z(B)$   & err. &    $u_Z(B)$   & err.  &    $u_Z(B)$   & err. \\
				& [m] &  [\%]     & [m] &  [\%]    & [m] &  [\%]  & [m] &  [\%] \\
				\hline
				\hline
				$4\!\times\!4$   & $-0.03425$  &  $5.12$ & $-0.034258$ & $5.10$ &  $-0.037036$ & $2.59$ &  $-0.037045$ & $2.62$\\
				$8\!\times\!8$   & $-0.03528$  &  $2.26$ & $-0.035284$ & $2.26$ &  $-0.036228$ & $0.35$ &  $-0.036233$ & $0.37$\\
				$16\!\times\!16$ & $-0.03585$  &  $0.70$ & $-0.035858$ & $0.67$ &  $-0.036102$ & $0.01$ &  $-0.036105$ & $0.01$\\
				$32\!\times\!32$ & -           &    -    & $-0.036070$ & $0.08$ &  $-0.036086$ & $0.04$ &  $-0.036078$ & $0.04$\\
				\hline
				\multicolumn{9}{c}{Reference solution: $u_Z(B)=-0.0361\,\text{m}$ (theory of deep shell, \cite{Scordelis})} \\ 
				\hline
			\end{tabularx} \nonumber
			\label{tab7}
			\caption{Deflection $u_Z(B)$ for the Scordelis-Lo roof problem, $B\approx[1.928,3,2.298]$}
		\end{table}
	\end{center}
	
	\begin{figure}
		\centering
		\begin{tikzpicture} 
		\begin{loglogaxis}[ 
		xlabel=\text{number of elements}, 
		ylabel= $\frac{u_Z(B)}{u_{Z,\rf}(B)}$, 
		x tick label style={
			/pgf/number format/1000 sep={}  
		},
		y tick label style={
			/pgf/number format/fixed,       
			/pgf/number format/precision=5  
		},
		scaled y ticks=false, 
		xmin= 16,
		xmax = 1074,
		legend style={at={(1.02,0.02)},anchor=south west}
		]     
		\addplot[line width=1.5pt, color=red, mark=triangle*, mark size=7pt] coordinates { 
			(16,   5.12E-2) 
			(64,   2.26E-2) 
			(256,  0.70E-2) 
		};   
		\addplot[line width=1.5pt, color=blue, mark=triangle*, mark size=4pt, dotted] coordinates { 
			(16,   5.10E-2) 
			(64,   2.26E-2) 
			(256,  0.67E-2) 
			(1024, 0.08E-2)
		};
		\addplot[line width=1.5pt, color=col1, mark=x, mark size=5pt, dashed] coordinates { 
			(16,   2.59E-2) 
			(64,   0.35E-2) 
			(256,  0.01E-2) 
			(1024, 0.04E-2)
		};  
		\addplot[line width=1.5pt, color=green, mark=square*, mark size=4pt] coordinates { 
			(16,   2.62E-2) 
			(64,   0.37E-2) 
			(256,  0.01E-2) 
			(1024, 0.04E-2)
		};    
		\legend{DKMQ24 \cite{Katili2015}, DKMQ24, DKMQ24$_1$+ , DKMQ24$_2$+}
		\end{loglogaxis} 
		\end{tikzpicture} 
		\caption{Relative error of deflection $u_Z(B)$ for the Scordelis-Lo roof problem}
		\label{fig:BENCH3a}
	\end{figure}
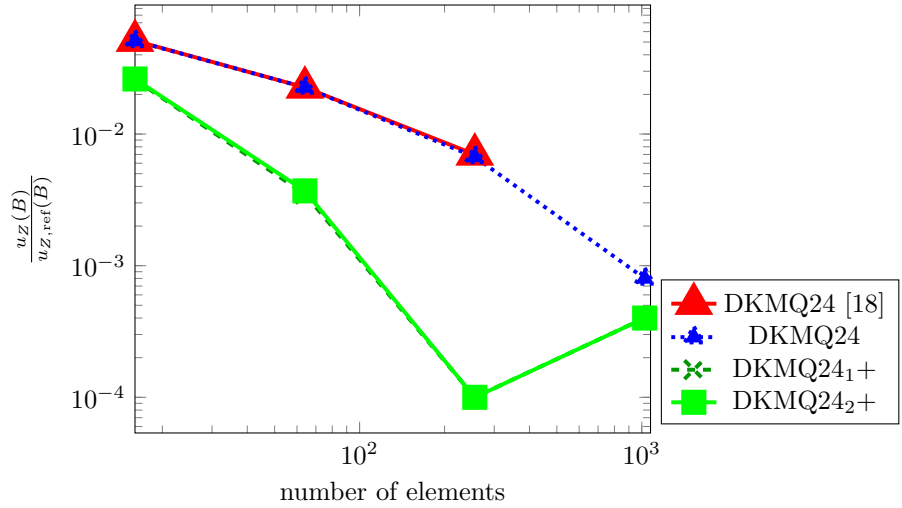

	
	\begin{center}
		\begin{table}
			\centering
			\begin{tabularx}{130mm}{ c *{8}{Y} }
				\hline
				Mesh  
				& \multicolumn{2}{c}{DKMQ24 \cite{Katili2015}}  
				& \multicolumn{2}{c}{DKMQ24}  
				& \multicolumn{2}{c}{DKMQ24$_1$+} 
				& \multicolumn{2}{c}{DKMQ24$_2$+}  \rule[-2mm]{0mm}{6mm}  \\
				& $u_Z(C)$   & err. &    $u_Z(C)$   & err. &    $u_Z(C)$   & err.  &    $u_Z(C)$   & err. \\
				& [m] &  [\%]     & [m] &  [\%]    & [m] &  [\%]  & [m] &  [\%] \\
				\hline
				\hline
				$4\!\times\!4$   & $0.00513$   &  $5.18$ & $0.005130$ & $5.17$ & $0.005599$ & $3.49$ & $0.005600$ & $3.51$ \\
				$8\!\times\!8$   & $0.00529$   &  $2.21$ & $0.005294$ & $2.15$ & $0.005447$ & $0.67$ & $0.005447$ & $0.69$ \\
				$16\!\times\!16$ & $0.00538$   &  $0.55$ & $0.005378$ & $0.58$ & $0.005419$ & $0.17$ & $0.005419$ & $0.17$ \\ 
				$32\!\times\!32$ & -           &    -    & $0.005407$ & $0.05$ & $0.005414$ & $0.07$ & $0.005415$ & $0.08$ \\ 
				\hline
				\multicolumn{9}{c}{Reference solution: $u_Z(C)=0.00541\,\text{m}$ (theory of deep shell, \cite{Scordelis})} \\ 
				\hline
			\end{tabularx} \nonumber
			\label{tab8}
			\caption{Deflection $u_Z(C)$ for the Scordelis-Lo roof problem, $C=[0,3,3]$}
		\end{table}
	\end{center}
	
	\begin{figure}
		\centering
		\begin{tikzpicture} 
		\begin{loglogaxis}[ 
		xlabel=\text{number of elements}, 
		ylabel= $\frac{u_Z(C)}{u_{Z,\rf}(C)}$, 
		x tick label style={
			/pgf/number format/1000 sep={}  
		},
		y tick label style={
			/pgf/number format/fixed,       
			/pgf/number format/precision=5  
		},
		scaled y ticks=false, 
		xmin= 16,
		xmax = 1074,
		legend style={at={(0.02,0.02)},anchor=south west}
		]     
		\addplot[line width=1.5pt, color=red, mark=triangle*, mark size=5pt] coordinates { 
			(16,   5.18E-2) 
			(64,   2.21E-2) 
			(256,  0.55E-2) 
		};   
		\addplot[line width=1.5pt, color=blue, mark=triangle*, mark size=4pt, dotted] coordinates { 
			(16,   5.17E-2) 
			(64,   2.15E-2) 
			(256,  0.58E-2) 
			(1024, 0.05E-2)
		};
		\addplot[line width=1.5pt, color=col1, mark=x, mark size=5pt, dashed] coordinates { 
			(16,   3.49E-2) 
			(64,   0.67E-2) 
			(256,  0.17E-2) 
			(1024, 0.07E-2)
		};  
		\addplot[line width=1.5pt, color=green, mark=square*, mark size=4pt] coordinates { 
			(16,   3.51E-2) 
			(64,   0.69E-2) 
			(256,  0.17E-2) 
			(1024, 0.08E-2)
		};    
		\legend{DKMQ24 \cite{Katili2015}, DKMQ24, DKMQ24$_1$+ , DKMQ24$_2$+}
		\end{loglogaxis} 
		\end{tikzpicture} 
		\caption{Relative error of deflection $u_Z(C)$ for the Scordelis-Lo roof problem}
		\label{fig:BENCH3b}
	\end{figure}
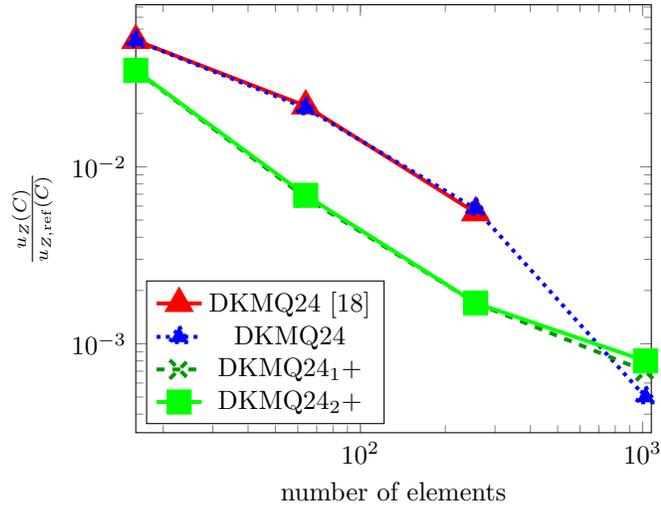
	
	
	\begin{center}
		\begin{table}
			\centering
			\begin{tabularx}{130mm}{ c *{8}{Y} }
				\hline
				Mesh  
				& \multicolumn{2}{c}{DKMQ24 \cite{Katili2015}}  
				& \multicolumn{2}{c}{DKMQ24}  
				& \multicolumn{2}{c}{DKMQ24$_1$+} 
				& \multicolumn{2}{c}{DKMQ24$_2$+}  \rule[-2mm]{0mm}{6mm}  \\
				& $m_x(C)$   & err. &    $m_x(C)$   & err. &    $m_x(C)$   & err.  &    $m_x(C)$   & err. \\
				& [N] &  [\%]     & [N] &  [\%]    & [N] &  [\%]  & [N] &  [\%] \\
				\hline
				\hline
				$4\!\times\!4$   & $1874$    &  $8.85$  & $1740.1$ & $15.36$ & $1930.7$ & $6.09$  & $1931.3$ & $6.07$ \\ 
				$8\!\times\!8$   & $1998$    &  $2.82$  & $1947.8$ & $5.26$  & $2023.4$ & $1.59$  & $2023.7$ & $1.57$ \\ 
				$16\!\times\!16$ & $2041$    &  $0.73$  & $2027.3$ & $1.40$  & $2048.0$ & $0.39$  & $2048.1$ & $0.38$ \\ 
				$32\!\times\!32$ & $2054$    &  $0.10$  & $2050.4$ & $0.27$  & $2054.5$ & $0.07$  & $2054.6$ & $0.07$ \\
				\hline
				\multicolumn{9}{c}{Reference solution: $m_x(C)=2056\,\text{N}$ (theory of shallow shell, \cite{Scordelis})} \\ 
				\hline
			\end{tabularx} \nonumber
			\label{tab9}
			\caption{Internal moment $m_x(C)$ for the Scordelis-Lo roof problem}
		\end{table}
	\end{center}

	\begin{figure}
		\centering
		\begin{tikzpicture} 
		\begin{loglogaxis}[ 
		xlabel=\text{number of elements}, 
		ylabel= $\frac{m_x(C)}{m_{x,\rf}(C)}$,
		x tick label style={
			/pgf/number format/1000 sep={}  
		},
		y tick label style={
			/pgf/number format/fixed,       
			/pgf/number format/precision=5  
		},
		scaled y ticks=false, 
		xmin= 16,
		xmax = 1074,
		legend style={at={(0.02,0.02)},anchor=south west}
		]     
		\addplot[line width=1.5pt, color=red, mark=triangle*, mark size=5pt] coordinates { 
			(16,   8.85E-2) 
			(64,   2.82E-2) 
			(256,  0.73E-2)
			(1024, 0.10E-2) 
		};   
		\addplot[line width=1.5pt, color=blue, mark=triangle*, mark size=4pt, dotted] coordinates { 
			(16,   15.36E-2) 
			(64,   5.26E-2) 
			(256,  1.40E-2) 
			(1024, 0.27E-2)
		};
		\addplot[line width=1.5pt, color=col1, mark=x, mark size=5pt, dashed] coordinates { 
			(16,   6.09E-2) 
			(64,   1.59E-2) 
			(256,  0.39E-2) 
			(1024, 0.07E-2)
		};  
		\addplot[line width=1.5pt, color=green, mark=square*, mark size=4pt] coordinates { 
			(16,   6.07E-2) 
			(64,   1.57E-2) 
			(256,  0.38E-2) 
			(1024, 0.07E-2)
		};    
		\legend{DKMQ24 \cite{Katili2015}, DKMQ24, DKMQ24$_1$+ , DKMQ24$_2$+}
		\end{loglogaxis} 
		\end{tikzpicture} 
		\caption{Relative error of internal moment $m_x(C)$ for the Scordelis-Lo roof problem}
		\label{fig:BENCH3c}
	\end{figure}
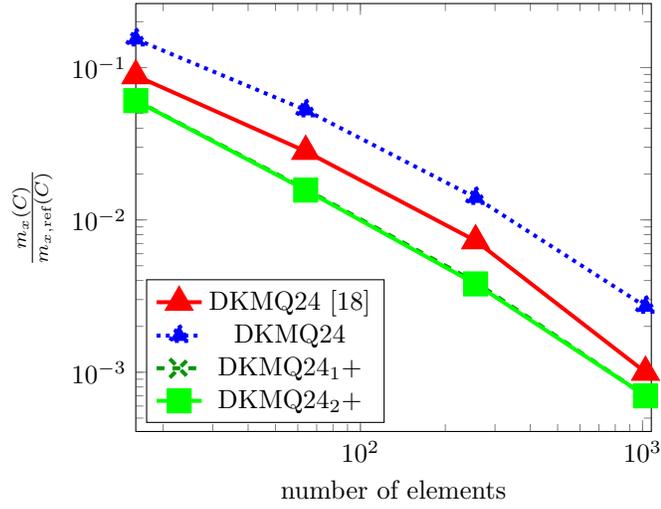
	
	
	\begin{center}
		\begin{table}
			\centering
			\begin{tabularx}{130mm}{ c *{8}{Y} }
				\hline
				Mesh  
				& \multicolumn{2}{c}{DKMQ24 \cite{Katili2015}}  
				& \multicolumn{2}{c}{DKMQ24}  
				& \multicolumn{2}{c}{DKMQ24$_1$+} 
				& \multicolumn{2}{c}{DKMQ24$_2$+}  \rule[-2mm]{0mm}{6mm}  \\
				& $n_y(B)$   & err. &    $n_y(B)$   & err. &    $n_y(B)$   & err.  &    $n_y(B)$   & err. \\
				& [Nm$^{-1}$] &  [\%]     & [Nm$^{-1}$] &  [\%]    & [Nm$^{-1}$] &  [\%]  & [Nm$^{-1}$] &  [\%] \\
				\hline
				\hline
				$4\!\times\!4$   & $498430$    & $22.24$  & $498741$ & $22.19$ & $549772$ & $14.23$ & $550517$ & $14.12$ \\
				$8\!\times\!8$   & $592680$    &  $7.54$  & $592745$ & $7.53$  & $607481$ & $5.23$  & $607709$ & $5.19$ \\
				$16\!\times\!16$ & $620750$    &  $3.16$  & $620909$ & $3.13$  & $624429$ & $2.59$  & $624490$ & $2.58$ \\
				$32\!\times\!32$ & $628280$    &  $1.98$  & $628908$ & $1.89$  & $629006$ & $1.87$  & $629044$ & $1.87$ \\  
				\hline
				\multicolumn{9}{c}{Reference solution: $n_y(B)=641000\,\text{Nm}^{-1}$ (theory of shallow shell, \cite{Scordelis})} \\ 
				\hline
			\end{tabularx} \nonumber
			\label{tab10}
			\caption{Internal force $n_y(B)$ for the Scordelis-Lo roof problem}
		\end{table}
	\end{center}
	
	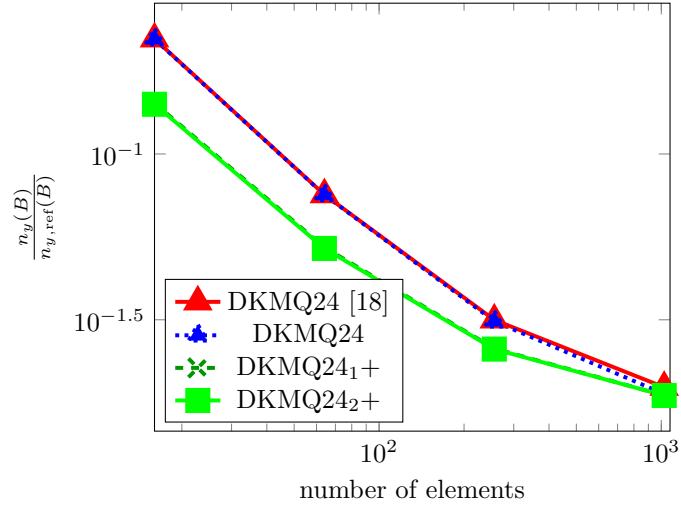
\begin{figure}
		\centering
		\begin{tikzpicture} 
		\begin{loglogaxis}[ 
		xlabel=\text{number of elements}, 
		ylabel= $\frac{n_y(B)}{n_{y,\rf}(B)}$, 
		x tick label style={
			/pgf/number format/1000 sep={}  
		},
		y tick label style={
			/pgf/number format/fixed,       
			/pgf/number format/precision=5  
		},
		scaled y ticks=false, 
		xmin= 16,
		xmax = 1074,
		legend style={at={(0.02,0.02)},anchor=south west}
		]     
		\addplot[line width=1.5pt, color=red, mark=triangle*, mark size=5pt] coordinates { 
			(16,   22.24E-2) 
			(64,   7.54E-2) 
			(256,  3.16E-2)
			(1024, 1.98E-2) 
		};   
		\addplot[line width=1.5pt, color=blue, mark=triangle*, mark size=4pt, dotted] coordinates { 
			(16,   22.19E-2) 
			(64,   7.53E-2) 
			(256,  3.13E-2) 
			(1024, 1.89E-2)
		};
		\addplot[line width=1.5pt, color=col1, mark=x, mark size=5pt, dashed] coordinates { 
			(16,   14.23E-2) 
			(64,   5.23E-2) 
			(256,  2.59E-2) 
			(1024, 1.87E-2)
		};  
		\addplot[line width=1.5pt, color=green, mark=square*, mark size=4pt] coordinates { 
			(16,   14.12E-2) 
			(64,   5.19E-2) 
			(256,  2.58E-2) 
			(1024, 1.87E-2)
		};    
		\legend{DKMQ24 \cite{Katili2015}, DKMQ24, DKMQ24$_1$+ , DKMQ24$_2$+}
		\end{loglogaxis} 
		\end{tikzpicture} 
		\caption{Relative error of internal force $n_y(B)$ for the Scordelis-Lo roof problem}
		\label{fig:BENCH3d}
	\end{figure}
	
	\begin{center}
		\begin{table}
			\centering
			\begin{tabularx}{160mm}{ c *{10}{Y} }
				\hline
				Mesh  
				& \multicolumn{2}{c}{MITC4+}  
				& \multicolumn{2}{c}{US-ATFHS8}
				& \multicolumn{2}{c}{DKMQ24}  
				& \multicolumn{2}{c}{DKMQ24$_1$+} 
				& \multicolumn{2}{c}{DKMQ24$_2$+}  \rule[-2mm]{0mm}{6mm}  \\
				& $u_Z(B)$   & err. &    $u_Z(B)$   & err. &    $u_Z(B)$   & err.  &    $u_Z(B)$   & err. &    $u_Z(B)$   & err.\\
				& [m] &  [\%]     & [m] &  [\%]    & [m] &  [\%]  & [m] &  [\%] & [m] &  [\%] \\
				\hline
				\hline
				$4\!\times\!4$   & $-0.31692$ & $4.80$ & $-0.31767$ & $5.05$ & $-0.30825$ & $1.93$ &  $-0.30832$ & $1.96$ &  $-0.30367$ & $0.42$ \\
				$8\!\times\!8$   & $-0.30391$ & $0.50$ & $-0.30443$ & $0.67$ & $-0.30180$ & $0.20$ &  $-0.30185$ & $0.18$ &  $-0.30069$ & $0.56$ \\
				$16\!\times\!16$ & $-0.30158$ & $0.27$ & $-0.30173$ & $0.22$ & $-0.30083$ & $0.52$ &  $-0.30085$ & $0.51$ &  $-0.30056$ & $0.61$ \\
				$32\!\times\!32$ & $-0.30113$ & $0.42$ &     -      &     -  & $-0.30071$ & $0.56$ &  $-0.30072$ & $0.56$ &  $-0.30065$ & $0.58$ \\
				\hline
				\multicolumn{11}{c}{Reference solution: $u_Z(B)=-0.3024\,\text{m}$ (\cite{MacNealHarder})} \\ 
				\hline
			\end{tabularx} \nonumber
			\label{tab10}
			\caption{Deflection $u_Z(B)$ for the Scordelis-Lo roof problem with modified parameters according to \cite{MacNealHarder} ($L=50\,\text{m},\,R=25\,\text{m},\,h=0.25\,\text{m},\,E=432\,\text{MPa},\nu=0,\,\rho=36\,\text{kg\,m}^{-3},\,g=10\,\text{m\,s}^{-2}$)}
		\end{table}
	\end{center}
	
	\begin{figure}
		\centering
		\begin{tikzpicture} 
		\begin{loglogaxis}[ 
		xlabel=\text{number of elements}, 
		ylabel= $\frac{u_Z(B)}{u_{Z,\rf}(B)}$, 
		x tick label style={
			/pgf/number format/1000 sep={}  
		},
		y tick label style={
			/pgf/number format/fixed,       
			/pgf/number format/precision=5  
		},
		scaled y ticks=false, 
		xmin= 16,
		xmax = 1074, 
		legend style={at={(1.02,0.02)},anchor=south west}
		]     
		\addplot[line width=1.5pt, color=purple, mark=+, mark size=5pt] coordinates {  
			(16,    4.80E-2) 
			(64,    0.50E-2) 
			(256,   0.27E-2) 
			(1024,  0.42E-2)
		};
		\addplot[line width=1.5pt, color=orange, mark=x, mark size=5pt] coordinates { 
			(16,    5.05E-2) 
			(64,    0.67E-2) 
			(256,   0.22E-2) 
		};
		\addplot[line width=1.5pt, color=blue, mark=triangle*, mark size=4pt, dotted] coordinates { 
			(16,    5.71E-2) 
			(64,    2.80E-2) 
			(256,   1.19E-2) 
			(1024,  0.60E-2)
		};
		\addplot[line width=1.5pt, color=col1, mark=x, mark size=5pt, dashed] coordinates { 
			(16,    1.93E-2) 
			(64,    0.20E-2) 
			(256,   0.52E-2) 
			(1024,  0.56E-2)
		}; 
		\addplot[line width=1.5pt, color=green, mark=square*, mark size=4pt] coordinates { 
			(16,    1.96E-2) 
			(64,    0.18E-2) 
			(256,   0.51E-2) 
			(1024,  0.56E-2)
		};    
		\legend{MITC4+, US-ATFHS8, DKMQ24, DKMQ24$_1$+ , DKMQ24$_2$+}
		\end{loglogaxis} 
		\end{tikzpicture} 
		\caption{Relative error of deflection $u_Z(B)$ for the Scordelis-Lo roof problem with modified parameters according to \cite{MacNealHarder} ($L=50\,\text{m},\,R=25\,\text{m},\,h=0.25\,\text{m},\,E=432\,\text{MPa},\nu=0,\,\rho=36\,\text{kg\,m}^{-3},\,g=10\,\text{m\,s}^{-2}$)}
		\label{fig:BENCH3ba}
	\end{figure}
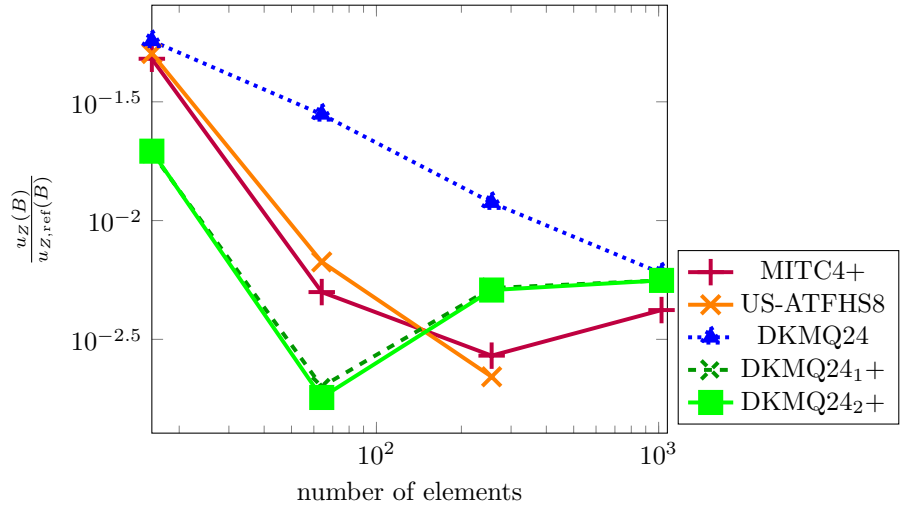
	
	\FloatBarrier
	
	\subsection{Pinched cylinder shell}
	\label{b4}
	
	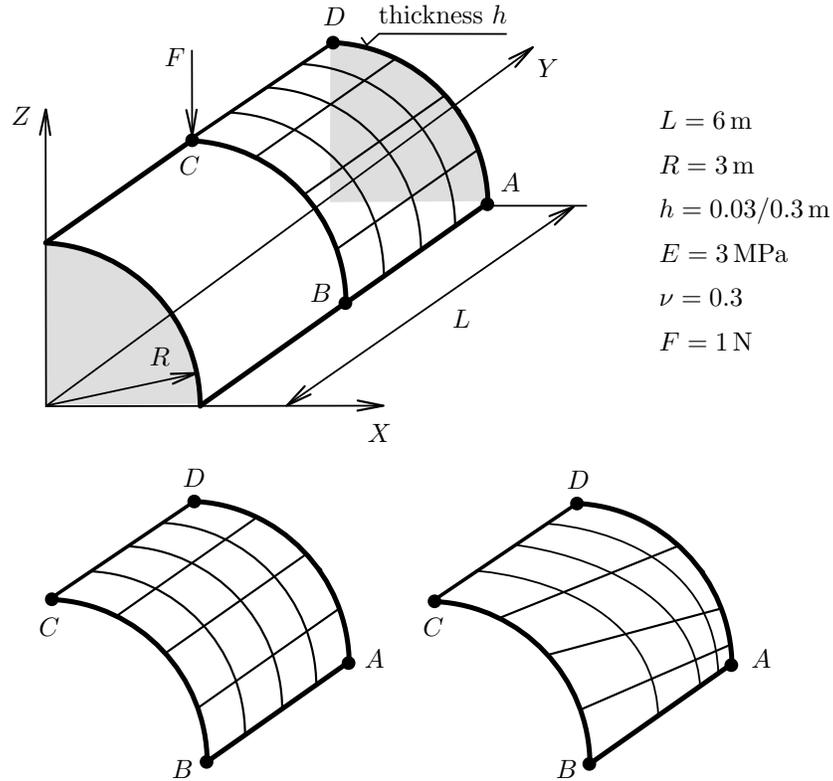
\begin{figure}[H]
		\centering
		\begin{overpic}[width=0.575\textwidth]{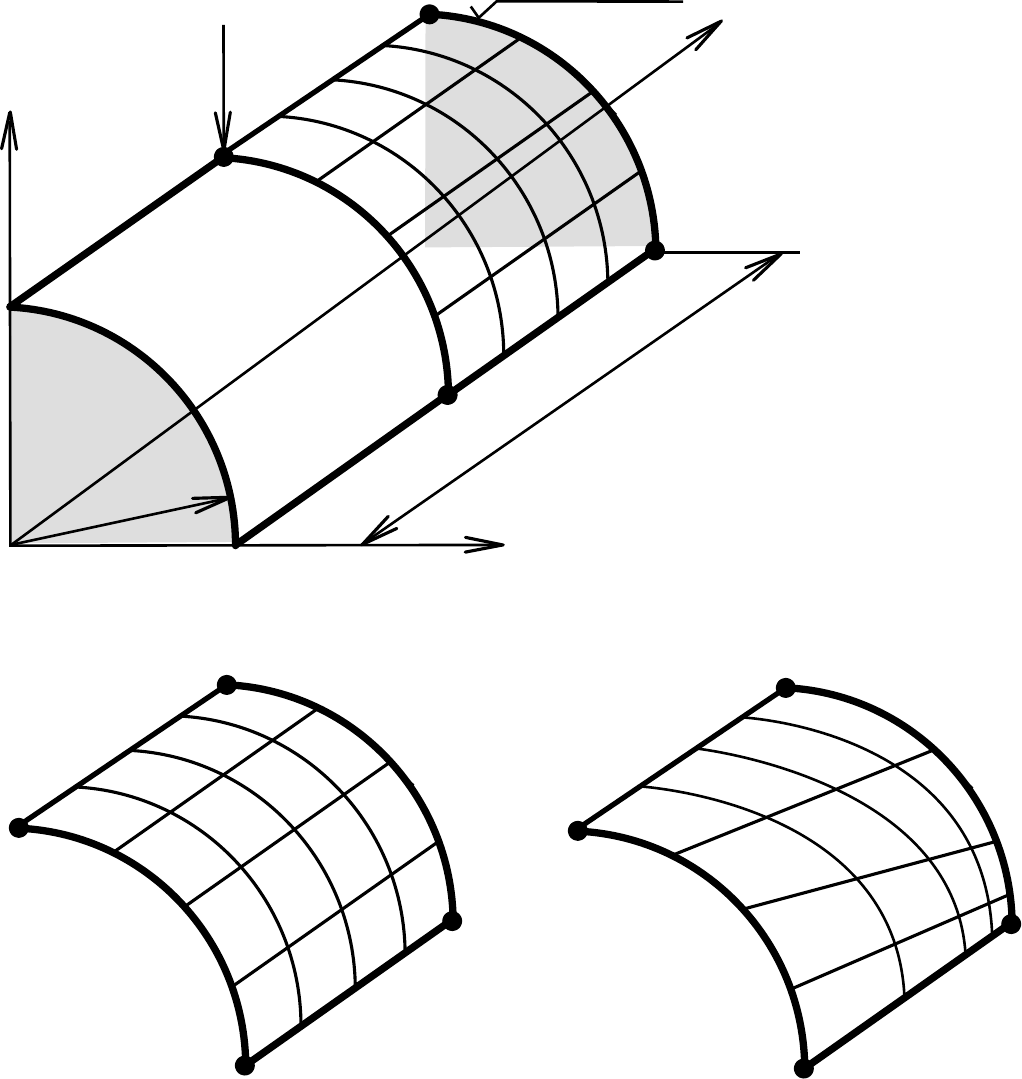} 
			\put(62.5,78.5){$A$}
			\put(36.75,63.5){$B$}
			\put(19,80.75){$C$}
			\put(38.5,101){$D$}
			\put(44.5,44.5){$X$}
			\put(67.5,94){$Y$}
			\put(-3.75,87.5){$Z$}
			\put(15,55){$R$}
			\put(17,95.5){$F$}
			\put(56,60){$L$}		
			\put(44.1,14){$A$}
			\put(18,-1){$B$}
			\put(0,18.25){$C$}
			\put(19.5,38.5){$D$}
			\put(96.5,14){$A$}
			\put(70,-1){$B$}
			\put(52,18.25){$C$}
			\put(71.5,38.25){$D$}
			\put(46,101.25){$\text{thickness}~h$}	
			\put(84,87){$L=6\,\text{m}$}
			\put(84,81){$R=3\,\text{m}$}
			\put(84,75){$h=0.03/0.3\,\text{m}$}
			\put(84,69){$E=3\,\text{MPa}$}
			\put(84,63){$\nu=0.3$}	
			\put(84,57){$F=1\,\text{N}$}
		\end{overpic}
		\caption{Pinched cylinder shell. Lower left: uniform mesh $4\times4$, lower right: distorted mesh  $4\times4$ with the aspect ratio 10 on sides AB and AD.}
		\label{PinchedCylinder}	
	\end{figure}
	\FloatBarrier
	The pinched cylinder, depicted in Figure \ref{PinchedCylinder}, is loaded by two opposite pinching forces $F$. The cylinder has a rigid diaphragm at both ends. Only one eighth of the model is calculated due to symmetries (to which the force $\frac{F}{4}$ is applied only). The boundary conditions are as follows: $u_X=u_Z=\varphi_Y=0$ on the side AD, on the remaining sides symmetric conditions are prescribed: $u_Z=\varphi_X=\varphi_Y=0$ on side AB, $u_Y=\varphi_X=\varphi_Z=0$ on the side BC and $u_X=\varphi_Y=\varphi_Z=0$ on side CD. Two different shell thicknesses are considered in our calculation, $h=0.03\,\text{m}$ and $h=0.3\,\text{m}$. We use Young's modulus $E=3\,\text{MPa}$ instead of the value $30000\,\text{MPa}$ used in \cite{Katili2015}.
	The uniform and distorted meshes are analysed. For the distorted mesh a bias factor, i.e. the ratio of the largest and smallest element edge, of 10 is used. The analytical solution for the case $h=0.03\,\text{m}$ is based on the Kirchhoff theory proposed by Lindberg, et al. \cite{Ref37}, which follows the solution of Fl\"{u}gge \cite{Ref38}, and equals to $u_Z(C)=-1.825\,\text{mm}$. The analytical solution for the case $h=0.3\,\text{m}$, taken from \cite{Ref39}, is based on the Mindlin theory and equals to $u_Z(C)=-0.01261\,\text{mm}$. However, as reported in \cite{Katili2015} and according to our numerical results, the numerical results converge to slightly larger values. Therefore, numerical results of the DKMQ24$_2$+ shell element on mesh $128\!\times\!128$ are taken as reference values. 
	The numerical results presented in Tables 9-12 reveal slow convergence of all the considered elements, which is in concordance with other authors \cite{Belytschko, MITC4Pinched} and is a consequence of the single point loading.
	
	When compared to the DKMQ24 \cite{Katili2015}  element, our modified elements DKMQ24$_1$+ and DKMQ24$_2$+ have almost identical results (for $h=0.03\,\text{m}$ at some cases  slightly better, at some slightly worse, for $h=0.3\,\text{m}$ slightly better at all test cases). When compared to the MITC4+ element, where we only have data for the case $h=0.03\,\text{m}$, the modified elements DKMQ24$_1$+ and DKMQ24$_2$+ exhibit lower errors in case of the uniform mesh, and slightly higher errors in case of the distorted mesh.
	
	
	
	

	\begin{center}
		\begin{table}
			\centering
			\begin{tabularx}{160mm}{ c *{10}{Y} }
				\hline
				Mesh  
				& \multicolumn{2}{c}{MITC4+}  
				& \multicolumn{2}{c}{DKMQ24 \cite{Katili2015}}  
				& \multicolumn{2}{c}{DKMQ24}  
				& \multicolumn{2}{c}{DKMQ24$_1$+} 
				& \multicolumn{2}{c}{DKMQ24$_2$+}  \rule[-2mm]{0mm}{6mm}  \\ 
				& $u_Z(C)$ & err. & $u_Z(C)$ & err. & $u_Z(C)$ & err. & $u_Z(C)$ & err. & $u_Z(C)$ & err. \\
				& [mm] &  [\%]    & [mm] &  [\%]    & [mm] &  [\%]     & [mm] &  [\%]    & [mm] &  [\%]  \\
				\hline
				\hline
				$4\!\times4\!$     &  $-0.71240$ & $60.03$ & $-1.1255$ & $39.26$  &    $-1.12416$   &  $39.33$    &  $-1.18399$  &  $36.10$  &  $-1.18474$  &  $36.06$      \\
				$8\!\times8\!$     &  $-1.37736$ & $24.15$ & $-1.7246$ &  $6.93$  &    $-1.72382$   &   $6.97$    &  $-1.74588$  &   $5.78$  &  $-1.74647$  &   $5.75$      \\
				$16\!\times16\!$   &  $-1.69944$ & $6.77$  & $-1.8599$ &  $0.37$  &    $-1.85952$   &   $0.35$    &  $-1.86512$  &   $0.65$  &  $-1.86557$  &   $0.68$      \\
				$32\!\times32\!$   &  $-1.80254$ & $1.20$  &     -     &     -    &    $-1.85688$   &   $0.21$    &  $-1.85939$  &   $0.34$  &  $-1.85951$  &   $0.35$      \\
				$64\!\times64\!$   &  $-1.83392$ & $0.49$  &     -     &     -    &    $-1.85151$   &   $0.08$    &  $-1.85444$  &   $0.08$  &  $-1.85445$  &   $0.08$      \\
				\hline
				\multicolumn{11}{c}{Reference solution: $u_Z(C)=-1.853\,\text{mm}$ (DKMQ24$_2$+ element on mesh $128\!\times\!128$)} \\ 
				\hline
			\end{tabularx} \nonumber
			\caption{Deflection $u_Z(C)$ for the pinched cylinder problem, thickness $h=0.03\,\text{m}$, uniform mesh}
		\end{table}
	\end{center}

	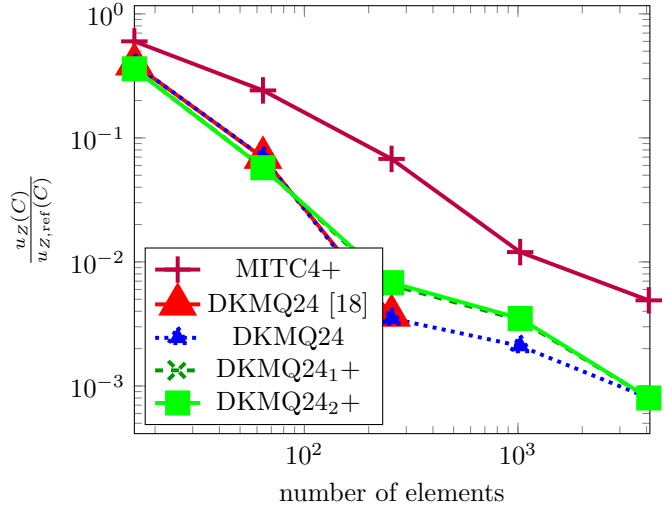
\begin{figure}
		\centering
		\begin{tikzpicture} 
		\begin{loglogaxis}[ 
		xlabel=\text{number of elements}, 
		ylabel= $\frac{u_Z(C)}{u_{Z,\rf}(C)}$, 
		x tick label style={
			/pgf/number format/1000 sep={}  
		},
		y tick label style={
			/pgf/number format/fixed,       
			/pgf/number format/precision=5  
		},
		scaled y ticks=false, 
		xmin= 16,
		xmax = 4146,
		legend style={at={(0.02,0.02)},anchor=south west}
		]    
		\addplot[line width=1.5pt, color=purple, mark=+, mark size=5pt] coordinates {
			(16,  60.03E-2) 
			(64,  24.15E-2) 
			(256,  6.77E-2) 
			(1024, 1.20E-2)
			(4096, 0.49E-2)
		};    
		\addplot[line width=1.5pt, color=red, mark=triangle*, mark size=7pt] coordinates { 
			(16,  39.26E-2) 
			(64,  6.93E-2) 
			(256, 0.37E-2) 
		};    
		\addplot[line width=1.5pt, color=blue, mark=triangle*, mark size=4pt, dotted] coordinates { 
			(16,  39.33E-2) 
			(64,  6.97E-2)  
			(256, 0.35E-2) 
			(1024,0.21E-2)
			(4096,0.08E-2)
		};    
		\addplot[line width=1.5pt, color=col1, mark=x, mark size=5pt, dashed] coordinates { 
			(16,  36.10E-2) 
			(64,  5.78E-2) 
			(256, 0.65E-2) 
			(1024,0.34E-2)
			(4096,0.08E-2)
		};    
		\addplot[line width=1.5pt, color=green, mark=square*, mark size=4pt] coordinates { 
			(16,  36.06E-2) 
			(64,  5.75E-2) 
			(256, 0.68E-2) 
			(1024,0.35E-2)
			(4096,0.08E-2)
		};    
		\legend{MITC4+, DKMQ24 \cite{Katili2015}, DKMQ24, DKMQ24$_1$+ , DKMQ24$_2$+}
		\end{loglogaxis} 
		\end{tikzpicture} 
		\caption{Relative error of deflection $u_Z(C)$ for the pinched cylinder problem, thickness $h=0.03\,\text{m}$, uniform mesh}
		\label{fig:BENCH4a}
	\end{figure}
	
	
	\begin{center}
		\begin{table}
			\centering
			\begin{tabularx}{160mm}{ c *{10}{Y} }
				\hline
				Mesh  
				& \multicolumn{2}{c}{MITC4+}  
				& \multicolumn{2}{c}{DKMQ24 \cite{Katili2015}}  
				& \multicolumn{2}{c}{DKMQ24}  
				& \multicolumn{2}{c}{DKMQ24$_1$+} 
				& \multicolumn{2}{c}{DKMQ24$_2$+}  \rule[-2mm]{0mm}{6mm}   \\
				& $u_Z(C)$ & err. & $u_Z(C)$ & err. & $u_Z(C)$ & err. & $u_Z(C)$ & err. & $u_Z(C)$ & err. \\
				& [mm] &  [\%]    & [mm] &  [\%]    & [mm] &  [\%]     & [mm] &  [\%]    & [mm] &  [\%]  \\
				\hline
				\hline
				$4\!\times4\!$     &  $-0.69215$ & $61.13$ & $-0.2397$ & $87.06$  &    $-0.24254$   &  $86.91$    &  $-0.24973$  &  $86.52$  &  $-0.25396$  &  $86.29$      \\
				$8\!\times8\!$     &  $-1.37499$ & $24.27$ & $-1.2326$ & $33.46$  &    $-1.23106$   &  $33.56$    &  $-1.19506$  &  $35.51$  &  $-1.19549$  &  $35.48$      \\
				$16\!\times16\!$   &  $-1.70090$ & $6.69$  & $-1.7144$ &  $7.48$  &    $-1.71369$   &   $7.52$    &  $-1.70823$  &   $7.81$  &  $-1.70857$  &   $7.79$      \\
				$32\!\times32\!$   &  $-1.80400$ & $1.12$  &     -     &     -    &    $-1.82068$   &   $1.74$    &  $-1.82087$  &   $1.73$  &  $-1.82104$  &   $1.72$      \\
				$64\!\times64\!$   &  $-1.83392$ & $0.49$  &     -     &     -    &    $-1.84221$   &   $0.58$    &  $-1.84481$  &   $0.44$  &  $-1.84484$  &   $0.44$      \\
				\hline
				\multicolumn{11}{c}{Reference solution: $u_Z(C)=-1.853\,\text{mm}$ (DKMQ24$_2$+ element on mesh $128\!\times\!128$)} \\ 
				\hline
			\end{tabularx} \nonumber
			\caption{Deflection $u_Z(C)$ for the pinched cylinder problem, thickness $h=0.03\,\text{m}$, distorted mesh}
		\end{table}
	\end{center}
	
	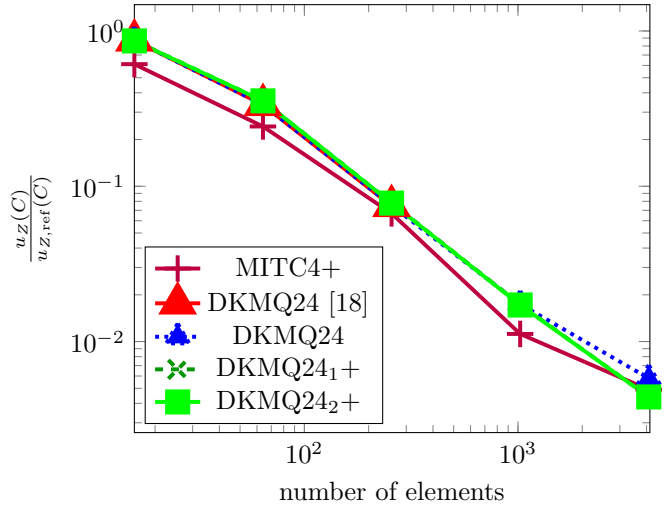
\begin{figure}
		\centering
		\begin{tikzpicture} 
		\begin{loglogaxis}[ 
		xlabel=\text{number of elements}, 
		ylabel= $\frac{u_Z(C)}{u_{Z,\rf}(C)}$, 
		x tick label style={
			/pgf/number format/1000 sep={}  
		},
		y tick label style={
			/pgf/number format/fixed,       
			/pgf/number format/precision=5  
		},
		scaled y ticks=false, 
		xmin= 16,
		xmax = 4146,
		legend style={at={(0.02,0.02)},anchor=south west}
		]    
		\addplot[line width=1.5pt, color=purple, mark=+, mark size=5pt] coordinates {
			(16,  61.13E-2) 
			(64,  24.27E-2) 
			(256, 6.69E-2) 
			(1024,1.12E-2)
			(4096,0.49E-2)
		};        
		\addplot[line width=1.5pt, color=red, mark=triangle*, mark size=7pt] coordinates { 
			(16,  87.06E-2) 
			(64,  33.46E-2) 
			(256, 7.48E-2) 
		};      
		\addplot[line width=1.5pt, color=blue, mark=triangle*, mark size=5pt, dotted] coordinates { 
			(16,  86.91E-2) 
			(64,  33.56E-2) 
			(256, 7.52E-2) 
			(1024,1.74E-2)
			(4096,0.58E-2)
		};        
		\addplot[line width=1.5pt, color=col1, mark=x, mark size=5pt, dashed] coordinates { 
			(16,  86.52E-2) 
			(64,  35.51E-2) 
			(256, 7.81E-2) 
			(1024,1.73E-2)
			(4096,0.44E-2)
		};      
		\addplot[line width=1.5pt, color=green, mark=square*, mark size=4pt] coordinates { 
			(16,  86.29E-2) 
			(64,  35.48E-2) 
			(256, 7.79E-2) 
			(1024,1.72E-2)
			(4096,0.44E-2)
		};    
		\legend{MITC4+, DKMQ24 \cite{Katili2015}, DKMQ24, DKMQ24$_1$+ , DKMQ24$_2$+}
		\end{loglogaxis} 
		\end{tikzpicture} 
		\caption{Relative error of deflection $u_Z(C)$ for the pinched cylinder problem, thickness $h=0.03\,\text{m}$, distorted mesh}
		\label{fig:BENCH4b}
	\end{figure}
	
	
	\begin{center}
		\begin{table}
			\centering
			\begin{tabularx}{148mm}{ c *{8}{Y} }
				\hline
				Mesh  
				& \multicolumn{2}{c}{DKMQ24 \cite{Katili2015}}  
				& \multicolumn{2}{c}{DKMQ24}  
				& \multicolumn{2}{c}{DKMQ24$_1$+} 
				& \multicolumn{2}{c}{DKMQ24$_2$+}  \rule[-2mm]{0mm}{6mm}   \\
				& $u_Z(C)$ & err. & $u_Z(C)$ & err. & $u_Z(C)$ & err. & $u_Z(C)$ & err.    \\
				& [mm] &  [\%]    & [mm] &  [\%]    & [mm] &  [\%]     & [mm] &  [\%]      \\
				\hline
				\hline
				$4\!\times4\!$     & $-0.011445$ & $19.97$ &   $-0.011435$ & $20.04$   &  $-0.011622$  &   $18.73$  &  $-0.011624$  &   $18.71$      \\
				$8\!\times8\!$     & $-0.012180$ & $14.83$ &   $-0.012125$ & $15.21$   &  $-0.012272$  &   $14.18$  &  $-0.012273$  &   $14.18$      \\
				$16\!\times16\!$   & $-0.012822$ & $10.34$ &   $-0.012638$ & $11.62$   &  $-0.012822$  &   $10.34$  &  $-0.012823$  &   $10.33$      \\
				$32\!\times32\!$   &      -      &     -   &   $-0.013099$ & $8.40$    &  $-0.013474$  &    $5.78$  &  $-0.013462$  &    $5.86$      \\
				$64\!\times64\!$   &      -      &     -   &   $-0.013505$ & $5.56$    &  $-0.014191$  &    $0.76$  &  $-0.014163$  &    $0.96$      \\
				\hline
				\multicolumn{9}{c}{Reference solution: $u_Z(C)=-0.01430\,\text{mm}$ (DKMQ24$_2$+ element on mesh $128\!\times\!128$)} \\ 
				\hline
			\end{tabularx} \nonumber
			\caption{Deflection $u_Z(C)$ for the pinched cylinder problem, thickness $h=0.3\,\text{m}$, uniform mesh}
		\end{table}
	\end{center}

	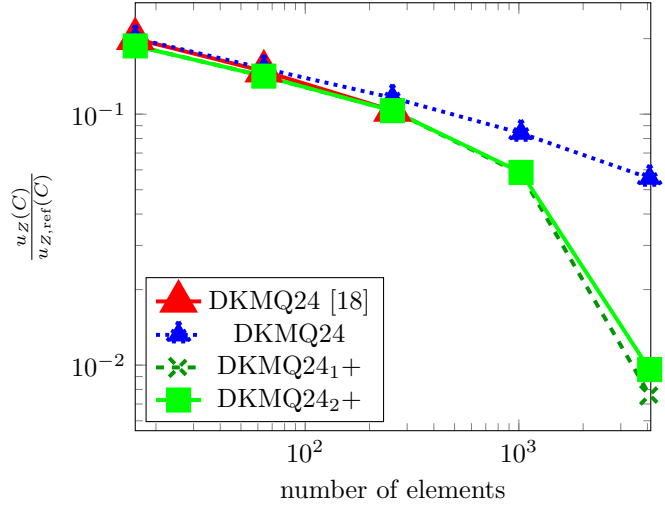
\begin{figure}
		\centering
		\begin{tikzpicture} 
		\begin{loglogaxis}[ 
		xlabel=\text{number of elements}, 
		ylabel= $\frac{u_Z(C)}{u_{Z,\rf}(C)}$, 
		x tick label style={
			/pgf/number format/1000 sep={}  
		},
		y tick label style={
			/pgf/number format/fixed,       
			/pgf/number format/precision=5  
		},
		scaled y ticks=false, 
		xmin= 16,
		xmax = 4146,
		legend style={at={(0.02,0.02)},anchor=south west}
		]          
		\addplot[line width=1.5pt, color=red, mark=triangle*, mark size=7pt] coordinates { 
			(16,  19.97E-2) 
			(64,  14.83E-2) 
			(256, 10.34E-2) 
		};         
		\addplot[line width=1.5pt, color=blue, mark=triangle*, mark size=5pt, dotted] coordinates { 
			(16,  20.04E-2)     
			(64,  15.21E-2) 
			(256, 11.62E-2) 
			(1024, 8.40E-2)
			(4096, 5.56E-2)
		};           
		\addplot[line width=1.5pt, color=col1, mark=x, mark size=5pt, dashed] coordinates { 
			(16,  18.73E-2) 
			(64,  14.18E-2) 
			(256, 10.34E-2) 
			(1024, 5.78E-2)
			(4096, 0.76E-2)
		};        
		\addplot[line width=1.5pt, color=green, mark=square*, mark size=4pt] coordinates { 
			(16,  18.71E-2) 
			(64,  14.18E-2) 
			(256, 10.33E-2) 
			(1024, 5.86E-2)
			(4096, 0.96E-2)
		};    
		\legend{DKMQ24 \cite{Katili2015}, DKMQ24, DKMQ24$_1$+ , DKMQ24$_2$+}
		\end{loglogaxis} 
		\end{tikzpicture} 
		\caption{Relative error of deflection $u_Z(C)$ for the pinched cylinder problem, thickness $h=0.3\,\text{m}$, uniform mesh}
		\label{fig:BENCH4c}
	\end{figure}
	

	\begin{center}
		\begin{table}
			\centering
			\begin{tabularx}{148mm}{ c *{8}{Y} }
				\hline
				Mesh  
				& \multicolumn{2}{c}{DKMQ24 \cite{Katili2015}}  
				& \multicolumn{2}{c}{DKMQ24}  
				& \multicolumn{2}{c}{DKMQ24$_1$+} 
				& \multicolumn{2}{c}{DKMQ24$_2$+}  \rule[-2mm]{0mm}{6mm}   \\
				& $u_Z(C)$ & err. & $u_Z(C)$ & err. & $u_Z(C)$ & err. & $u_Z(C)$ & err.    \\
				& [mm] &  [\%]    & [mm] &  [\%]    & [mm] &  [\%]     & [mm] &  [\%]      \\
				\hline
				\hline
				$4\!\times4\!$     & $-0.008617$ & $39.74$ &   $-0.008675$  &  $39.33$   &  $-0.008857$  &   $38.07$  &  $-0.008858$  &   $38.05$      \\
				$8\!\times8\!$     & $-0.011443$ & $19.98$ &   $-0.011444$  &  $19.97$   &  $-0.011767$  &   $17.71$  &  $-0.011768$  &   $17.70$      \\
				$16\!\times16\!$   & $-0.012508$ & $12.53$ &   $-0.012428$  &  $13.09$   &  $-0.012748$  &   $10.85$  &  $-0.012750$  &   $10.84$      \\
				$32\!\times32\!$   &      -      &     -   &   $-0.013039$  &   $8.82$   &  $-0.013472$  &    $5.79$  &  $-0.013464$  &    $5.85$      \\
				$64\!\times64\!$   &      -      &     -   &   $-0.013499$  &   $5.60$   &  $-0.014190$  &    $0.77$  &  $-0.014161$  &    $0.97$      \\
				\hline
				\multicolumn{9}{c}{Reference solution: $u_Z(C)=-0.01430\,\text{mm}$ (DKMQ24$_2$+ element on mesh $128\!\times\!128$)} \\ 
				\hline
			\end{tabularx} \nonumber
			\caption{Deflection $u_Z(C)$ for the pinched cylinder problem, thickness $h=0.3\,\text{m}$, distorted mesh}
		\end{table}
	\end{center}

	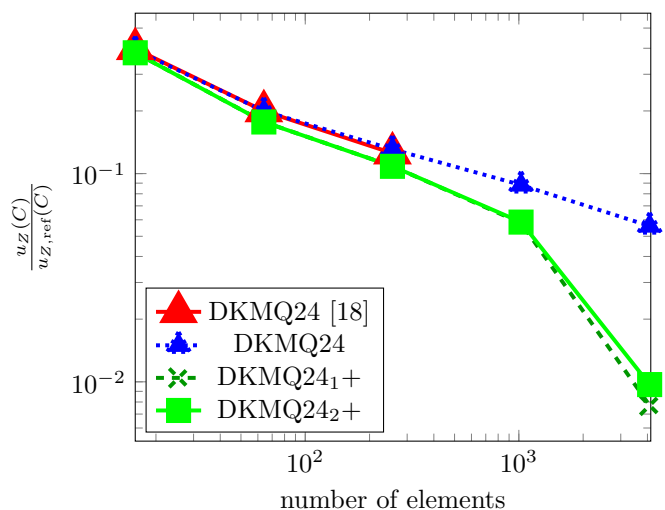
\begin{figure}
		\centering
		\begin{tikzpicture} 
		\begin{loglogaxis}[ 
		xlabel=\text{number of elements}, 
		ylabel= $\frac{u_Z(C)}{u_{Z,\rf}(C)}$, 
		x tick label style={
			/pgf/number format/1000 sep={}  
		},
		y tick label style={
			/pgf/number format/fixed,       
			/pgf/number format/precision=5  
		},
		scaled y ticks=false, 
		xmin= 16,
		xmax = 4146,
		legend style={at={(0.02,0.02)},anchor=south west}
		]          
		\addplot[line width=1.5pt, color=red, mark=triangle*, mark size=7pt] coordinates { 
			(16,  39.74E-2) 
			(64,  19.98E-2) 
			(256, 12.53E-2) 
		};         
		\addplot[line width=1.5pt, color=blue, mark=triangle*, mark size=5pt, dotted] coordinates { 
			(16,   39.33E-2) 
			(64,   19.97E-2) 
			(256,  13.09E-2) 
			(1024, 8.82E-2)
			(4096, 5.60E-2)
		};           
		\addplot[line width=1.5pt, color=col1, mark=x, mark size=5pt, dashed] coordinates { 
			(16,   38.07E-2) 
			(64,   17.71E-2) 
			(256,  10.85E-2) 
			(1024,  5.79E-2)
			(4096,  0.77E-2)
		};        
		\addplot[line width=1.5pt, color=green, mark=square*, mark size=4pt] coordinates { 
			(16,   38.05E-2) 
			(64,   17.70E-2) 
			(256,  10.84E-2) 
			(1024, 5.85E-2)
			(4096, 0.97E-2)
		};    
		\legend{DKMQ24 \cite{Katili2015}, DKMQ24, DKMQ24$_1$+ , DKMQ24$_2$+}
		\end{loglogaxis} 
		\end{tikzpicture} 
		\caption{Relative error of deflection $u_Z(C)$ for the pinched cylinder problem, thickness $h=0.3\,\text{m}$, distorted mesh}
		\label{fig:BENCH4d}
	\end{figure}
	
	\FloatBarrier
	
	\subsection{Hyperbolic paraboloid shell}
	\label{b5}
	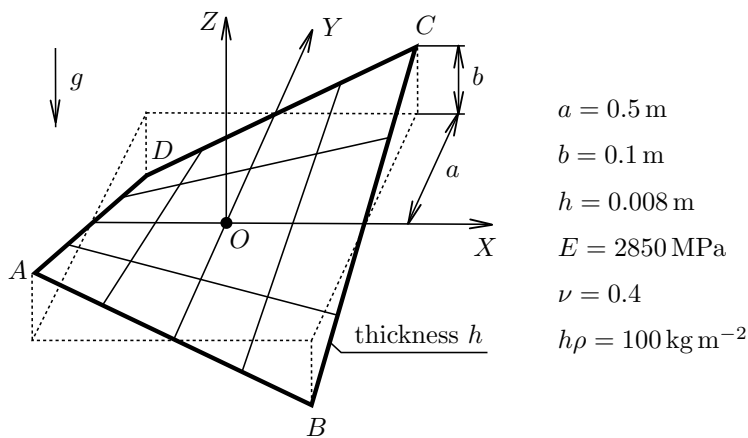
\begin{figure}[H]
		\centering
		\begin{overpic}[width=0.38\textwidth]{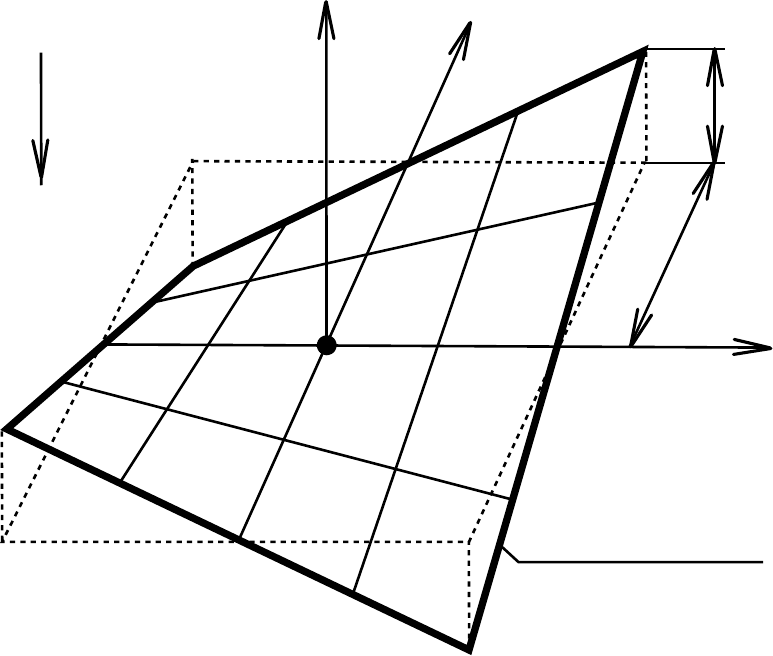} 
			\put(8.5,70){$g$} 
			\put(43,34.6){$O$}
			\put(-5,28){$A$}
			\put(59.5,-6){$B$}
			\put(83.2,81){$C$}
			\put(26,54){$D$}
			\put(95.8,33){$X$}
			\put(63,80.6){$Y$}
			\put(36.4,81.8){$Z$}
			\put(89.8,50){$a$}
			\put(95.5,70){$b$}
			\put(70,14){thickness $h$}
			\put(114,63){$a=0.5\,\text{m}$}
			\put(114,53){$b=0.1\,\text{m}$}
			\put(114,43){$h=0.008\,\text{m}$}
			\put(114,33){$E=2850\,\text{MPa}$}
			\put(114,23){$\nu=0.4$}
			\put(114,13){$h\rho=100\,\text{kg}\,\text{m}^{-2}$}
		\end{overpic}
		\caption{Hyperbolic paraboloid shell with mesh $4\times 4$}
		\label{HyperbolicShell}
	\end{figure}
	\FloatBarrier
	The hyperbolic paraboloid shell (hypar shell) problem is depicted in Figure \ref{HyperbolicShell}. The shell middle surface is described by $Z=\frac{b}{a^2}XY$. The shell is loaded by self-weight, where the gravitational acceleration is considered to be $g=10\,\text{m}\, \text{s}^{-2}$. The shell is fully clamped on all sides ABCD, i.e. $u_X=u_Y=u_Z=\varphi_X=\varphi_Y=\varphi_Z=0$. The analytical solution is not available (\cite{Brebia}, p.245), therefore we use the numerical result on the fine mesh as reference. 
	Our modified elements DKMQ24$_1$+ and DKMQ24$_2$+ provide almost identical results at this benchmark when compared to the shell element DKMQ24 \cite{Katili2015}.
	
	\begin{center}
		\begin{table}
			\centering
			\begin{tabularx}{140mm}{ c *{8}{Y} }
				\hline
				Mesh           &  \multicolumn{2}{c}{DKMQ24 \cite{Katili2015}}  & \multicolumn{2}{c}{DKMQ24} & \multicolumn{2}{c}{DKMQ24$_1$+} & \multicolumn{2}{c}{DKMQ24$_2$+} \rule[-2mm]{0mm}{6mm}  \\
				&$u_Z(O)$   & err. &    $u_Z(O)$   & err.&    $u_Z(O)$   & err. &    $u_Z(O)$   & err. \\
				& [mm] &  [\%]    & [mm] &  [\%]    & [mm] &  [\%]    & [mm] &  [\%]   \\
				\hline
				\hline
				$2\!\times\!2$     & $0.35900$ & $44.52$ &   $0.35889$ & $44.47$  & $0.35874$ & $44.41$  & $0.35891$ & $44.48$  \\
				$4\!\times\!4$     & $0.27074$ &  $8.99$ &   $0.26938$ &  $8.44$  & $0.27064$ &  $8.94$  & $0.27055$ &  $8.91$  \\
				$8\!\times\!8$     & $0.25023$ &  $0.73$ &   $0.25030$ &  $0.76$  & $0.25013$ &  $0.69$  & $0.25019$ &  $0.71$  \\
				$16\!\times\!16$   & $0.24864$ &  $0.09$ &   $0.24866$ &  $0.10$  & $0.24859$ &  $0.07$  & $0.24868$ &  $0.10$  \\
				$32\!\times\!32$   &     -     &   -     &   $0.24847$ &  $0.02$  & $0.24845$ &  $0.01$  & $0.24846$ &  $0.02$  \\
				\hline                                                                                        
				\multicolumn{9}{c}{Reference solution: $u_Z(O)=0.248417\,\text{mm}$ (DKMQ24$_2$+ element on mesh $128\!\times\!128$)} \\ 
				\hline
			\end{tabularx} \nonumber
			\caption{Deflection $u_Z(O)$ for the hyperbolic paraboloid shell problem} 
		\end{table}
	\end{center}
	
	\FloatBarrier
	
	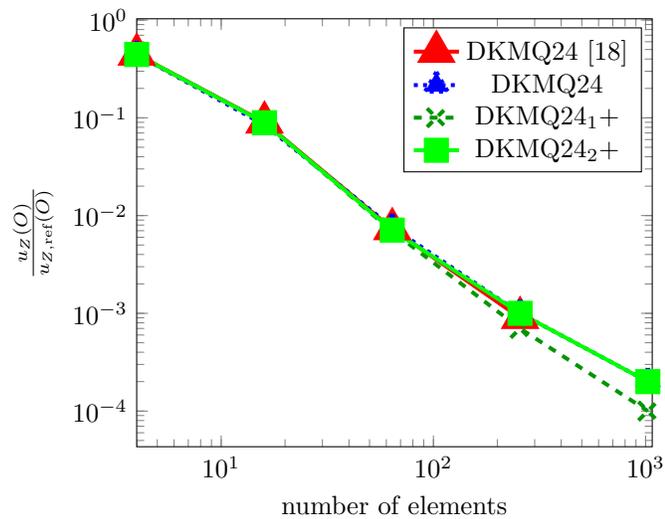
\begin{figure}
		\centering
		\begin{tikzpicture} 
		\begin{loglogaxis}[ 
		xlabel=\text{number of elements}, 
		ylabel= $\frac{u_Z(O)}{u_{Z,\rf}(O)}$, 
		x tick label style={
			/pgf/number format/1000 sep={}  
		},
		y tick label style={
			/pgf/number format/fixed,       
			/pgf/number format/precision=5  
		},
		scaled y ticks=false, 
		xmin= 4,
		xmax = 1074,
		legend style={at={(0.98,0.98)},anchor=north east}
		] 
		\addplot[line width=1.5pt, color=red, mark=triangle*, mark size=7pt] coordinates { 
			(4,   44.52E-2) 
			(16,  8.99E-2) 
			(64,  0.73E-2) 
			(256, 0.09E-2) 
		};    
		\addplot[line width=1.5pt, color=blue, mark=triangle*, mark size=5pt, dotted] coordinates { 
			(4,   44.47E-2) 
			(16,  8.44E-2) 
			(64,  0.76E-2) 
			(256, 0.10E-2) 
			(1024,0.02E-2)
		};    
		\addplot[line width=1.5pt, color=col1, mark=x, mark size=5pt, dashed] coordinates { 
			(4,   44.41E-2) 
			(16,  8.94E-2) 
			(64,  0.69E-2) 
			(256, 0.07E-2) 
			(1024,0.01E-2)
		};    
		\addplot[line width=1.5pt, color=green, mark=square*, mark size=4pt] coordinates { 
			(4,   44.48E-2) 
			(16,  8.91E-2) 
			(64,  0.71E-2) 
			(256, 0.10E-2) 
			(1024,0.02E-2)
		};  
		\legend{DKMQ24 \cite{Katili2015}, DKMQ24, DKMQ24$_1$+ , DKMQ24$_2$+}
		\end{loglogaxis} 
		\end{tikzpicture} 
		\caption{Relative error of deflection $u_Z(O)$ for the hyperbolic paraboloid shell problem}
		\label{fig:BENCH5a}
	\end{figure}
	
	\FloatBarrier
	
	\subsection{Twisted beam}
	\label{b6}
	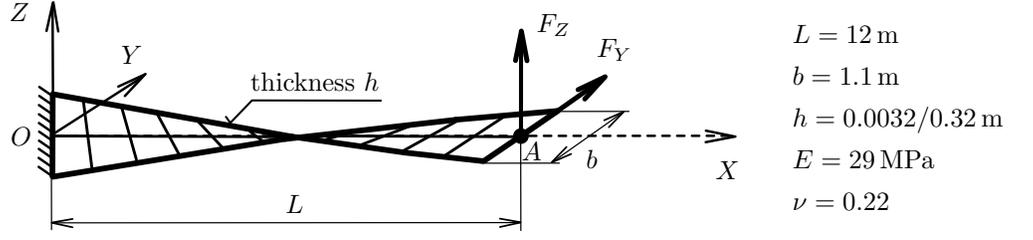
\begin{figure}[H]
		\centering
		\begin{overpic}[width=0.575\textwidth]{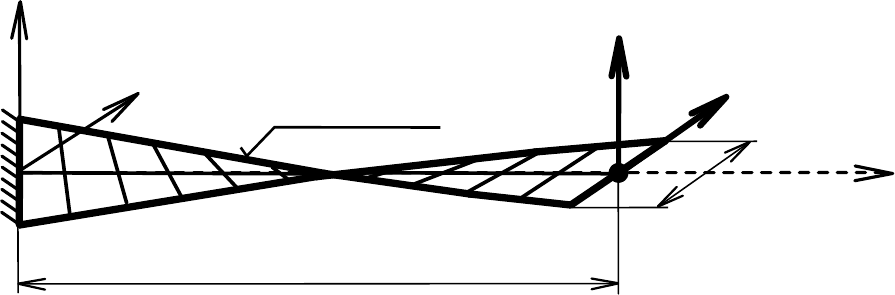} 
			\put(35.5,2.6){$L$}
			\put(-3.8,12.3){$O$}
			\put(69.3,10.1){$A$}
			\put(78.5,9){$b$}
			\put(80,25){$F_Y$}
			\put(71.4,28.5){$F_Z$}
			\put(97,7.5){$X$}
			\put(12,24){$Y$}
			\put(-4,30){$Z$}
			\put(30.5,20.1){thickness $h$}
			\put(108,27){$L=12\,\text{m}$}
			\put(108,21){$b=1.1\,\text{m}$}
			\put(108,15){$h=0.0032/0.32\,\text{m}$}
			\put(108,9){$E=29\,\text{MPa}$}
			\put(108,3){$\nu=0.22$}		    
		\end{overpic}
		\caption{Twisted beam with mesh $2\times 12$}
		\label{twisted}
	\end{figure}
	\FloatBarrier
	
	A shell strip of length $L$, width $b$ and thickness $h$ is twisted in the undeformed configuration by angle $\frac{\pi}{2}$ between points $O$ and $A$ (see Figure \ref{twisted}). The strip is clamped at $X\!=\!0$ and loaded by a single force $[0,F_Y, F_Z]$ at point $A$. The shell middle surface is described by: \mbox{$[X(\varphi,s),Y(\varphi,s),Z(\varphi,s)]=\left[\frac{2\varphi L}{\pi},~s\sin\varphi,~s\cos\varphi\right],$} $0\le \varphi \le\frac{\pi}{2},-\frac{b}{2}\le s \le\frac{b}{2}$.
	The example was introduced by MacNeal and Harder \cite{MacNealHarder} and represents an important test example at which many flat quadrilateral shell elements fail to converge. This is in accordance with our numerical experience, which shows that the DKMQ24/DKMQ24$_1$+/DKMQ24$_2$+ elements fail to converge if the coupling term, including derivatives of the element normal with respect to reference coordinates in Eq. (\ref{eqNormalDerivatives}, \ref{kappa}), which is identicaly zero for flat elements, is missing. The analytical solution based on the beam theory (which ignores the deformation in the transverse direction) is taken from \cite{TwistedBelytschko}.
	Four different settings of loading forces and shell thicknesses were tested, see Tables 17-20. 
	
	For the case $h=0.0032\,\text{m}$, when compared to  the DKMQ24 \cite{Katili2015} element, our modified elements DKMQ24$_1$+ and DKMQ24$_2$+ provide smaller errors at all mesh sizes except at the coarsest one. When comparing the same case $h=0.0032\,\text{m}$ with the MITC4+ element, our modified elements DKMQ24$_1$+ and DKMQ24$_2$+ provide smaller errors on the finer meshes. The smallest error at mesh size $4\times24$ is given by the \US element. For the case $h=0.0032\,\text{m}$ our modified elements DKMQ24$_1$+ and DKMQ24$_2$+ are superior to all other elements in all test cases.
	
	
	
	\FloatBarrier
	
	\begin{center}
		\begin{table}
			\centering
			\begin{tabularx}{160mm}{ c *{12}{Y} }
				\hline
				Mesh  & \multicolumn{2}{c}{MITC4+} & \multicolumn{2}{c}{US-ATFHS8}  &  \multicolumn{2}{c}{DKMQ24 \cite{Katili2015}}  & \multicolumn{2}{c}{DKMQ24} & \multicolumn{2}{c}{DKMQ24$_1$+} & \multicolumn{2}{c}{DKMQ24$_2$+} \rule[-2mm]{0mm}{6mm}  \\
				&    $u_Y(A)$   & err. &    $u_Y(A)$   & err. &    $u_Y(A)$   & err. &    $u_Y(A)$   & err.&    $u_Y(A)$   & err.&    $u_Y(A)$   & err. \\
				& [m] &  [\%]    & [m] &  [\%]    & [m] &  [\%]  & [m] &  [\%]  & [m] &  [\%]  & [m] &  [\%]    \\
				\hline
				\hline
				$2\!\times\!12$ &$5.2718$&$0.30$&$5.6197$&$6.92$ & $5.200$ &  $1.07$  & $5.1798$ &  $1.45$    &  $4.3494$   &  $17.25$ &  $5.1128$   &  $2.72$    \\
				$4\!\times\!24$ &$5.2381$&$0.34$&$5.2665$&$0.20$ & $5.217$ &  $0.74$  & $5.2162$ &  $0.76$    &  $4.9987$   &  $4.90$  &  $5.2289$   &  $0.52$    \\
				$8\!\times\!48$ &$5.2423$&$0.26$& -      &  -    & $5.238$ &  $0.34$  & $5.2409$ &  $0.29$    &  $5.2334$   &  $0.43$  &  $5.2500$   &  $0.11$    \\
				$16\!\times\!96$&$5.2444$&$0.22$& -      &  -    & -       &    -     & $5.2479$ &  $0.15$    &  $5.2499$   &  $0.12$  &  $5.2505$   &  $0.10$    \\
				\hline
				\multicolumn{13}{c}{Reference solution: $u_Y(A)=5.256\,\text{m}$ (\cite{TwistedBelytschko})} \\ 
				\hline
			\end{tabularx} \nonumber
			\caption{Deflection $u_Y(A)$ for the twisted beam problem with $h=0.0032\,\text{m},~F_Y=0.001\,\text{N},~F_Z=0$}
		\end{table}
	\end{center}
	
	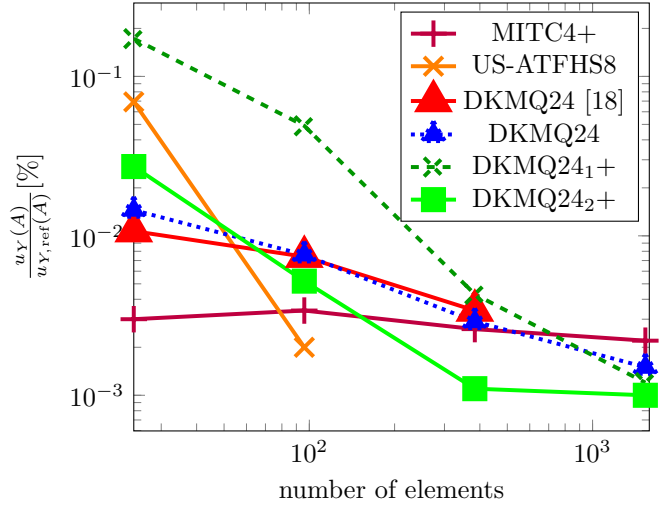
\begin{figure}
		\centering
		\begin{tikzpicture} 
		\begin{loglogaxis}[ 
		xlabel=\text{number of elements}, 
		ylabel= $\frac{u_Y(A)}{u_{Y,\rf}(A)}$\text{[\%]},
		x tick label style={
			/pgf/number format/1000 sep={}  
		},
		y tick label style={
			/pgf/number format/fixed,       
			/pgf/number format/precision=5  
		},
		scaled y ticks=false, 
		xmin= 24,
		xmax = 1586,
		legend style={at={(0.98,0.98)},anchor=north east}
		]  
		\addplot[line width=1.5pt, color=purple, mark=+, mark size=5pt] coordinates { 
			(24,    0.30E-2) 
			(96,    0.34E-2) 
			(384,   0.26E-2) 
			(1536,  0.22E-2)
		};    
		\addplot[line width=1.5pt, color=orange, mark=x, mark size=5pt] coordinates { 
			(24,    6.92E-2) 
			(96,    0.20E-2) 
		};    
		\addplot[line width=1.5pt, color=red, mark=triangle*, mark size=7pt] coordinates { 
			(24,    1.07E-2) 
			(96,    0.74E-2) 
			(384,   0.34E-2) 
		};
		\addplot[line width=1.5pt, color=blue, mark=triangle*, mark size=5pt, dotted] coordinates { 
			(24,    1.45E-2) 
			(96,    0.76E-2) 
			(384,   0.29E-2) 
			(1536,  0.15E-2)
		}; 
		\addplot[line width=1.5pt, color=col1, mark=x, mark size=5pt, dashed] coordinates { 
			(24,    17.25E-2) 
			(96,     4.90E-2) 
			(384,    0.43E-2) 
			(1536,   0.12E-2)
		};
		\addplot[line width=1.5pt, color=green, mark=square*, mark size=4pt] coordinates {    
			(24,    2.72E-2) 
			(96,    0.52E-2) 
			(384,   0.11E-2) 
			(1536,  0.10E-2)
		};
		\legend{MITC4+, US-ATFHS8, DKMQ24 \cite{Katili2015}, DKMQ24, DKMQ24$_1$+, DKMQ24$_2$+}
		\end{loglogaxis} 
		\end{tikzpicture} 
		\caption{Relative error of deflection $u_Y(A)$ for the twisted beam problem with $h=0.0032\,\text{m},~F_Y=0.001\,\text{N},~F_Z=0$}
		\label{fig:BENCH6a}
	\end{figure}	
	
	
	\begin{center}
		\begin{table}
			\centering
			\begin{tabularx}{160mm}{ c *{12}{Y} }
				\hline
				Mesh  & \multicolumn{2}{c}{MITC4+} & \multicolumn{2}{c}{US-ATFHS8}  &  \multicolumn{2}{c}{DKMQ24 \cite{Katili2015}}  & \multicolumn{2}{c}{DKMQ24} & \multicolumn{2}{c}{DKMQ24$_1$+} & \multicolumn{2}{c}{DKMQ24$_2$+} \rule[-2mm]{0mm}{6mm}  \\
				&    $u_Z(A)$   & err. &    $u_Z(A)$   & err. &    $u_Z(A)$   & err. &    $u_Z(A)$   & err.&    $u_Z(A)$   & err.&    $u_Z(A)$   & err. \\
				& [m] &  [\%]    & [m] &  [\%]    & [m] &  [\%]  & [m] &  [\%]  & [m] &  [\%]  & [m] &  [\%]    \\
				\hline
				\hline
				$2\!\times\!12$   &  $1.2777$ & $1.26$  &  $1.3762$ &  $6.35$  &  $1.274$    &  $1.55$  &  $1.2690$   &  $1.93$    &  $1.2644$   &  $2.28$ &  $1.2729$   &  $1.63$    \\
				$4\!\times\!24$   &  $1.2874$ & $0.51$  &  $1.2953$ &  $0.10$  &  $1.285$    &  $0.70$  &  $1.2859$   &  $0.62$    &  $1.2877$   &  $0.49$ &  $1.2898$   &  $0.32$    \\
				$8\!\times\!48$   &  $1.2908$ & $0.25$  &  -        &  -       &  $1.290$    &  $0.31$  &  $1.2912$   &  $0.22$    &  $1.2925$   &  $0.11$ &  $1.2926$   &  $0.11$    \\
				$16\!\times\!96$  &  $1.2917$ & $0.18$  &  -        &  -       &  -          &    -     &  $1.2926$   &  $0.11$    &  $1.2930$   &  $0.08$ &  $1.2930$   &  $0.08$    \\
				\hline
				\multicolumn{13}{c}{Reference solution: $u_Z(A)=1.294\,\text{m}$ (\cite{TwistedBelytschko})} \\ 
				\hline
			\end{tabularx} \nonumber
			\caption{Deflection $u_Z(A)$ for the twisted beam problem with $h=0.0032\,\text{m},~F_Y=0,~F_Z=0.001\,\text{N}$}
		\end{table}
	\end{center}
	
	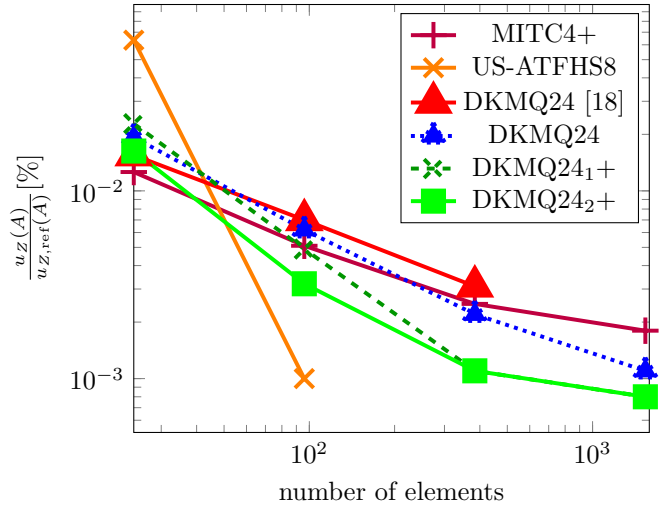
\begin{figure}
		\centering
		\begin{tikzpicture} 
		\begin{loglogaxis}[ 
		xlabel=\text{number of elements}, 
		ylabel= $\frac{u_Z(A)}{u_{Z,\rf}(A)}$\text{[\%]},
		x tick label style={
			/pgf/number format/1000 sep={}  
		},
		y tick label style={
			/pgf/number format/fixed,       
			/pgf/number format/precision=5  
		},
		scaled y ticks=false, 
		xmin= 24,
		xmax = 1586,
		legend style={at={(0.98,0.98)},anchor=north east}
		]  
		\addplot[line width=1.5pt, color=purple, mark=+, mark size=5pt] coordinates { 
			(24,    1.26E-2) 
			(96,    0.51E-2) 
			(384,   0.25E-2) 
			(1536,  0.18E-2)
		};    
		\addplot[line width=1.5pt, color=orange, mark=x, mark size=5pt] coordinates { 
			(24,    6.35E-2) 
			(96,    0.10E-2) 
		};    
		\addplot[line width=1.5pt, color=red, mark=triangle*, mark size=7pt] coordinates { 
			(24,    1.55E-2) 
			(96,    0.70E-2) 
			(384,   0.31E-2) 
		};
		\addplot[line width=1.5pt, color=blue, mark=triangle*, mark size=5pt, dotted] coordinates { 
			(24,    1.93E-2) 
			(96,    0.62E-2) 
			(384,   0.22E-2) 
			(1536,  0.11E-2)
		}; 
		\addplot[line width=1.5pt, color=col1, mark=x, mark size=5pt, dashed] coordinates { 
			(24,     2.28E-2) 
			(96,     0.49E-2) 
			(384,    0.11E-2) 
			(1536,   0.08E-2)
		};
		\addplot[line width=1.5pt, color=green, mark=square*, mark size=4pt] coordinates {    
			(24,    1.63E-2) 
			(96,    0.32E-2) 
			(384,   0.11E-2) 
			(1536,  0.08E-2)
		};
		\legend{MITC4+, US-ATFHS8, DKMQ24 \cite{Katili2015}, DKMQ24, DKMQ24$_1$+, DKMQ24$_2$+}
		\end{loglogaxis} 
		\end{tikzpicture} 
		\caption{Relative error of deflection $u_Z(A)$ for the twisted beam problem with $h=0.0032\,\text{m},~F_Y=0,~F_Z=0.001\,\text{N}$}
		\label{fig:BENCH6b}
	\end{figure}
	

	\begin{center}
		\begin{table}
			\centering
			\begin{tabularx}{160mm}{ c *{10}{Y} }
				\hline
				Mesh  & \multicolumn{2}{c}{MITC4+}  &  \multicolumn{2}{c}{DKMQ24 \cite{Katili2015}}  & \multicolumn{2}{c}{DKMQ24} & \multicolumn{2}{c}{DKMQ24$_1$+} & \multicolumn{2}{c}{DKMQ24$_2$+} \rule[-2mm]{0mm}{6mm}  \\
				&    $u_Y(A)$   & err. &    $u_Y(A)$   & err. &    $u_Y(A)$   & err.&    $u_Y(A)$   & err.&    $u_Y(A)$   & err. \\
				& [m] &  [\%]    & [m] &  [\%]  & [m] &  [\%]  & [m] &  [\%]  & [m] &  [\%]    \\
				\hline
				\hline
				$2\!\times\!12$   &  $5.3925$ & $0.58$ &  $5.393$    &  $0.57$  &    $5.3762$   &  $0.88$    &  $5.4403$   &  $0.30$ &  $5.4411$   &  $0.32$    \\
				$4\!\times\!24$   &  $5.4023$ & $0.40$ &  $5.403$    &  $0.39$  &    $5.4013$   &  $0.42$    &  $5.4197$   &  $0.08$ &  $5.4198$   &  $0.08$    \\
				$8\!\times\!48$   &  $5.4066$ & $0.32$ &  $5.410$    &  $0.24$  &    $5.4116$   &  $0.23$    &  $5.4167$   &  $0.13$ &  $5.4167$   &  $0.13$    \\
				$16\!\times\!96$  &  $5.4083$ & $0.29$ &  -          &    -     &    $5.4151$   &  $0.16$    &  $5.4165$   &  $0.14$ &  $5.4165$   &  $0.14$    \\
				\hline
				\multicolumn{11}{c}{Reference solution: $u_Y(A)=5.424\,\text{m}$ (\cite{TwistedBelytschko})} \\ 
				\hline
			\end{tabularx} \nonumber
			\caption{Deflection $u_Y(A)$ for the twisted beam problem with $h=0.32\,\text{m},~F_Y=1000\,\text{N},~F_Z=0$}
		\end{table}
	\end{center}
	
	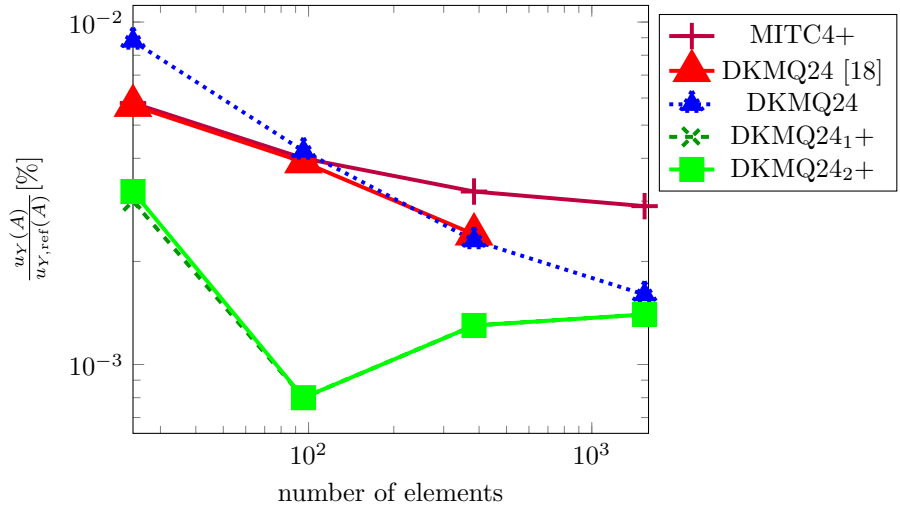
\begin{figure}
		\centering
		\begin{tikzpicture} 
		\begin{loglogaxis}[ 
		xlabel=\text{number of elements}, 
		ylabel= $\frac{u_Y(A)}{u_{Y,\rf}(A)}$\text{[\%]},
		x tick label style={
			/pgf/number format/1000 sep={}  
		},
		y tick label style={
			/pgf/number format/fixed,       
			/pgf/number format/precision=5  
		},
		scaled y ticks=false, 
		xmin= 24,
		xmax = 1586,
		legend style={at={(1.02,0.98)},anchor=north west}
		]  
		\addplot[line width=1.5pt, color=purple, mark=+, mark size=5pt] coordinates { 
			(24,    0.58E-2) 
			(96,    0.40E-2) 
			(384,   0.32E-2) 
			(1536,  0.29E-2)
		};      
		\addplot[line width=1.5pt, color=red, mark=triangle*, mark size=7pt] coordinates { 
			(24,    0.57E-2) 
			(96,    0.39E-2) 
			(384,   0.24E-2) 
		};
		\addplot[line width=1.5pt, color=blue, mark=triangle*, mark size=5pt, dotted] coordinates { 
			(24,    0.88E-2) 
			(96,    0.42E-2) 
			(384,   0.23E-2) 
			(1536,  0.16E-2)
		}; 
		\addplot[line width=1.5pt, color=col1, mark=x, mark size=5pt, dashed] coordinates { 
			(24,     0.30E-2) 
			(96,     0.08E-2) 
			(384,    0.13E-2) 
			(1536,   0.14E-2)
		};
		\addplot[line width=1.5pt, color=green, mark=square*, mark size=4pt] coordinates {    
			(24,    0.32E-2) 
			(96,    0.08E-2) 
			(384,   0.13E-2) 
			(1536,  0.14E-2)
		};
		\legend{MITC4+, DKMQ24 \cite{Katili2015}, DKMQ24, DKMQ24$_1$+, DKMQ24$_2$+}
		\end{loglogaxis} 
		\end{tikzpicture} 
		\caption{Relative error of deflection $u_Y(A)$ for the twisted beam problem with $h=0.32\,\text{m},~F_Y=1000\,\text{N},~F_Z=0$}
		\label{fig:BENCH6c}
	\end{figure}
	
	
	\begin{center}
		\begin{table}
			\centering
			\begin{tabularx}{160mm}{ c *{10}{Y} }
				\hline
				Mesh  & \multicolumn{2}{c}{MITC4+}  &  \multicolumn{2}{c}{DKMQ24 \cite{Katili2015}}  & \multicolumn{2}{c}{DKMQ24} & \multicolumn{2}{c}{DKMQ24$_1$+} & \multicolumn{2}{c}{DKMQ24$_2$+} \rule[-2mm]{0mm}{6mm}  \\
				&    $u_Z(A)$   & err. &    $u_Z(A)$   & err. &    $u_Z(A)$   & err.&    $u_Z(A)$   & err.&    $u_Z(A)$   & err. \\
				& [m] &  [\%]    & [m] &  [\%]  & [m] &  [\%]  & [m] &  [\%]  & [m] &  [\%]    \\
				\hline
				\hline
				$2\!\times\!12$   &  $1.7245$ & $1.68$ &  $1.624$    &  $7.41$  &    $1.6198$   &  $7.65$    &  $1.7486$   &  $0.31$ &  $1.7487$   &  $0.30$    \\
				$4\!\times\!24$   &  $1.7428$ & $0.64$ &  $1.711$    &  $2.45$  &    $1.7109$   &  $2.46$    &  $1.7503$   &  $0.21$ &  $1.7503$   &  $0.21$    \\
				$8\!\times\!48$   &  $1.7479$ & $0.35$ &  $1.740$    &  $0.80$  &    $1.7409$   &  $0.75$    &  $1.7514$   &  $0.15$ &  $1.7514$   &  $0.15$    \\
				$16\!\times\!96$  &  $1.7493$ & $0.27$ &  -          &    -     &    $1.7491$   &  $0.28$    &  $1.7518$   &  $0.13$ &  $1.7518$   &  $0.13$    \\
				\hline
				\multicolumn{11}{c}{Reference solution: $u_Z(A)=1.754\,\text{m}$ (\cite{TwistedBelytschko})} \\ 
				\hline
			\end{tabularx} \nonumber
			\caption{Deflection $u_Z(A)$ for the twisted beam problem with $h=0.32\,\text{m},~F_Y=0,~F_Z=1000\,\text{N}$}
		\end{table}
	\end{center}

	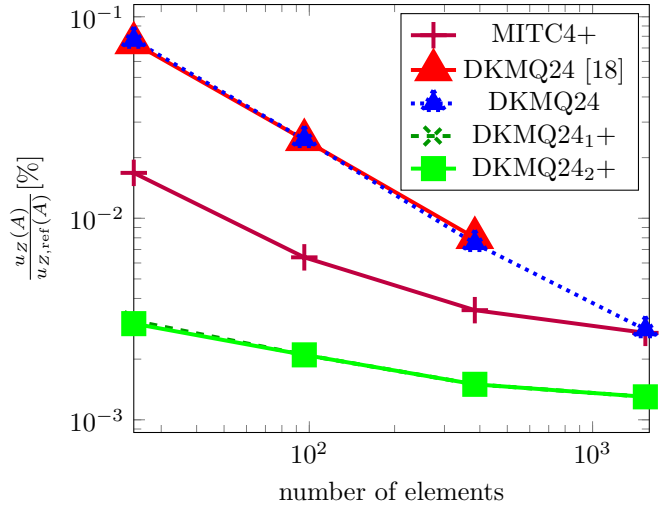
\begin{figure}
		\centering
		\begin{tikzpicture} 
		\begin{loglogaxis}[ 
		xlabel=\text{number of elements}, 
		ylabel= $\frac{u_Z(A)}{u_{Z,\rf}(A)}$\text{[\%]},
		x tick label style={
			/pgf/number format/1000 sep={}  
		},
		y tick label style={
			/pgf/number format/fixed,       
			/pgf/number format/precision=5  
		},
		scaled y ticks=false, 
		xmin= 24,
		xmax = 1586,
		legend style={at={(0.98,0.98)},anchor=north east}
		]  
		\addplot[line width=1.5pt, color=purple, mark=+, mark size=5pt] coordinates { 
			(24,    1.68E-2) 
			(96,    0.64E-2) 
			(384,   0.35E-2) 
			(1536,  0.27E-2)
		};      
		\addplot[line width=1.5pt, color=red, mark=triangle*, mark size=7pt] coordinates { 
			(24,    7.41E-2) 
			(96,    2.45E-2) 
			(384,   0.80E-2) 
		};
		\addplot[line width=1.5pt, color=blue, mark=triangle*, mark size=5pt, dotted] coordinates { 
			(24,    7.65E-2) 
			(96,    2.46E-2) 
			(384,   0.75E-2) 
			(1536,  0.28E-2)
		}; 
		\addplot[line width=1.5pt, color=col1, mark=x, mark size=5pt, dashed] coordinates { 
			(24,     0.31E-2) 
			(96,     0.21E-2) 
			(384,    0.15E-2) 
			(1536,   0.13E-2)
		};
		\addplot[line width=1.5pt, color=green, mark=square*, mark size=4pt] coordinates {    
			(24,    0.30E-2) 
			(96,    0.21E-2) 
			(384,   0.15E-2) 
			(1536,  0.13E-2)
		};
		\legend{MITC4+, DKMQ24 \cite{Katili2015}, DKMQ24, DKMQ24$_1$+, DKMQ24$_2$+}
		\end{loglogaxis} 
		\end{tikzpicture} 
		\caption{Relative error of deflection $u_Z(A)$ for the twisted beam problem with $h=0.32\,\text{m},~F_Y=0,~F_Z=1000\,\text{N}$}
		\label{fig:BENCH6d}
	\end{figure}
	
	\FloatBarrier
	
	\subsection{Raasch's hook}
	\label{b7}
	\begin{figure}[H]
		\centering
		\begin{overpic}[width=0.4\textwidth]{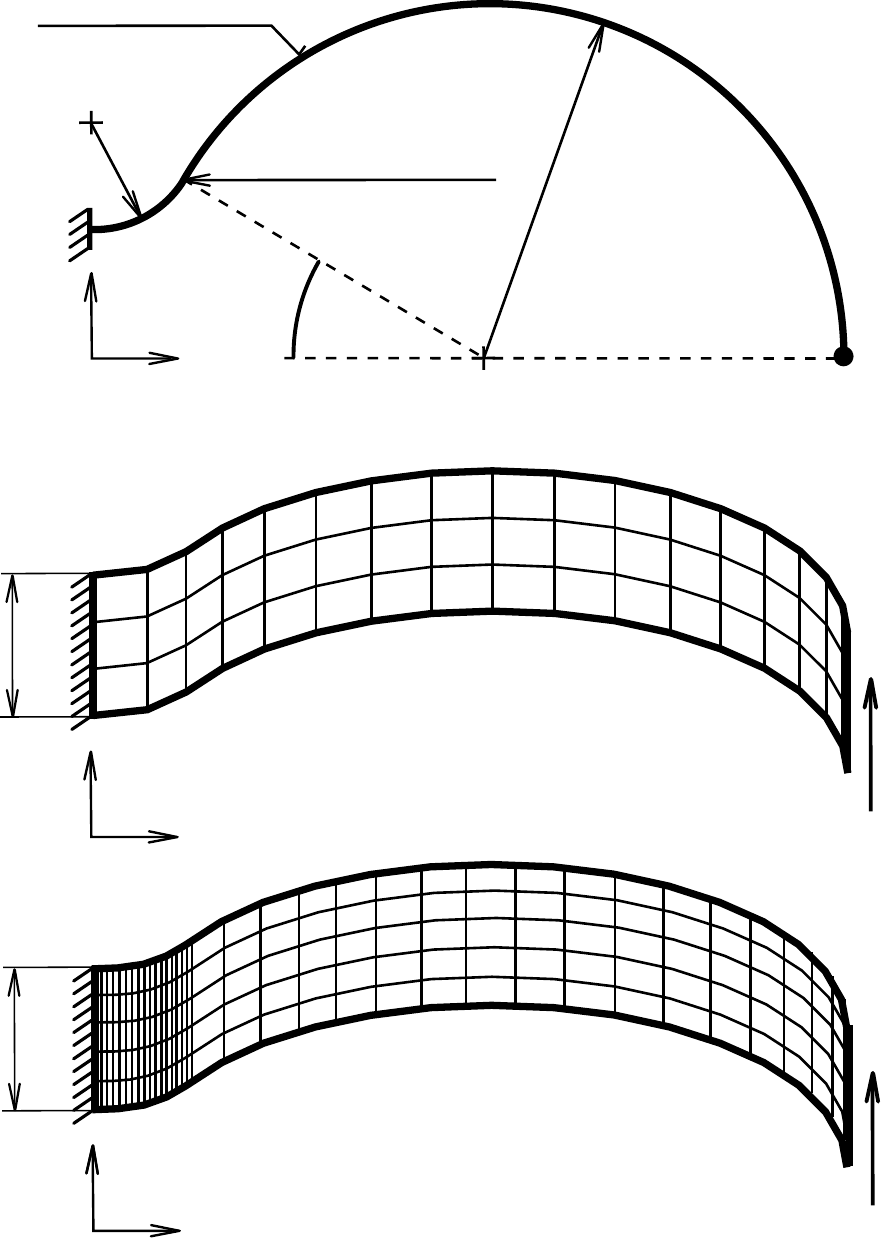} 
			\put(10,87.25){$R_1$}
			\put(19.5,86.5){$\text{inflect.point}$}
			\put(45,82){$R_2$}
			\put(-3.5,14.25){$w$}
			\put(-3.5,47){$w$}
			\put(26,72.5){$30^{\circ}$}		    
			\put(70,70){$A$}		    
			\put(3.25,99.1){$\text{thickness}~h$}		    
			\put(12.5,66){$X$}
			\put(12.5,27.75){$X$}
			\put(12.5,-4){$X$}
			\put(2.25,75){$Y$}
			\put(2,5){$Z$}
			\put(2,37){$Z$}
			\put(72,6){$F$}
			\put(72,38){$F$}
			\put(82,84.5){$R_1=14\,\text{m}$}
			\put(82,74.5){$R_2=46\,\text{m}$}
			\put(82,64.5){$w=20\,\text{m}$}
			\put(82,54.5){$h=2\,\text{m}$}
			\put(82,44.5){$E=3300\,\text{Pa}$}
			\put(82,34.5){$\nu=0.35$}	
			\put(82,24.5){$F=1\,\text{N}$}
		\end{overpic}
		\caption{Raasch's hook with mesh $17\times3$ (above) and $34\times5$ (below)}
		\label{RaaschsHook}
	\end{figure}
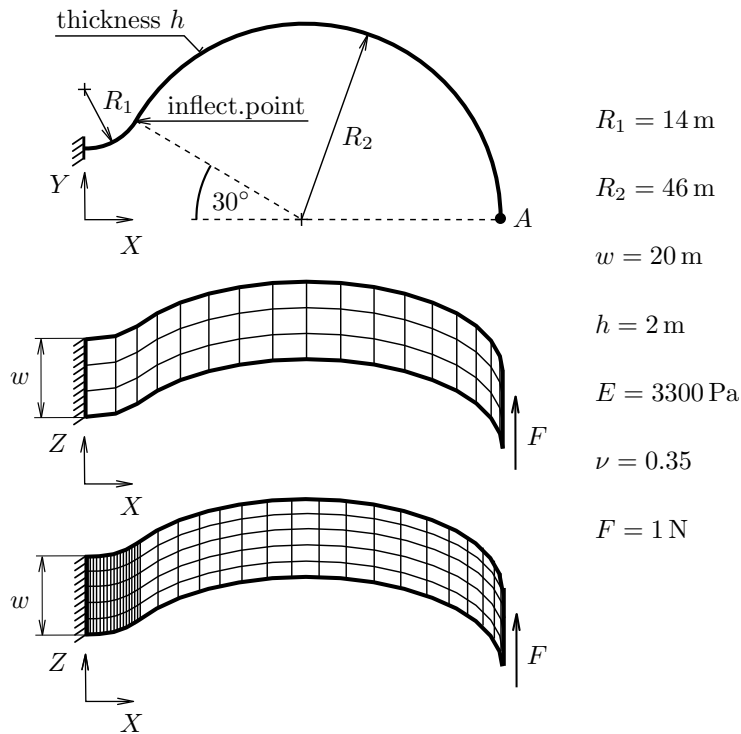
	\FloatBarrier
	
	Raasch's hook problem was first reported by Ingo Raasch of the BMW corporation in Germany (see Figure \ref{RaaschsHook}). It is a challenging problem because of the inherent coupling between three modes of deformation: bending, extension, and twist. Many shell elements with drilling dofs exhibit divergence of the tip deflection at this benchmark (\cite{Raasch}).
	Except of the coarsest mesh $17\times 3$, which has uniformly distributed elements, all other meshes we used have, in accordance with \cite{MITC4Plus3}, the same number of elements before and after the inflection point, see Figure \ref{RaaschsHook}.
	
	The results are shown in Table \ref{tableRaasch}. We note that a sufficient size of the penalty constant is needed in order to obtain convergence at this benchmark. 
	
	Our implementation of the DKMQ24 element exhibits large errors at this benchmark and, moreover, diverges to infinity when the mesh gets finer. An increase of the penalty constant $c$ does not yield expected convergence, moreover, an increase of the penalty constant negatively influences the error at other benchmarks. For further discussion on this topic see Section 4. The DKMQ24$_2$+ element provides lower errors when compared to the DKMQ24$_1$+ element.  Moreover, the DKMQ24$_2$+ element provides a fast convergence to a fixed value, which is however slightly higher then the reference value (by $0.3$\%). The DKMQ24$_2$+ element and the MITC4+ element provide similar errors, however, we note that the MITC4+ element is tested with a slightly different setting ($\nu=0.3$ instead of $\nu=0.35$) and, therefore, converging to a different target value ($4.82482\,\text{m}$ instead of $5.027\,\text{m}$). 
	
	
	
	
	\begin{center}
		\begin{table}
			\centering
			\begin{tabularx}{120mm}{ c *{6}{Y} }
				\hline
				Mesh           &  \multicolumn{2}{c}{DKMQ24}  & \multicolumn{2}{c}{DKMQ24$_1$+} & \multicolumn{2}{c}{DKMQ24$_2$+} \rule[-2mm]{0mm}{6mm}  \\
				&    $u_Z(A)$   & err. &    $u_Z(A)$   & err. &    $u_Z(A)$   & err. \\
				& [m] &  [\%]    & [m] &  [\%]    & [m] &  [\%]    \\
				\hline
				\hline
				$17\!\times\!3$     &  $2.5852$ & $48.57$    &    $4.6973$   &  $6.56$    &  $4.6997$     &  $6.51$     \\
				$34\!\times\!5$     &  $2.3730$ & $52.79$    &    $4.7895$   &  $4.72$    &  $4.8070$     &  $4.38$     \\
				$68\!\times\!10$    &  $3.5430$ & $29.52$    &    $4.9537$   &  $1.46$    &  $4.9629$     &  $1.27$     \\
				$136\!\times\!20$   &  $5.9169$ & $17.70$    &    $5.0712$   &  $0.88$    &  $5.0435$     &  $0.33$     \\
				$272\!\times\!40$   &  $8.1284$ & $61.70$    &    $5.1054$   &  $1.56$    &  $5.0496$     &  $0.45$     \\
				$544\!\times\!80$   &  $9.1194$ & $81.41$    &    $5.1146$   &  $1.74$    &  $5.0436$     &  $0.33$     \\
				\hline
				\multicolumn{7}{c}{Reference solution: $u_Z(A)=5.027\,\text{m}$ (\cite{Kemp})} \\ 
				\hline
			\end{tabularx} \nonumber
			\caption{Deflection $u_Z(A)$ for Raasch's hook}
			\label{tableRaasch} 
		\end{table}
	\end{center}
	
	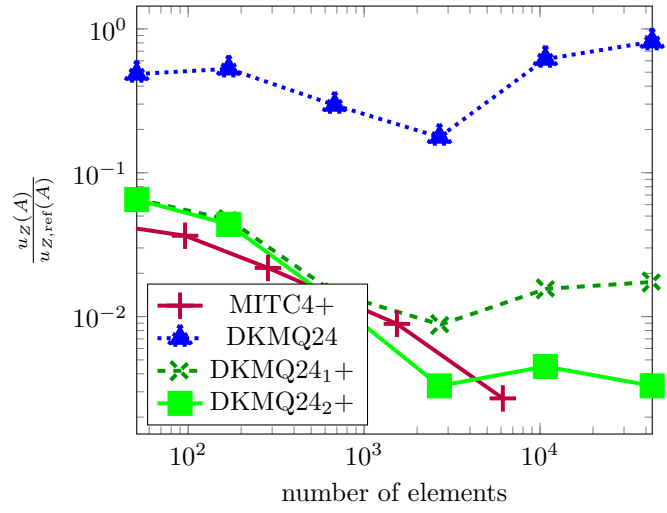
\begin{figure}
		\centering
		\begin{tikzpicture} 
		\begin{loglogaxis}[ 
		xlabel=\text{number of elements}, 
		ylabel= $\frac{u_Z(A)}{u_{Z,\rf}(A)}$,  
		x tick label style={
			/pgf/number format/1000 sep={}  
		},
		y tick label style={
			/pgf/number format/1000 sep={},
			/pgf/number format/fixed,       
			/pgf/number format/precision=5  
		},
		scaled y ticks=false, 
		xmin= 51,
		xmax = 43410, 
		legend style={at={(0.02,0.02)},anchor=south west}
		] 
		\addplot[line width=1.5pt, color=purple, mark=+, mark size=5pt] coordinates { 
			(24,  4.69E-2)
			(96,  3.65E-2) 
			(284, 2.18E-2) 
			(1536, 0.89E-2) 
			(6144, 0.27E-2) 
		}; 
		\addplot[line width=1.5pt, color=blue, mark=triangle*, mark size=5pt, dotted] coordinates { 
			(51,  48.57E-2)
			(170, 52.79E-2) 
			(680,  29.52E-2) 
			(2680, 17.70E-2)
			(10720,61.70E-2) 
			(43360, 81.41E-2)  
		};
		\addplot[line width=1.5pt, color=col1, mark=x, mark size=5pt, dashed] coordinates { 
			(51,  6.56E-2)
			(170, 4.72E-2) 
			(680,  1.46E-2) 
			(2720, 0.88E-2) 
			(10840, 1.56E-2) 
			(43360, 1.74E-2) 
		};
		\addplot[line width=1.5pt, color=green, mark=square*, mark size=4pt] coordinates { 
			(51,    6.51E-2)
			(170,   4.38E-2) 
			(680,   1.27E-2) 
			(2680,  0.33E-2) 		
			(10720, 0.45E-2)
			(43360, 0.33E-2)  
		};
		\legend{MITC4+, DKMQ24, DKMQ24$_1$+, DKMQ24$_2$+} 
		\end{loglogaxis} 
		\end{tikzpicture} 
		\caption{Relative error of deflection $u_Z(A)$ for Raasch's hook. The MITC4+ element results, taken from \cite{MITC4Plus3}, are valid for the different seting $\nu=0.3$ having the different target value $4.82482\,\text{m}$.}
		\label{fig:BENCH9a}
	\end{figure}
	
	\FloatBarrier

	\subsection{Partly clamped hyperbolic paraboloid shell}
	\label{b8}
	\begin{figure}[H]
		\centering
		\begin{overpic}[width=0.55\textwidth]{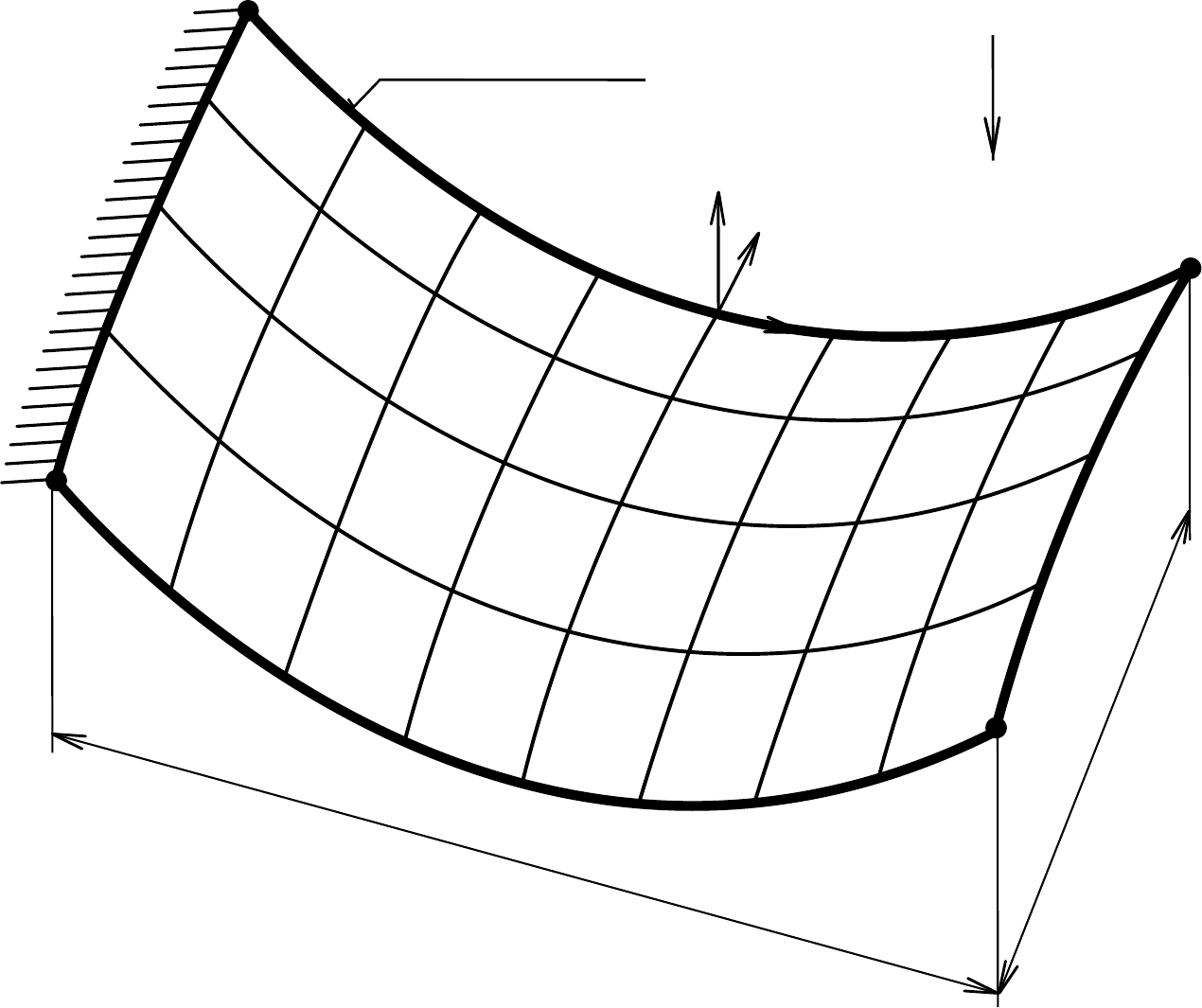} 
			\put(66,57.25){$X$}
			\put(64,63.75){$Y$}
			\put(55.5,65.75){$Z$}
			\put(0,38.75){$A$}
			\put(84.5,21){$B$}
			\put(101,62){$C$}
			\put(21,84.75){$D$}
			\put(85,75.25){$g$}
			\put(32.5,78.25){$\text{thickness}~h$}		    
			\put(42,13.5){$L$}
			\put(90,29){$\frac{L}{2}$}
			\put(110,63){$L=1\,\text{m}$}
			\put(110,57){$h=0.001/~0.0001\,\text{m}$}
			\put(110,51){$\rho=360/\,3.6\,\text{kg}\,\text{m}^{-3}$}
			\put(110,45){$g=1\,\text{m}\,\text{s}^{-2}$}
			\put(110,39){$E=200\,\text{GPa}$}
			\put(110,33){$\nu=0.3$}	
			\put(110,27){$p=1\,\text{Pa}$}
		\end{overpic}
		\caption{Partly clamped hyperbolic paraboloid shell with uniform mesh $4\times8$}
		\label{PartiallyClampedHyperbolicParaboloid}
	\end{figure}
	
	\begin{figure}[H]
		\centering
		\begin{overpic}[width=0.4\textwidth]{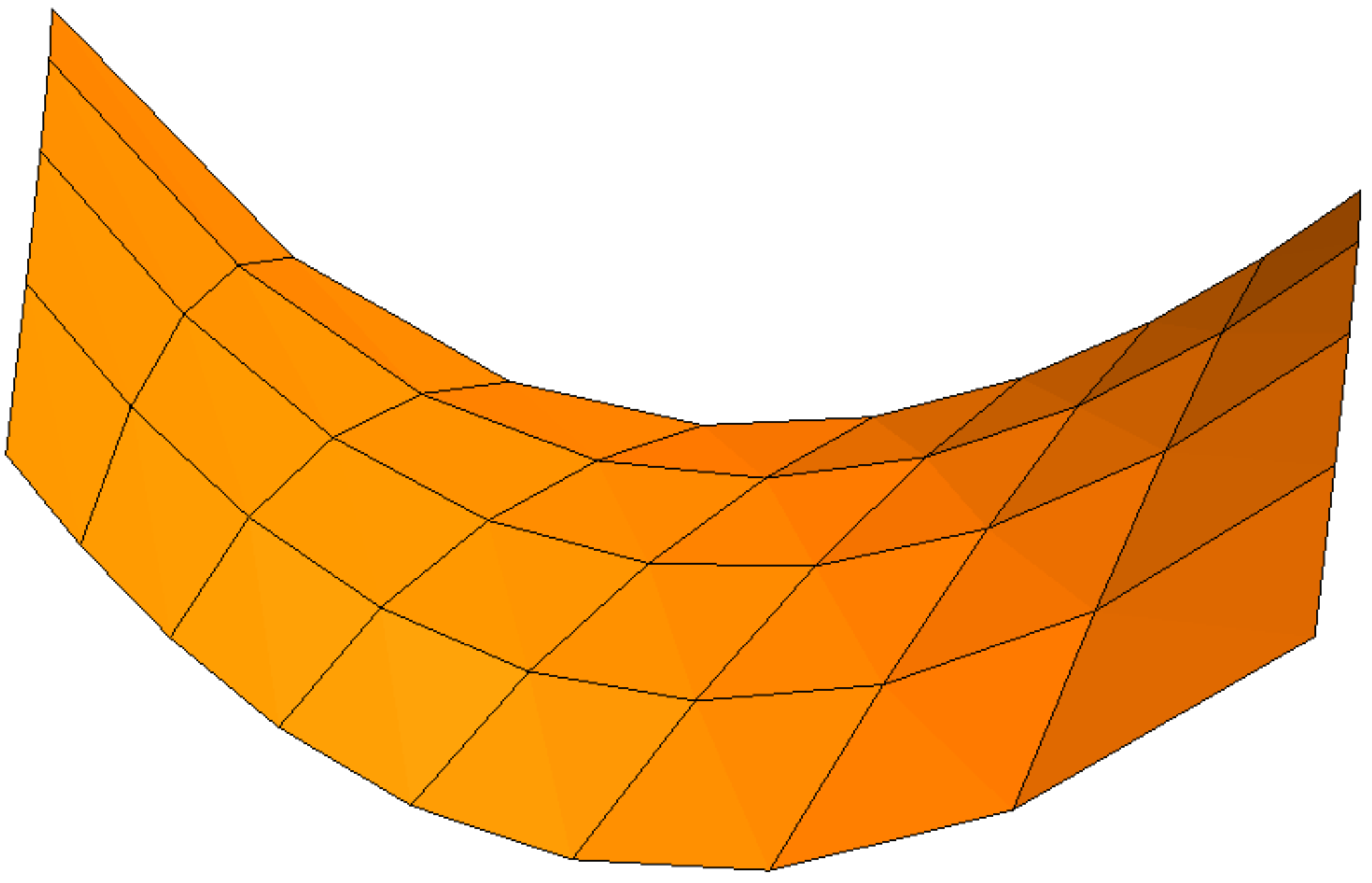} 
		\end{overpic}
		\caption{Partly clamped hyperbolic paraboloid shell -- distorted mesh $4\times8$ with the aspect ratio 3 (in the XY plane) on sides AB and CD}
		\label{PartiallyClampedHyperbolicParaboloidDistMesh}
	\end{figure}
	\FloatBarrier
	
	The partly clamped hyperbolic paraboloid shell benchmark, taken from \cite{MITC4Plus3}, whose geometry is described by equation 
	\begin{align}
	Z=X^2-Y^2,~X\in\left[-\frac{1}{2},\frac{1}{2}\right],~Y\in\left[-\frac{1}{2},0\right]
	\end{align}
	is loaded by self-weight (see Figure \ref{PartiallyClampedHyperbolicParaboloid}). On side CD, symmetry conditions are prescribed ($u_Y=\varphi_X=\varphi_Z=0$), the side DA is fully clamped ($u_X=u_Y=u_Z=\varphi_X=\varphi_Y=\varphi_Z=0$) and other sides are free. Two different shell thicknesses are considered: $h/L=0.001$ and $h/L=0.0001$. The corresponding load magnitudes are scaled by using different densities $\rho=360\,\text{kg}\,\text{m}^{-3}$ and  $\rho=3.6\,\text{kg}\,\text{m}^{-3}$, respectively. The reference solution is taken from \cite{MITC4Plus3}.
	
	Let us first comment the results on uniform meshes. For the case $h=0.001\,\text{m}$ the modified elements DKMQ24$_1$+ and DKMQ24$_2$+ are superior to all other elements, except of the coarsest mesh $4\times8$. For the case $h=0.0001\,\text{m}$ the modified elements DKMQ24$_1$+ and DKMQ24$_2$+ provide smaller errors when compared to the MITC4+ element only at finer meshes ($32\times 64$ and $64 \times 128$). Let us now comment the results on distorted meshes. For the case $h=0.001\,\text{m}$ the modified elements DKMQ24$_1$+ and DKMQ24$_2$+ when compared to the MITC4+ element, are superior at finer meshes ($32\times64$ and $64\times128$), comparable at the middle-sized mesh $16\times32$ and worse on coarser meshes. For the case $h=0.0001\,\text{m}$  the MITC4+ elements shows lowest errors at all mesh sizes, which demonstrated its good ability to resist to a mesh distorsion. In all the tested cases we see that the modified elements DKMQ24$_1$+ and DKMQ24$_2$+ provide significantly lower errors when compared to the DKMQ24 element, which demonstrated better behaviour of both elements at membrane-dominated problems.

	\begin{center}
		\begin{table}
			\centering
			\begin{tabularx}{160mm}{ c *{10}{Y} }
				\hline
				Mesh  
				& \multicolumn{2}{c}{MITC4+} 
				& \multicolumn{2}{c}{US-ATFHS8}    
				& \multicolumn{2}{c}{DKMQ24}  
				& \multicolumn{2}{c}{DKMQ24$_1$+} 
				& \multicolumn{2}{c}{DKMQ24$_2$+}  \rule[-2mm]{0mm}{6mm}  \\
				& $u_Z(C)$ & err. & $u_Z(C)$ & err. & $u_Z(C)$ & err. & $u_Z(C)$ & err. & $u_Z(C)$ & err. \\
				& [mm] &  [\%]    & [mm] &  [\%]    & [mm] &  [\%]  & [mm] &  [\%]    & [mm] &  [\%]    \\
				\hline
				\hline
				$4\!\times\!8$   & $-0.27433$ & $4.68$ & $-0.28875$ & $0.33$ &  $-0.14864$ & $48.35$   & $-0.17433$ & $39.43$ &   $-0.18961$ & $34.12$ \\
				$8\!\times\!16$   & $-0.27499$ & $4.45$ & $-0.27896$ & $3.07$ &  $-0.22297$ & $22.53$   & $-0.28434$ & $1.20$  &   $-0.28502$ & $0.97$  \\
				$16\!\times\!32$ & $-0.28095$ & $2.38$ & $-0.28230$ & $1.91$ &  $-0.25861$ & $10.14$   & $-0.28877$ & $0.34$  &   $-0.28883$ & $0.36$  \\
				$32\!\times\!64$ & $-0.28509$ & $0.94$ & -          &  -     &  $-0.27085$ & $5.89$    & $-0.28799$ & $0.07$  &   $-0.28800$ & $0.07$  \\
				$64\!\times\!128$ & $-0.28691$ & $0.31$ & -          &  -     &  $-0.27472$ & $4.55$    & $-0.28782$ & $0.01$  &   $-0.28783$ & $0.01$  \\
				\hline
				\multicolumn{11}{c}{Reference solution: $u_Z(C)=-0.28780\,\text{mm}$ (\cite{MITC4Plus3})} \\ 
				\hline
			\end{tabularx} \nonumber
			\caption{Deflection $u_Z(C)$ for the partly clamped hyperbolic paraboloid shell with $h=0.001\,\text{m},~\rho=360\,\text{kg}\,\text{m}^{-3}$, uniform mesh}
		\end{table}
	\end{center}

	\begin{figure}
		\centering
		\begin{tikzpicture} 
		\begin{loglogaxis}[ 
		xlabel=\text{number of elements}, 
		ylabel= $\frac{u_Z(C)}{u_{Z,\rf}(C)}$,  
		x tick label style={
			/pgf/number format/1000 sep={}  
		},
		y tick label style={
			/pgf/number format/1000 sep={},
			/pgf/number format/fixed,       
			/pgf/number format/precision=5  
		},
		scaled y ticks=false, 
		xmin= 32,
		xmax = 8242,
		legend style={at={(0.02,0.02)},anchor=south west}
		] 
		\addplot[line width=1.5pt, color=purple, mark=+, mark size=5pt] coordinates {
			(32,  4.68E-2) 
			(128,  4.45E-2) 
			(512, 2.38E-2) 
			(2048,0.94E-2)
			(8192,0.31E-2)
		};
		\addplot[line width=1.5pt, color=orange, mark=x, mark size=5pt] coordinates { 
			(32,  0.33E-2) 
			(128,  3.07E-2) 
			(512, 1.91E-2) 
		};  
		\addplot[line width=1.5pt, color=blue, mark=triangle*, mark size=5pt, dotted] coordinates { 
			(32,  48.35E-2) 
			(128,  22.53E-2) 
			(512, 10.14E-2) 
			(2048, 5.89E-2)
			(8192, 4.55E-2)
		};
		\addplot[line width=1.5pt, color=col1, mark=x, mark size=5pt, dashed] coordinates { 
			(32,  39.43E-2) 
			(128,  1.20E-2) 
			(512,  0.34E-2) 
			(2048, 0.07E-2)
			(8192, 0.01E-2)
		};
		\addplot[line width=1.5pt, color=green, mark=square*, mark size=4pt] coordinates { 
			(32,  34.12E-2) 
			(128,  0.97E-2) 
			(512, 0.36E-2) 
			(2048,0.07E-2)
			(8192,0.01E-2)
		};
		\legend{MITC4+, US-ATFHS8, DKMQ24, DKMQ24$_1$+, DKMQ24$_2$+} 
		\end{loglogaxis} 
		\end{tikzpicture} 
		\caption{Deflection $u_Z(C)$ for the partly clamped hyperbolic paraboloid shell with $h=0.001\,\text{m},~\rho=360\,\text{kg}\,\text{m}^{-3}$, uniform mesh} 
		\label{fig:BENCH8a}
	\end{figure}
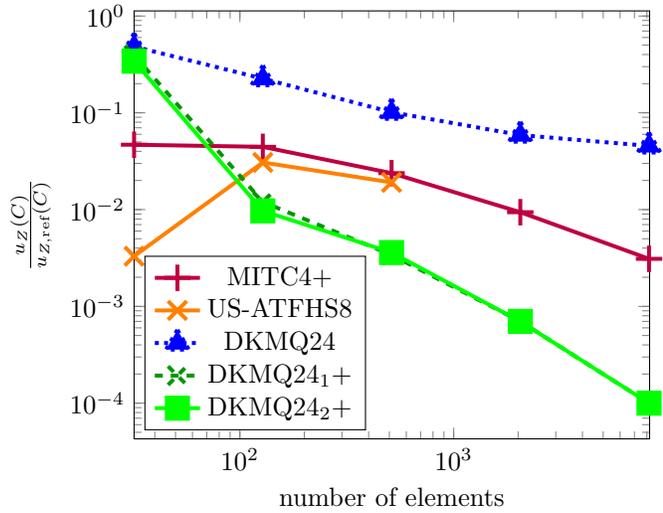
	
	\begin{center}
		\begin{table}
			\centering
			\begin{tabularx}{160mm}{ c *{8}{Y} }
				\hline
				Mesh  & \multicolumn{2}{c}{MITC4+}  & \multicolumn{2}{c}{DKMQ24} & \multicolumn{2}{c}{DKMQ24$_1$+} 
				& \multicolumn{2}{c}{DKMQ24$_2$+}  \rule[-2mm]{0mm}{6mm}  \\
				& $u_Z(C)$ & err. & $u_Z(C)$ & err. & $u_Z(C)$ & err. & $u_Z(C)$ & err. \\
				& [mm] &  [\%]    & [mm] &  [\%]    & [mm] &  [\%]  & [mm] &  [\%]      \\
				\hline
				\hline
				$4\!\times\!8$   & $-0.24333$ & $2.00$ & $-0.01083$ & $95.46$ &  $-0.00758$ & $96.82$ & $-0.02796$ & $88.28$  \\
				$8\!\times\!16$   & $-0.23443$ & $1.73$ & $-0.01316$ & $94.48$ &  $-0.09337$ & $60.86$ & $-0.11258$ & $52.81$  \\
				$16\!\times\!32$ & $-0.23324$ & $2.23$ & $-0.03860$ & $83.82$ &  $-0.22284$ &  $6.59$ & $-0.22400$ &  $6.10$  \\
				$32\!\times\!64$ & $-0.23541$ & $1.32$ & $-0.09919$ & $58.42$ &  $-0.23892$ &  $0.15$ & $-0.23900$ &  $0.18$  \\
				$64\!\times\!128$ & $-0.23730$ & $0.53$ & $-0.17178$ & $27.99$ &  $-0.23891$ &  $0.15$ & $-0.23892$ &  $0.15$  \\
				\hline
				\multicolumn{9}{c}{Reference solution: $u_Z(C)=-0.23856\,\text{mm}$ (\cite{MITC4Plus3})} \\ 
				\hline
			\end{tabularx} \nonumber
			\caption{Deflection $u_Z(C)$ for the partly clamped hyperbolic paraboloid shell with $h=0.0001\,\text{m},~\rho=3.6\,\text{kg}\,\text{m}^{-3}$, uniform mesh}
		\end{table}
	\end{center}

	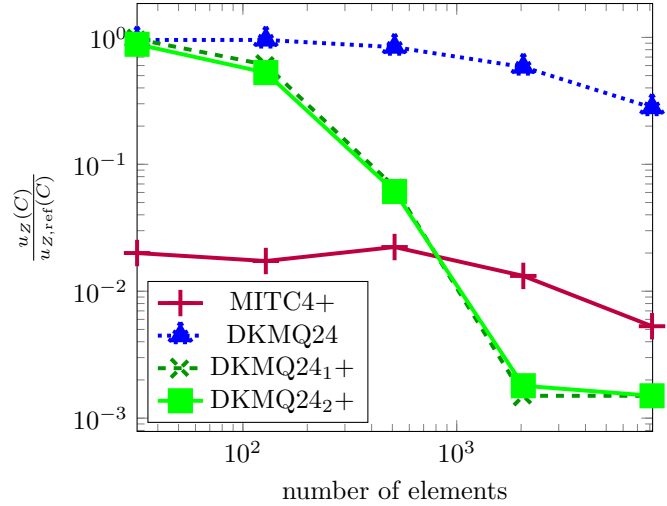
\begin{figure}
		\centering
		\begin{tikzpicture} 
		\begin{loglogaxis}[ 
		xlabel=\text{number of elements}, 
		ylabel= $\frac{u_Z(C)}{u_{Z,\rf}(C)}$,  
		x tick label style={
			/pgf/number format/1000 sep={}  
		},
		y tick label style={
			/pgf/number format/1000 sep={},
			/pgf/number format/fixed,       
			/pgf/number format/precision=5  
		},
		scaled y ticks=false, 
		xmin= 32,
		xmax = 8242,
		legend style={at={(0.02,0.02)},anchor=south west}
		] 
		\addplot[line width=1.5pt, color=purple, mark=+, mark size=5pt] coordinates {
			(32,  2.0E-2) 
			(128,  1.73E-2) 
			(512, 2.23E-2) 
			(2048,1.32E-2)
			(8192,0.53E-2)
		};  
		\addplot[line width=1.5pt, color=blue, mark=triangle*, mark size=5pt, dotted] coordinates { 
			(32,  95.46E-2) 
			(128,  95.48E-2) 
			(512, 83.82E-2) 
			(2048,58.42E-2)
			(8192,27.99E-2)
		};
		\addplot[line width=1.5pt, color=col1, mark=x, mark size=5pt, dashed] coordinates { 
			(32,  96.82E-2) 
			(128,  60.86E-2) 
			(512,  6.59E-2)  
			(2048, 0.15E-2)
			(8192, 0.15E-2)
		};
		\addplot[line width=1.5pt, color=green, mark=square*, mark size=4pt] coordinates { 
			(32,  88.28E-2) 
			(128,  52.81E-2) 
			(512, 6.10E-2) 
			(2048,0.18E-2)
			(8192,0.15E-2)
		};
		\legend{MITC4+, DKMQ24, DKMQ24$_1$+, DKMQ24$_2$+}
		\end{loglogaxis} 
		\end{tikzpicture} 
		\caption{Relative error of deflection $u_Z(C)$ for the partly clamped hyperbolic paraboloid shell with $h=0.0001\,\text{m},~\rho=3.6\,\text{kg}\,\text{m}^{-3}$, uniform mesh}
		\label{fig:BENCH8b}
	\end{figure}
	
	
	\begin{center}
		\begin{table}
			\centering
			\begin{tabularx}{160mm}{ c *{8}{Y} }
				\hline
				Mesh  
				& \multicolumn{2}{c}{MITC4+} 
				& \multicolumn{2}{c}{DKMQ24}  
				& \multicolumn{2}{c}{DKMQ24$_1$+} 
				& \multicolumn{2}{c}{DKMQ24$_2$+}  \rule[-2mm]{0mm}{6mm}  \\ 
				& $u_Z(C)$ & err. & $u_Z(C)$ & err. & $u_Z(C)$ & err. & $u_Z(C)$ & err.  \\
				& [mm] &  [\%]    & [mm] &  [\%]    & [mm] &  [\%]  & [mm] &  [\%]       \\
				\hline
				\hline
				$4\!\times\!8$   & $-0.21755$ & $24.41$ &    $-0.06171$ & $78.56$   & $-0.07755$ & $73.05$ &   $-0.07956$ & $72.35$ \\
				$8\!\times\!16$  &  $-0.27692$ & $3.78$ &    $-0.18916$ & $34.27$   & $-0.22173$ & $22.96$ &   $-0.22324$ & $22.43$ \\
				$16\!\times\!32$ &  $-0.28504$ & $0.96$ &     $-0.26034$ & $9.54$   &  $-0.28385$ & $1.37$ &    $-0.28403$ & $1.31$ \\
				$32\!\times\!64$ &  $-0.28670$ & $0.38$ &     $-0.27294$ & $5.16$   &  $-0.28781$ & $0.00$ &    $-0.28784$ & $0.01$ \\
				$64\!\times\!128$&  $-0.28745$ & $0.12$ &     $-0.27582$ & $4.16$   &  $-0.28784$ & $0.02$ &    $-0.28785$ & $0.02$ \\
				\hline
				\multicolumn{9}{c}{Reference solution: $u_Z(C)=-0.28780\,\text{mm}$ (\cite{MITC4Plus3})} \\ 
				\hline
			\end{tabularx} \nonumber
			\caption{Deflection $u_Z(C)$ for the partly clamped hyperbolic paraboloid shell with $h=0.001\,\text{m},~\rho=360\,\text{kg}\,\text{m}^{-3}$, distorted mesh}
		\end{table}
	\end{center}

	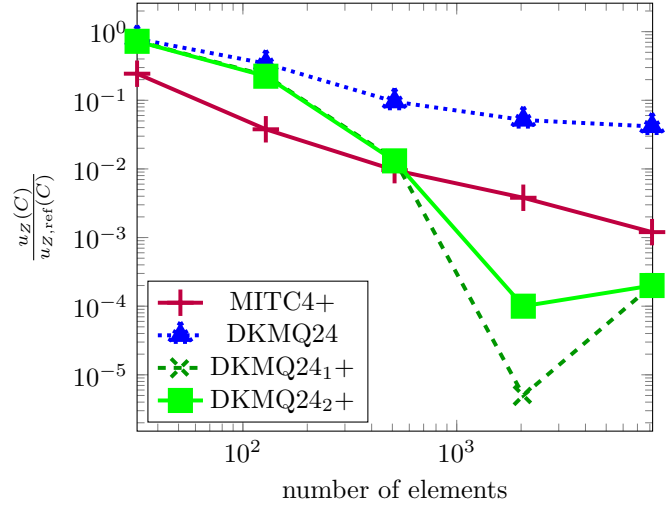
\begin{figure}
		\centering
		\begin{tikzpicture} 
		\begin{loglogaxis}[ 
		xlabel=\text{number of elements}, 
		ylabel= $\frac{u_Z(C)}{u_{Z,\rf}(C)}$,  
		x tick label style={
			/pgf/number format/1000 sep={}  
		},
		y tick label style={
			/pgf/number format/1000 sep={},
			/pgf/number format/fixed,       
			/pgf/number format/precision=5  
		},
		scaled y ticks=false, 
		xmin= 32,
		xmax = 8242,
		legend style={at={(0.02,0.02)},anchor=south west}
		] 
		\addplot[line width=1.5pt, color=purple, mark=+, mark size=5pt] coordinates {
			(32,   24.41E-2) 
			(128,  3.78E-2) 
			(512,  0.96E-2) 
			(2048, 0.38E-2)
			(8192, 0.12E-2)
		}; 
		\addplot[line width=1.5pt, color=blue, mark=triangle*, mark size=5pt, dotted] coordinates { 
			(32,   78.56E-2) 
			(128,  34.27E-2) 
			(512,  9.54E-2) 
			(2048, 5.16E-2)
			(8192, 4.16E-2)
		};
		\addplot[line width=1.5pt, color=col1, mark=x, mark size=5pt, dashed] coordinates { 
			(32,   73.05E-2) 
			(128,  22.96E-2) 
			(512,  1.37E-2) 
			(2048, 0.0005E-2)
			(8192, 0.02E-2)
		};
		\addplot[line width=1.5pt, color=green, mark=square*, mark size=4pt] coordinates { 
			(32,   72.35E-2) 
			(128,  22.43E-2) 
			(512,  1.31E-2) 
			(2048, 0.01E-2)
			(8192, 0.02E-2)
		};
		\legend{MITC4+, DKMQ24, DKMQ24$_1$+, DKMQ24$_2$+} 
		\end{loglogaxis} 
		\end{tikzpicture} 
		\caption{Deflection $u_Z(C)$ for the partly clamped hyperbolic paraboloid shell with $h=0.001\,\text{m},~\rho=360\,\text{kg}\,\text{m}^{-3}$, distorted mesh} 
		\label{fig:BENCH8c}
	\end{figure}
	
	\begin{center}
		\begin{table}
			\centering
			\begin{tabularx}{160mm}{ c *{8}{Y} }
				\hline
				Mesh  & \multicolumn{2}{c}{MITC4+}  & \multicolumn{2}{c}{DKMQ24} & \multicolumn{2}{c}{DKMQ24$_1$+} 
				& \multicolumn{2}{c}{DKMQ24$_2$+}  \rule[-2mm]{0mm}{6mm}  \\
				& $u_Z(C)$ & err. & $u_Z(C)$ & err. & $u_Z(C)$ & err. & $u_Z(C)$ & err. \\
				& [mm] &  [\%]    & [mm] &  [\%]    & [mm] &  [\%]  & [mm] &  [\%]      \\
				\hline
				\hline
				$4\!\times\!8$    & $-0.12648$ & $46.98$  &    $-0.00120$ & $99.50$   & $-0.00161$ & $99.33$ &   $-0.00173$ & $99.27$ \\
				$8\!\times\!16$   & $-0.21096$ & $11.57$  &    $-0.00751$ & $96.85$   & $-0.01770$ & $92.58$ &   $-0.01807$ & $92.42$ \\
				$16\!\times\!32$  & $-0.23703$ & $0.64$   &    $-0.03676$ & $84.59$   & $-0.10256$ & $57.01$ &   $-0.10314$ & $56.77$ \\
				$32\!\times\!64$  & $-0.23856$ & $0.00$   &    $-0.11378$ & $52.31$   & $-0.21383$ & $10.37$ &   $-0.21411$ & $10.25$ \\
				$64\!\times\!128$ & $-0.23844$ & $0.0005$ &    $-0.18645$ & $21.84$   &  $-0.23724$ & $0.55$ &    $-0.23726$ & $0.54$ \\
				\hline
				\multicolumn{9}{c}{Reference solution: $u_Z(C)=-0.23856\,\text{mm}$ (\cite{MITC4Plus3})} \\ 
				\hline
			\end{tabularx} \nonumber
			\caption{Deflection $u_Z(C)$ for the partly clamped hyperbolic paraboloid shell with $h=0.0001\,\text{m},~\rho=3.6\,\text{kg}\,\text{m}^{-3}$, distorted mesh}
		\end{table}
	\end{center}

	\begin{figure}
		\centering
		\begin{tikzpicture} 
		\begin{loglogaxis}[ 
		xlabel=\text{number of elements}, 
		ylabel= $\frac{u_Z(C)}{u_{Z,\rf}(C)}$,  
		x tick label style={
			/pgf/number format/1000 sep={}  
		},
		y tick label style={
			/pgf/number format/1000 sep={},
			/pgf/number format/fixed,       
			/pgf/number format/precision=5  
		},
		scaled y ticks=false, 
		xmin= 32,
		xmax = 8242,
		legend style={at={(0.02,0.02)},anchor=south west}
		] 
		\addplot[line width=1.5pt, color=purple, mark=+, mark size=5pt] coordinates {
			(32,   46.98E-2) 
			(128,  11.57E-2) 
			(512,  0.64E-2) 
			(2048, 0.01E-2)
			(8192, 0.05E-2)
		};  
		\addplot[line width=1.5pt, color=blue, mark=triangle*, mark size=5pt, dotted] coordinates { 
			(32,   99.50E-2) 
			(128,  96.85E-2) 
			(512,  84.59E-2) 
			(2048, 52.31E-2)
			(8192, 21.84E-2)
		};
		\addplot[line width=1.5pt, color=col1, mark=x, mark size=5pt, dashed] coordinates { 
			(32,   99.33E-2) 
			(128,  92.58E-2) 
			(512,  57.01E-2) 
			(2048, 10.37E-2)
			(8192,  0.55E-2)
		};
		\addplot[line width=1.5pt, color=green, mark=square*, mark size=4pt] coordinates { 
			(32,   99.27E-2) 
			(128,  92.42E-2) 
			(512,  56.77E-2) 
			(2048, 10.25E-2)
			(8192,  0.54E-2)
		};
		\legend{MITC4+, DKMQ24, DKMQ24$_1$+, DKMQ24$_2$+}
		\end{loglogaxis} 
		\end{tikzpicture} 
		\caption{Relative error of deflection $u_Z(C)$ for the partly clamped hyperbolic paraboloid shell with $h=0.0001\,\text{m},~\rho=3.6\,\text{kg}\,\text{m}^{-3}$, distorted mesh}
		\label{fig:BENCH8d}
	\end{figure}
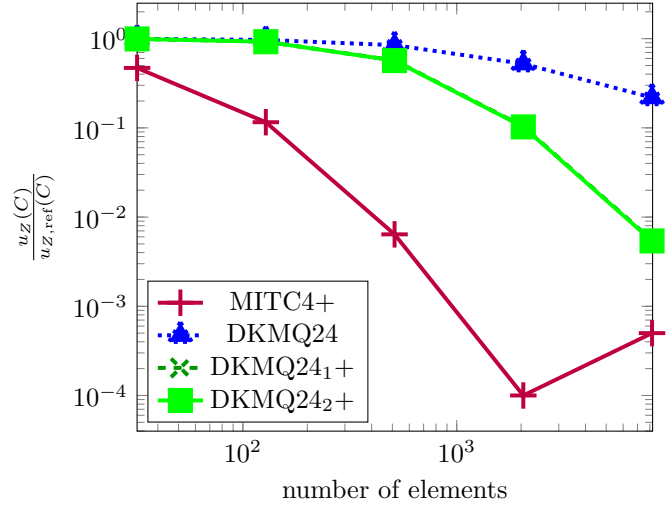

	\FloatBarrier

	\subsection{Square clamped plate under uniform loading with distorted mesh}
	\label{b9}
	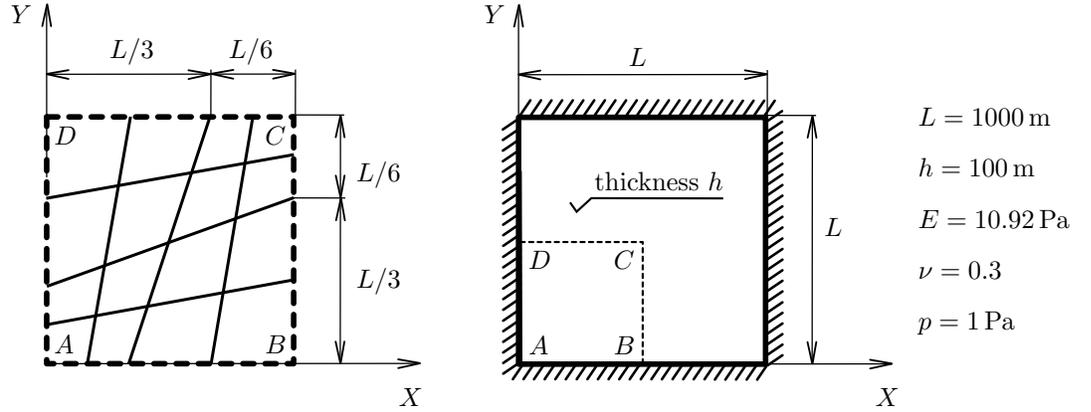
\begin{figure}[H]
		\centering
		\begin{overpic}[width=0.7\textwidth]{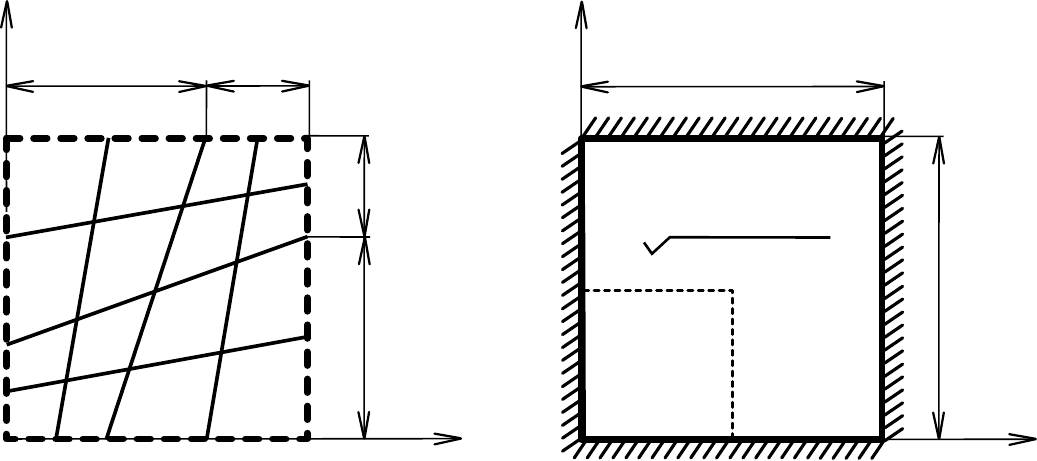} 
			\put(37,11){$L/3$}
			\put(37,23.5){$L/6$}
			\put(8,38){$L/3$}
			\put(22,38){$L/6$}
			\put(69,37){$L$}
			\put(92,17){$L$}
			\put(1.5,3.25){$A$}
			\put(26.3,3.25){$B$}
			\put(26.3,27.75){$C$}
			\put(1.5,27.75){$D$}
			\put(42,-3){$X$}
			\put(98,-3){$X$}
			\put(-3.6,42){$Y$}
			\put(52,42){$Y$}
			\put(57.25,3.15){$A$}
			\put(67.25,3.15){$B$}
			\put(67.25,13.15){$C$}
			\put(57.25,13.15){$D$}
			\put(65,22.4){$\text{thickness}~h$}		    
			\put(103,30){$L=1000\,\text{m}$}
			\put(103,24){$h=100\,\text{m}$}
			\put(103,18){$E=10.92\,\text{Pa}$}
			\put(103,12){$\nu=0.3$}	
			\put(103,6){$p=1\,\text{Pa}$}
		\end{overpic}
		\caption{Square clamped plate under uniform loading with distorted mesh $4\times4$.}
		\label{SquaredClampedPlate}
	\end{figure}
	\FloatBarrier
	
	A fully clamped pressure loaded square plate with distorted mesh is considered (see Figure (\ref{SquaredClampedPlate})). Only one quarter is computed due to symmetries. Boundary and symmetry conditions are as follows: $u_X=u_Y=u_Z=\varphi_X=\varphi_Y=\varphi_Z=0$ on sides AB and AD, $u_X=\varphi_Y=\varphi_Z=0$ on side BC and $u_Y=\varphi_X=\varphi_Z=0$ on side CD. This benchmark problem tests sensitivity to mesh distortion. The reference solution is taken from \cite{Srinivas}\footnote{The value $62.83$ given in \cite{Katili2015}, taken from \cite{Srinivas}, is not the deformation $u_Z$ as stated, but $\frac{Gu_Z}{hp}$.}. This benchmark is taken as an example of pure bending, where our modifications of the membrane part of the DKMQ24 element does not play any role.  Despite of that, our modified elements DKMQ24$_1$+ and DKMQ24$_2$+ provide significantly lower errors in all test cases, which is due to the nodal moment corrections.  
	
	\begin{center}
		\begin{table}
			\centering
			\begin{tabularx}{\tabulkywidth}{ c *{8}{Y} }
				\hline
				Mesh           &  \multicolumn{2}{c}{DKMQ24 \cite{Katili2015}}  & \multicolumn{2}{c}{DKMQ24} & \multicolumn{2}{c}{DKMQ24$_1$+}  & \multicolumn{2}{c}{DKMQ24$_2$+} \rule[-2mm]{0mm}{6mm}  \\
				&    $u_Z(C)$   & err. &    $u_Z(C)$   & err. &    $u_Z(C)$   & err. &    $u_Z(C)$   & err.\\
				& [m] &  [\%]    & [m] &  [\%]    & [m] &  [\%]  & [m] &  [\%]      \\
				\hline
				\hline
				$2\!\times2\!$      &  $-1920.3$    &  $28.37$ &    $-1910.8$   &  $27.73$   &  $-1725.3$    & $15.33$ &  $-1725.3$    & $15.33$    \\
				$4\!\times\!4$      &  $-1616.1$    &  $8.03$  &    $-1906.6$   &  $7.39$    &  $-1560.0$    &  $4.27$ &  $-1560.0$    &  $4.27$    \\
				$8\!\times\!8$      &  $-1538.0$    &  $2.81$  &    $-1529.2$   &  $2.22$    &  $-1517.6$    &  $1.45$ &  $-1517.6$    &  $1.45$    \\
				$16\!\times\!16$    &  $-1519.0$    &  $1.54$  &    $-1510.6$   &  $0.98$    &  $-1507.7$    &  $0.79$ &  $-1507.7$    &  $0.79$    \\
				$32\!\times\!32$    &   -           &    -     &    $-1506.1$   &  $0.68$    &  $-1505.4$    &  $0.63$ &  $-1505.4$    &  $0.63$    \\
				\hline
				\multicolumn{9}{c}{Reference solution: $u_Z(C)=-1495.95\,\text{m}$ (\cite{Srinivas})} \\ 
				\hline
			\end{tabularx} \nonumber
			\caption{Deflection $u_Z(C)$ for the plate with distorted mesh}
		\end{table}
	\end{center}

	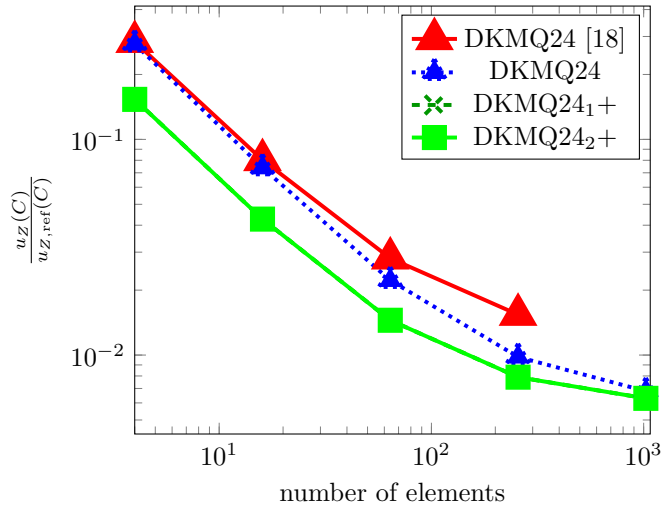
\begin{figure}
		\centering
		\begin{tikzpicture} 
		\begin{loglogaxis}[ 
		xlabel=\text{number of elements}, 
		ylabel= $\frac{u_Z(C)}{u_{Z,\rf}(C)}$,  
		x tick label style={
			/pgf/number format/1000 sep={}  
		},
		y tick label style={
			/pgf/number format/1000 sep={},
			/pgf/number format/fixed,       
			/pgf/number format/precision=5  
		},
		scaled y ticks=false, 
		xmin=4,
		xmax = 1074,
		legend style={at={(0.98,0.98)},anchor=north east}
		] 
		\addplot[line width=1.5pt, color=red, mark=triangle*, mark size=7pt] coordinates { 
			(4,  28.37E-2) 
			(16, 8.03E-2) 
			(64, 2.81E-2) 
			(256, 1.54E-2)
		};  
		\addplot[line width=1.5pt, color=blue, mark=triangle*, mark size=5pt, dotted] coordinates { 
			(4,  27.73E-2) 
			(16, 7.39E-2) 
			(64, 2.22E-2) 
			(256, 0.98E-2)
			(1024,0.68E-2)
		};
		\addplot[line width=1.5pt, color=col1, mark=x, mark size=5pt, dashed] coordinates { 
			(4,  15.33E-2)   
			(16, 4.27E-2) 
			(64, 1.45E-2) 
			(256, 0.79E-2) 
			(1024,0.63E-2)
		};
		\addplot[line width=1.5pt, color=green, mark=square*, mark size=4pt] coordinates {
			(4,  15.33E-2)   
			(16, 4.27E-2) 
			(64, 1.45E-2) 
			(256, 0.79E-2) 
			(1024,0.63E-2)
		};
		\legend{DKMQ24 \cite{Katili2015}, DKMQ24, DKMQ24$_1$+, DKMQ24$_2$+}
		\end{loglogaxis} 
		\end{tikzpicture} 
		\caption{Relative error of deflection $u_Z(C)$ for the square clamped plate under uniform loading with distorted mesh}
		\label{fig:BENCH9a}
	\end{figure}
	
	\FloatBarrier
	

	\section{Conclusions}
	The new quadrangle shell element DKMQ24$_2$+ with four nodes and six degrees of freedom per node, based on the DKMQ24 quadrangle shell element introduced by Katili et al. \cite{Katili2015}, with improved membrane behaviour was introduced. The element DKMQ24$_2$+ uses quadratic enrichment of the in-plane approximations with the help of drilling rotations, a selective reduced integration of shear terms, a proportional scaling of the penalty coefficient to the shell thickness $h$ and nodal moment corrections. Another variant of the element formulation called DKMQ24$_1$+ uses a constant penalty coefficient without any scaling and otherwise is identical with the DKMQ24$_2$+ variant. According to our numerical comparisons, the DKMQ24$_2$+ element variant is superior to the DKMQ24$_1$+ element variant.
	
	The new quadrangle shell element DKMQ24$_2$+, when compared to the DKMQ24 element by Katili \cite{Katili2015}, exhibits lower errors on majority of the tested benchmark problems, in particular on pure membrane problems and on shell problems where in-plane loading is dominant. Moreover, the presented element provides convergent behaviour at Raasch's hook benchmark.
	
	The DKMQ24$_2$+ shell element, when compared to the MITC4+ and \US elements, exhibits in many cases a lower level of error. The error is in general higher on coarser meshes and smaller on finer meshes. The resistance to distored meshes of the DKMQ24$_2$+ shell element is not as high as of the MITC4+ and \US shell elements.
	
	
	The DKMQ24$_2$+ shell element passes all our convergence tests, without presence of neither shear locking nor membrane locking and converges in all considered cases to the reference solution. Moreover, all our modifications, except of the static condensation, do not introduce any additional computational costs when compared to the DKMQ24 shell element. To conclude, the DKMQ24$_2$+ shell element could be considered as a replacement of the DKMQ24 shell element. 
	
	
	
	
	\section{Acknowledgements}
	
	The authors gratefully acknowledge the support of the FWF Austrian Science Fund through the START project Y 1093 as well as the support through the doctoral college Computational Design at TU Wien.
	\section{Appendix A -- Static condensation}
	We consider the following system of equations defined on a single element  $e$ 
	\begin{align}
	\left[
	\begin{array}{cc}
	\bm{K}_{\m\m}      & \bm{K}_{\m\n}  \\
	\bm{K}_{\m\n}^{\T} & \bm{K}_{\n\n}  \\
	\end{array}
	\right]\left[
	\begin{array}{c}
	\bm{q}_{\m} \\ \bm{q}_{\n} 
	\end{array}
	\right]&=
	\left[
	\begin{array}{c}
	\bm{f}_{\m} \\ \bm{0}
	\end{array}
	\right],
	\end{align}
	where
	\begin{align}
	\bm{K}_{\m\m}&=\int_{e} \bm{B}_{\m,3\times24}^{\T} \bm{D}_{\m} \bm{B}_{\m,3\times24}\,\text{d}A, \\
	\bm{K}_{\m\n}&=\int_{e} \bm{B}_{\m,3\times24}^{\T} \bm{D}_{\m} \bm{B}_{\n,3\times2}\,\text{d}A, \\
	\bm{K}_{\n\n}&=\int_{e} \bm{B}_{\n,3\times2}^{\T} \bm{D}_{\m} \bm{B}_{\n,3\times2}\,\text{d}A.
	\end{align}
	If we eliminate the unknown $\bm{q}_{\n}$, we get the eliminated system of equations of the form 
	\begin{align}
	\tilde{\bm{K}}_{\m\m}\bm{q}_{\m}&=\bm{f}_{\m},
	\end{align}
	where
	\begin{align}
	\tilde{\bm{K}}_{\m\m} = \bm{K}_{\m\m}-\bm{K}_{\m\n}\bm{K}_{\n\n}^{-1}\bm{K}_{\m\n}^{\T},
	\end{align}
	From the programming point of view the following equivalent form can be advantageous 
	\begin{align}
	\tilde{\bm{K}}_{\m\m}=\int_{e}\tilde{\bm{B}}^{\T}\bm{D}_{\m}\tilde{\bm{B}}\,\text{d}A,~~\tilde{\bm{B}}=\bm{B}_{\m}-\bm{B}_{\n}\bm{K}_{\n\n}^{-1}\bm{K}_{\m\n}^{\T},
	\end{align}
	because the condensation is implemented only on the strain matrix level.

\end{document}